\documentstyle{article}

\def\GDS{\mbox{\bf GDS}}
\def\DS{\mbox{\bf DS}}
\def\QPN{\mbox{\bf QPN}}
\def\PN{\mbox{\bf PN}}
\def\QPNN{$\QPN^{\neg}$}
\def\QDS{\mbox{\bf QDS}}
\def\GQDS{\mbox{\bf GQDS}}
\def\QDSN{$\QDS^{\neg P}$}
\def\QMDS{\mbox{\bf QMDS}}
\def\GQMDS{\mbox{\bf GQMDS}}
\def\QMPN{\mbox{\bf QMPN}}
\def\QMPNN{$\QMPN^\neg$}

\def\eL{\mbox{$\cal L$}}
\def\eLn{\mbox{$\cal L_{\neg}$}}
\def\Pe{\mbox{$\cal P$}}

\def\mX{{\mathbf{x}}}
\def\mY{{\mathbf{y}}}
\def\mU{{\mathbf{u}}}
\def\mV{{\mathbf{v}}}
\def\mZ{{\mathbf{z}}}

\def\mj{{\mathbf{1}}}

\def\kon{\wedge}
\def\cirk{\,{\raisebox{.3ex}{\tiny $\circ$}}\,}
\def\pl{\!+\!}
\def\mn{\!-\!}

\def\Dk{\hat{\Delta}}
\def\Dd{\check{\Delta}}
\def\Sk{\hat{\Sigma}}
\def\Sd{\check{\Sigma}}
\def\Dkp{\hat{\Delta}'}
\def\Ddp{\check{\Delta}'}
\def\Skp{\hat{\Sigma}'}
\def\Sdp{\check{\Sigma}'}

\def\str{\rightarrow}
\def\rts{\leftarrow}

\def\ks{\mbox{\footnotesize$\;\xi\;$}}
\def\prop#1#2{\vspace{2ex} \noindent{\sc #1.} {\it #2} \par \vspace{2ex}}
\def\dkz{\noindent{\sc Proof. }}
\def\qed{\hfill $\dashv$}

\def\HDS{\vrule width0pt height2.3ex
depth1.05ex\displaystyle}

\def\f#1#2{{{\HDS #1}\over{\HDS #2}}}

\def\b#1#2{\stackrel{\raisebox{-2pt}{\mbox{\tiny $#1$}}}
{\raisebox{0pt}{$b$}}^{\raisebox{-7pt}{\scriptsize $#2$}}}

\def\c#1{\stackrel{\raisebox{-2pt}{\mbox{\tiny $\,#1$}}}
{\raisebox{0pt}{$c$}}}

\def\th#1#2{\stackrel{\raisebox{-2pt}{\mbox{\tiny $\,#1$}}}
{\raisebox{0pt}{$\theta$}}^{\raisebox{-7pt}{\scriptsize $#2$}}}

\def\kst{\raisebox{1pt}{\mbox{\tiny$\xi$}}}

\def\mix{{\mbox{\it mix}}\,}
\def\koc{{\raisebox{-.2ex}{$\Box$}}}

\def\encircle#1{\begin{picture}(12,7)
\put(5,2){\circle{10}} \put(5.5,0){\makebox(0,0)[b]{$#1$}}
\end{picture}}

\def\HDSS{\vrule width0pt height3ex
depth1.75ex\displaystyle}

\def\afrac#1{{\phantom{\HDSS #1}\atop{\HDSS #1}}}

\def\fS#1#2{{{\HDSS #1}\over{\HDSS #2}}}

\begin{document}

\title{Coherence in Linear Predicate Logic}
\author{\small {\sc Kosta Do\v sen} and {\sc Zoran Petri\' c}
\\[1ex]
{\small Mathematical Institute, SANU}\\[-.5ex]
{\small Knez Mihailova 35, p.f.\ 367, 11001 Belgrade,
Serbia}\\[-.5ex]
{\small email: \{kosta, zpetric\}@mi.sanu.ac.yu}}
\date{}
\maketitle

\begin{abstract}
\noindent Coherence with respect to Kelly-Mac Lane graphs is
proved for categories that correspond to the multiplicative
fragment without constant propositions of classical linear
first-order predicate logic without or with mix. To obtain this
result, coherence is first established for categories that
correspond to the multiplicative conjunction-disjunction fragment
with first-order quantifiers of classical linear logic, a fragment
lacking negation. These results extend results of \cite{DP04} and
\cite{DP07}, where coherence was established for categories of the
corresponding fragments of propositional classical linear logic,
which are related to proof nets, and which could be described as
star-autonomous categories without unit objects.
\end{abstract}

\vspace{.3cm}


\vspace{.5ex}


\setcounter{section}{-1}

\section{Introduction} The goal of this paper is to prove
coherence for categories that correspond to the multiplicative
fragment without propositional constants (nullary connectives) of
classical linear first-order predicate logic without or with mix.
(In this fragment the modal operators ! and ? are left out.) The
propositional logic corresponding to this fragment is the fragment
of linear logic caught by proof nets. Coherence for categories
that correspond to this propositional logic, called
\emph{proof-net categories}, was proved in \cite{DP07}, where it
is also demonstrated that the notion of proof-net category is the
right notion of star-autonomous category without unit objects (and
where references to related work may be found; see \cite{BS04} for
a general categorial introduction to propositional linear logic).
The notion of proof-net category is here extended with assumptions
concerning first-order quantifiers, and this yields the notion of
category that corresponds to the fragment of linear predicate
logic mentioned in the first sentence. We prove coherence for
these categories. Coherence in \cite{DP07} and here is understood
as in Kelly's and Mac Lane's coherence result of \cite{KML71} for
symmetric monoidal closed categories. It is coherence with respect
to the same kind of graphs.

Such coherence results are very useful, because they enable us to
decide easily equality of canonical arrows. The coherence results
of this paper are interesting also for general proof theory. They
deal with a plausible notion of identity of proofs in logic (see
\cite{DP04}, Chapter 1).

The addition of first-order quantification to categories treated
in \cite{DP07} does not bring anything new with respect to the
graphs. Before they involved propositional letters, and now they
involve predicate letters. Individual variables do not bring
anything to these graphs. All the new arrows for quantifiers have
identity graphs (see Section 1.5). Although it seems we have made
an inessential addition (based on a trivial adjunction; see
Section 1.4), we don't know how to reduce simply the coherence
result proved in this paper to the coherence result of
\cite{DP07}. The proofs in this paper extend and modify those of
\cite{DP07}, but they require considerable additional effort.

We omit the multiplicative propositional constants from our
treatment because they present special problems for coherence
(though their addition too may be based on a trivial adjunction;
see \cite{DP04}, Section 7.9). These problems are comparable to
those that the unit object of symmetric monoidal closed categories
makes for coherence of these categories (see \cite{KML71}). We
believe that before attacking such problems one should first
settle more tractable matters.

It may seem that in the absence of the multiplicative
propositional constant $\top$ we will not be able express that a
formula is a theorem with sequents of the form ${A\vdash B}$, for
$A$ and $B$ formulae. These are the sequents of categorial proof
theory, and of this paper. For a theorem $B$ we should have a
derivable sequent ${\top\vdash B}$, but in the absence of $\top$
we shall have instead the sequents ${A\vdash A\kon B}$ derivable
for every formula $A$.

We omit also the additive (lattice) connectives from our
treatment. They would lead to the same kind of problem for
coherence that arises for classical or intuitionistic
conjunctive-disjunctive logic with quantifiers added. This is a
challenging matter, with which we intend to deal on another
occasion. The present paper should lay the ground for this future
work.

The categories corresponding to the fragment of linear predicate
logic that we cover are here presented equationally, in an
axiomatic, regular and surveyable manner. These axiomatic
equations should correspond to the combinatorial building blocks
of identity of proofs in this fragment of logic, as in knot theory
the Reidemeister moves are the combinatorial building blocks of
identity of knots and links (see \cite{BZ85}).

In this paper, our approach to categories corresponding to
first-order predicate logic is quite syntactical. We deal mainly
with freely generated categories, which are a categorial
presentation of syntax. Objects are formulae, and arrows are
proofs, or deductions, i.e.\ equivalence classes of derivations.
At the level of objects, our first-order language is quite
standard. After these freely generated categories are introduced,
other concrete categories belonging to the classes in which our
categories are free may be taken as models---of proofs, rather
than formulae. The only models of this kind that we consider in
this paper are categories whose arrows are graphs. Our coherence
results may be understood as completeness results with respect to
these models. What is shown complete is the axiomatization of
equality between arrows in the freely generated category. Previous
treatments of first-order quantifiers in categorial logic, which
started with the work of Lawvere (see Section 1.4, and references
given there), are less syntactical in spirit than ours.

As we said above, the proofs in this paper are based on proofs to
be found in \cite{DP07}, which are themselves based on proofs in
\cite{DP04}. We will eschew repeating this previously published
material with all its details, and so our paper will not be
self-contained. To make it self-contained would yield a rather
sizable book, overlapping excessively with \cite{DP07} and
\cite{DP04}. We suppose the reader is acquainted up to a point
with \cite{DP04} (at least Sections 3.2-3, 7.6-8 and 8.4) and with
\cite{DP07} (at least Chapters 2 and 6). Although to avoid
unnecessary lengthy repetitions we sometimes presuppose the reader
knows the previous material, and we make only remarks concerning
additions and changes, we have in general strived to make our text
as self-contained as possible, so the reader can get an idea of
what we do from this text only.

In the first part of the paper we deal with categories that
correspond to the multiplicative conjunction-disjunction fragment
with first-order quantifiers of classical linear logic. Coherence
proved for these categories extends results proved for the
corresponding propositional fragment in \cite{DP04} (Chapter 7).
When we add the multiplicative propositional constants to the
corresponding propositional categories we may obtain the linearly
distributive categories of \cite{CS97} (see also \cite{BS04},
Section~2), or categories with more equations presupposed, which
coherence requires (see \cite{DP04}, Section 7.9).

In this first part we introduce in detail the categorial notions
brought by quantifiers. The greatest novelty here may be the
treatment accorded to renaming of free individual variables, which
is not usually considered as a primitive rule of inference (see
Sections 1.2, 1.8 and 2.2). Taking this renaming as primitive
enables us to have categorial axioms that are regular, surveyable
and easy to handle. We prove in the first part a categorial
cut-elimination result, which says that every arrow is
\emph{equal} to a cut-free one. The proof of this result requires
a preparation involving change of individual variables. This
preparation gives a categorial form to ideas of Gentzen and
Kleene.

In the second part of the paper we add negation (the only
connective missing from the first part), and proceed to prove
coherence following the direction of \cite{DP07} (Chapter 2). In
the third part, we add the mix principle, and indicate what
adjustments should be made in the proofs of the previously
obtained results in order to obtain coherence also in the presence
of mix. The exposition in the first part, where we introduce new
matters concerning quantifiers, is in general more detailed than
the exposition in the second and third part, where we rely even
more heavily on previously published results, and where we suppose
that the reader has already acquired some dexterity. A brief
concluding section points to future work.

\section{Coherence of \QDS}
\subsection{The language \eL}
Let $\cal P$ be a set whose elements we call \emph{predicate
letters}, for which we use $P$, $R,\ldots$, sometimes with
indices. To every member of $\cal P$ we assign a natural number
$n\geq 0$, called its \emph{arity}. For every $n\geq 0$, we assume
that we have infinitely many predicate letters in $\cal P$ of
arity $n$.

To build the first-order language \eL\ generated by $\cal P$, we
assume that we have infinitely many \emph{individual variables},
which we call simply \emph{variables}, and for which we use $x$,
$y$, $z$, $u$, $v,\ldots$, sometimes with indices. Let $\mX_n$
stand for a sequence of variables of length $n\geq 0$. The
\emph{atomic formulae} of \eL\ are all of the form $P\mX_n$ for
$P$ a member of $\cal P$ of arity $n$. We assume throughout this
paper that $\!\ks\!\in\{\kon,\vee\}$ and
$Q\in\{\forall,\exists\}$. The symbols $\kon$ and $\vee$ are used
here for the multiplicative conjunction and disjunction
connectives (for which $\otimes$ and the inverted ampersand are
used in \cite{G87}).

The \emph{formulae} of \eL\ are defined inductively by the
following clauses:
\begin{itemize}
\item[] every atomic formula is a formula;\vspace{-1.5ex} \item[]
if $A$ and $B$ are formulae, then ${(A\ks B)}$ is a
formula;\vspace{-1.5ex}\item[] if $A$ is a formula and $x$ is a
variable, then $Q_xA$ is a formula.
\end{itemize}
As usual, we will omit the outermost parentheses of formulae and
take them for granted. We call $Q_x$ a \emph{quantifier prefix}.
(The advantage of the quantifier prefix $\forall_x$ over the more
usual $\forall x$ is that in ${\forall_xx=x}$ we do not need
parentheses, or a dot before $x=x$, for which a need is felt in
${\forall xx=x}$; in this paper, where we use the schematic letter
$Q$ for quantifiers, we want the quantifier prefix $Q_x$ to be
clearly distinguished from a formula $Px$.) For formulae we use
$A$, $B$, $C,\ldots$, sometimes with indices.

The notions of free and bound occurrences of variables in a
formula are understood as usual, and, as usual, we say that $x$
\emph{is free in} $A$ when there is at least one free occurrence
of $x$ in $A$. We say that $x$ \emph{is bound in} $A$ when $Q_x$
occurs in $A$ (though the quantifier prefix need not bind any
occurrence of $x$, as in $\forall_xPy$).

The variable $y$ is said to be free for substitution for $x$ in
$A$ when no free occurrence of $x$ in $A$ is in a subformula of
$A$ of the form $Q_yB$. We write $A^x_y$ for the result of
uniformly substituting $y$ for the free occurrences of $x$ in $A$,
provided that, as usual, $y$ is free for substitution for $x$ in
$A$ (this notation may be found in \cite{P65}). If this proviso is
not satisfied, then $A^x_y$ is not defined.

\subsection{The category \QDS}
The category \QDS, which we introduce in this section, corresponds
to the multiplicative conjunction-disjunction fragment with
first-order quantifiers of classical linear logic. This category
extends with quantifiers (this is where \textbf{Q} comes from) the
propositional category \DS\ of \cite{DP04} (Section 7.6).

The objects of the category \QDS\ are the formulae of \eL. To
define the arrows of \QDS, we define first inductively a set of
expressions called the \emph{arrow terms} of \QDS. Every arrow
term will have a \emph{type}, which is an ordered pair of formulae
of \eL. We write $f\!:A\vdash B$ when the arrow term $f$ is of
type $(A,B)$. Here $A$ is the \emph{source}, and $B$ the
\emph{target} of $f$. For arrow terms we use $f$, $g$, $h,\ldots$,
sometimes with indices. Intuitively, the arrow term $f$ is the
code of a derivation of the conclusion $B$ from the premise $A$
(which explains why we write $\vdash$ instead of $\str$).

For all formulae $A$, $B$ and $C$ of \eL, for every variable $x$,
and for all formulae $D$ of \eL\ in which $x$ is not free, the
following \emph{primitive arrow terms}:
\begin{tabbing}
\centerline{$\mj_A\!: A\vdash A$,}
\\[1ex]
\mbox{\hspace{.2em}}\= $\hat{b}^{\str}_{A,B,C}\,$\= : \=
$A\kon(B\kon C)\vdash (A\kon B)\kon C$,\quad \=
$\check{b}^{\str}_{A,B,C}\,$\= : \= $A\vee(B\vee C)\vdash (A\vee
B)\vee C$,\kill \> $\hat{b}^{\str}_{A,B,C}\,$\> : \> $A\kon(B\kon
C)\vdash (A\kon B)\kon C$,\> $\check{b}^{\str}_{A,B,C}\,$\> : \>
$A\vee(B\vee C)\vdash (A\vee B)\vee C$,
\\*[1ex]
\> $\hat{b}^{\rts}_{A,B,C}\,$\> : \> $(A\kon B)\kon C\vdash A\kon
(B\kon C)$,\> $\check{b}^{\rts}_{A,B,C}\,$\> : \> $(A\vee B)\vee
C\vdash A\vee (B\vee C)$,
\\[1ex]
\> \> $\hat{c}_{A,B}\!\!$\' : \> $A\kon B\vdash B\kon A$,\> \>
$\check{c}_{A,B}\!\!$\' : \> $B\vee A\vdash A\vee B$,
\\[1ex]
\centerline{$d_{A,B,C}\!: A\kon(B\vee C)\vdash (A\kon B)\vee C$,}
\\[2ex]
\> \> $\iota^{\forall\!_x}_A\!\!$\' : \> $\forall_xA\vdash A$,\>
\> $\iota^{\exists_x}_A\!\!$\' : \> $A\vdash\exists_xA$,
\\*[1ex]
\> \> $\gamma^{\forall\!_x}_D\!\!$\' : \> $D\vdash\forall_xD$,\>
\> $\gamma^{\exists_x}_D\!\!$\' : \> $\exists_xD\vdash D$,
\\[1ex]
\> \> $\check{\theta}^{\forall\!_x\str}_{A,D}\!\!$\' : \>
$\forall_x(A\vee D)\vdash \forall_xA\vee D$,\> \>
$\hat{\theta}^{\exists_x\rts}_{A,D}\!\!$\' : \> $\exists_xA\kon
D\vdash\exists_x(A\kon D)$
\end{tabbing}
are arrow terms. (Intuitively, these are the axioms of our logic
with the codes of their trivial derivations.)

Next we have the following inductive clauses:
\begin{itemize}
\item[]if ${f\!:A\vdash B}$ and ${g\!:B\vdash C}$ are arrow
terms,\\ then ${(g\cirk f)\!:A\vdash C}$ is an arrow
term;\vspace{-1ex} \item[]if ${f_1\!:A_1\vdash B_1}$ and
${f_2\!:A_2\vdash B_2}$ are arrow terms,\\ then ${(f_1\ks
f_2)\!:A_1\ks A_2\vdash B_1\ks B_2}$ is an arrow
term;\vspace{-1ex} \item[]if ${f\!:A\vdash B}$ is an arrow term,\\
then ${Q_xf\!:Q_xA\vdash Q_xB}$ and ${[f]^x_y\!:A^x_y\vdash
B^x_y}$ are arrow terms,
\end{itemize}
provided $A^x_y$ and $B^x_y$ are defined. (Intuitively, the
operations on arrow terms $\cirk$, $\!\ks\!$, $Q_x$ and
${[\;\;]^x_y}$ are codes of the rules of inference of our logic.)
This defines the arrow terms of \QDS. As we do with formulae, we
will omit the outermost parentheses of arrow terms.

The types of the arrow terms $\iota^{\forall\!_x}_A$ and
$\gamma^{\forall\!_x}_D$ are related to the logical principles of
universal \emph{instantiation} and universal \emph{generalization}
respectively (this is where $\iota$ and $\gamma$ come from). The
logical principle related to the type of $\iota^{\exists\!_x}_A$,
and not of $\gamma^{\exists\!_x}_D$, is sometimes called
\emph{existential generalization}, but for the sake of duality we
use $\iota$ and $\gamma$ with the existential quantifier as with
the universal quantifier.

The logical principles of the types of
$\check{\theta}^{\forall\!_x\str}_{A,D}$ and
$\hat{\theta}^{\exists_x\rts}_{A,D}$ are \emph{distributivity}
principles. The first, which is the intuitionistically spurious
\emph{constant domain} principle, is the converse of distribution
of disjunction over universal quantification, and the second is
distribution of conjunction over existential quantification. We
define below arrow terms with the converse types, which are both
intuitionistically valid (cf.\ also the end of the section).

With $Q_x$ and ${[\;\;]^x_y}$ we are given infinite families of
operations, indexed by variables. We call ${[\;\;]^x_y}$
\emph{renaming of free variables}, or for short just
\emph{renaming}. The operations $Q_x$ and $\!\ks\!$ are total, but
composition $\!\cirk\!$ and renaming are not total operations on
arrow terms. The result ${[f]^x_y}$ of applying renaming
${[\;\;]^x_y}$ to ${f\!:A\vdash B}$ is defined iff $A^x_y$ and
$B^x_y$ are defined.

Note that renaming is not substitution. The arrow terms
${[\mj_A]^x_y\!:A^x_y\vdash A^x_y}$ and
${\mj_{A^x_y}\!:A^x_y\vdash A^x_y}$ are different arrow terms. The
renaming operation is in the object language of arrow terms, while
the substitution operation $\,^x_y$ of $A^x_y$ is not in the
object language of formulae \eL, but only in the metalanguage.

Note that ${[g\cirk f]^x_y}$ may be defined though ${[f]^x_y}$ and
${[g]^x_y}$ are not defined (for example, with $f$ being
$\iota^{\exists_y}_{Rxy}\!:Rxy\vdash \exists_yRxy$ and $g$ being
${\iota^{\exists_x}_{\exists_yRxy}\!:\exists_yRxy\vdash
\exists_x\exists_yRxy}$, where ${(\exists_yRxy)^x_y}$ is not
defined). Note also that ${[f]^x_y}$ and ${[g]^x_y}$ may be
defined and composable without ${[g\cirk f]^x_y}$ being defined
(for example, with $f$ being
${\iota^{\forall_y}_{Py}\!:\forall_yPy\vdash Py}$ and $g$ being
${\iota^{\exists_x}_{Px}\!:Px\vdash \exists_xPx}$, where ${g\cirk
f}$ is not defined).

Renaming ${[\;\;]^x_y}$ is usually implicitly considered in proof
theory as a derivable rule when it is applied to ${f\!:A\vdash B}$
with $x$ not free either in $A$ or in $B$. For example, for $x$
not free in $D$, we have
\[
\f{\f{D\vdash B}{D\vdash \forall_xB}}{D\vdash B^x_y}
\]
It is also implicit in a derivable rule, which we have in the
presence of implication:
\[
\f{\f{\f{A\vdash B}{\vdash A\str B}}{\vdash\forall_x(A\str
B)}}{\f{\vdash A^x_y\str B^x_y}{A^x_y\vdash B^x_y}}
\]
We assume renaming here in the absence of implication. Renaming
corresponds to a structural rule of logic, in Gentzen's
terminology.

Next we define inductively the set of \emph{equations of} \QDS,
which are expressions of the form ${f=g}$, where $f$ and $g$ are
arrow terms of \QDS\ of the same type. These equations hold
whenever both sides are defined. For example, in the equation
$\mbox{($\forall\gamma$ {\it nat})}$ below we assume that both
sides are defined, which introduces the proviso that $x$ is free
neither in $A$ nor $B$. An analogous proviso is introduced already
by $({\mbox{{\it cat}~2}})$ below, where we assume that $f$ and
$g$, as well as $g$ and $h$, are composable. We will always assume
these provisos, but we will usually not mention them explicitly.
These tacit provisos are carried by the conventions of the
notation for arrow terms and conditions concerning substitution of
variables in formulae.

Intuitively, these equations should catch a plausible notion of
identity of proofs, proofs being understood as equivalence classes
of derivations. Coherence results should justify our calling this
notion of equality plausible. A justification may also be provided
by purely syntactical results, like cut elimination, and other
similar normal-form results. The two justifications may, but need
not, coincide (see \cite{DP04}, Chapter~1). In this paper, we
concentrate on coherence results for the justification, but as a
tool for demonstrating this coherence we establish cut-elimination
and normal-form results. The latter results also provide a partial
justification: they show the sufficiency of the equations assumed.
We do not consider here (like in \cite{D99}) the question whether
all these equations are also necessary for these or related
syntactical results. (Such a question should first be precisely
phrased.)

In the long list of axiomatic equations below, only the
quantificational equations and the renaming equations at the end
are new. The preceding propositional \DS\ equations are taken from
\cite{DP04} and \cite{DP07} (Section 2.1), while the first two
categorial equations are, of course, omnipresent. We stipulate
first that all the instances of ${f=f}$ and of the following
equations are equations of \QDS:

\begin{tabbing}
\hspace{1em}\emph{categorial equations}:\\*[1ex]
\mbox{\hspace{2em}}\= $({\mbox{{\it cat}~1}})$\quad\quad\= $f\cirk
\mj_A=\mj_B\cirk f=f\!:A\vdash B$,
\\*[1ex]
\> $({\mbox{{\it cat}~2}})$\> $h\cirk (g\cirk f)=(h\cirk g)\cirk
f$,
\\[3.5ex]
\hspace{1em}\emph{\DS\ equations}:\\*[1ex]
 \> $(\!\ks\, 1)$\> $\mj_A\ks\mj_B=\mj_{A\kst B}$,
\\*[1ex]
\> $(\!\ks\, 2)$\> $(g_1\cirk f_1)\ks(g_2\cirk f_2)=(g_1\ks
g_2)\cirk(f_1\ks f_2)$,
\\[2ex]
\hspace{1em}for $f\!:A\vdash D$, $g\!:B\vdash E$ and $h\!:C\vdash
F$,
\\*[1ex]
\> $\mbox{($\b{\xi}{\str}$ {\it nat})}$\>  $((f\ks g)\ks
h)\cirk\!\b{\xi}{\str}_{A,B,C}\:=\:
\b{\xi}{\str}_{D,E,F}\!\cirk(f\ks (g\ks h))$,
\\[1ex]
\> $\mbox{($\hat{c}$ {\it nat})}$\> $(g\kon f)\cirk
\hat{c}_{A,B}=\hat{c}_{D,E}\!\cirk(f\kon g)$,
\\*[1ex]
\> $\mbox{($\check{c}$ {\it nat})}$\> $(g\vee
f)\cirk\check{c}_{B,A}=\check{c}_{E,D}\!\cirk(f\vee g)$,
\\[1ex]
\> $\mbox{($d$ {\it nat})}$\> $((f\kon g)\vee h)\cirk d_{A,B,C} =
d_{D,E,F}\cirk(f\kon (g\vee h))$,
\\[1.5ex]
\pushtabs \hspace{2em}$(\b{\xi}{}\b{\xi}{})$\hspace{2em}\=
$\b{\xi}{\str}_{A,B,C}\!\cirk\!\b{\xi}{\rts}_{A,B,C}\;
=\mj_{(A\kst B)\kst C}$,\quad\quad\=
$\b{\xi}{\rts}_{A,B,C}\!\cirk\!\b{\xi}{\str}_{A,B,C}\;=\mj_{A\kst(B\kst
C)}$,
\\[1ex]
\hspace{2em}$(\b{\xi}{}\!5)$\> $\b{\xi}{\rts}_{A,B,C\kst
D}\!\cirk\! \b{\xi}{\rts}_{A\kst
B,C,D}\;=(\mj_A\:\ks\b{\xi}{\rts}_{B,C,D})\cirk\!\b{\xi}{\rts}_{A,B\kst
C,D}\!\cirk(\b{\xi}{\rts}_{A,B,C}\ks\:\mj_D)$,
\\[2ex]
\hspace{2em}$(\hat{c}\hat{c})$\> $\hat{c}_{B,A}\cirk
\hat{c}_{A,B}=\mj_{A\kon B}$,
\\*[1ex]
\hspace{2em}$(\check{c}\check{c})$\>$\check{c}_{A,B}\cirk
\check{c}_{B,A}=\mj_{A\vee B}$,
\\[1.5ex]
\hspace{2em}$(\hat{b}\hat{c})$\> $(\mj_B\:\kon
\hat{c}_{C,A})\cirk\! \hat{b}^{\rts}_{B,C,A}\cirk \hat{c}_{A,B\kon
C}\cirk
\hat{b}^{\rts}_{A,B,C}\cirk(\hat{c}_{B,A}\kon\:\mj_C)\,$\=$=
\hat{b}^{\rts}_{B,A,C}$,
\\*[1ex]
\hspace{2em}$(\check{b}\check{c})$\> $(\mj_B\:\vee
\check{c}_{A,C})\cirk \check{b}^{\rts}_{B,C,A}\cirk
\check{c}_{B\vee C,A} \cirk
\check{b}^{\rts}_{A,B,C}\cirk(\check{c}_{A,B}\vee \mj_C)$\>$=
\check{b}^{\rts}_{B,A,C}$,
\\[1.5ex]
\hspace{2em}$(d \kon)$\> $(\hat{b}^{\rts}_{A,B,C}\vee \mj_D)\cirk
d_{A\kon B,C,D}\,$\=$=d_{A,B\kon C,D}\cirk(\mj_A\kon
d_{B,C,D})\cirk\hat{b}^{\rts}_{A,B,C\vee D}$,
\\[1ex]
\hspace{2em}$(d \vee)$\> $d_{D,C,B\vee A}\cirk(\mj_D \kon
\check{b}^{\rts}_{C,B,A})$\>$=\check{b}^{\rts}_{D\kon
C,B,A}\cirk(d_{D,C,B}\vee\mj_A)\cirk d_{D,C\vee B,A}$,
\\[1.5ex]
\hspace{1em}for $d^R_{C,B,A}=_{df}\check{c}_{C,B\kon A} \cirk
(\hat{c}_{A,B} \vee\:\mj_C)\cirk d_{A,B,C}\cirk(\mj_A\kon
\check{c}_{B,C})\cirk \hat{c}_{C\vee B,A}:$
\\*[.5ex]
\` $(C\vee B)\kon A\vdash C\vee(B\kon A)$,
\\[1ex]
\hspace{2em}$(d\hat{b})$\> $d^R_{A\kon
B,C,D}\cirk(d_{A,B,C}\kon\mj_D)\,$\=$=d_{A,B,C\kon
D}\cirk(\mj_A\kon d^R_{B,C,D})\cirk \hat{b}^{\rts}_{A,B\vee C,D}$,
\\*[1ex]
\hspace{2em}$(d\check{b})$\>$(\mj_D\vee d_{C,B,A})\cirk
d^R_{D,C,B\vee A}$\>$=\check{b}^{\rts}_{D,C\kon
B,A}\!\cirk(d^R_{D,C,B}\vee\mj_A)\cirk d_{D\vee C,B,A}$, \poptabs
\\[4.5ex]
\hspace{1em}\emph{quantificational equations}:\\*[1ex]
\hspace{2em}\= $(Q1)$\hspace{2.5em}\= $Q_x\mj_A=\mj_{Q\!_xA}$,
\\*[1ex]
\> $(Q2)$\> $Q_x(g\cirk f)=Q_xg\cirk Q_xf$,
\\[1.5ex]
\hspace{1em}for $f\!:A\vdash B$,\\*[1ex] \> $\mbox{($\forall\iota$
{\it nat})}$\quad\> $f\cirk
\iota^{\forall\!_x}_A=\iota^{\forall\!_x}_B\cirk\forall_xf$,\hspace{5em}\=
$\mbox{($\exists\iota$ {\it nat})}$\quad\=
$\exists_xf\cirk\iota^{\exists_x}_A=\iota^{\exists_x}_B\cirk f$,
\\*[1ex]
\> $\mbox{($\forall\gamma$ {\it nat})}$\>
$\forall_xf\cirk\gamma^{\forall\!_x}_A=\gamma^{\forall\!_x}_B\cirk
f$,\> $\mbox{($\exists\gamma$ {\it nat})}$\>
$f\cirk\gamma^{\exists_x}_A=\gamma^{\exists_x}_B\cirk\exists_xf$,
\\[1.5ex]
\> $(\forall\beta)$\>
$\iota^{\forall\!_x}_A\cirk\gamma^{\forall\!_x}_A=\mj_A$,\>
$(\exists\beta)$\>
$\gamma^{\exists_x}_A\cirk\iota^{\exists_x}_A=\mj_A$,
\\*[1ex]
\> $(\forall\eta)$\>
$\forall_x\iota^{\forall\!_x}_A\cirk\gamma^{\forall\!_x}_{\forall\!_xA}=\mj_{\forall\!_xA}$,\>
$(\exists\eta)$\>
$\gamma^{\exists_x}_{\exists_xA}\cirk\exists_x\iota^{\exists_x}_A=\mj_{\exists_xA}$,
\\[2ex]
\hspace{1em}for \= $\check{\theta}^{\forall\!_x\rts}_{A,D}\:$\=
$=_{df}\forall_x(\iota^{\forall\!_x}_A\vee\mj_D)\cirk\gamma^{\forall\!_x}_{\forall\!_xA\vee
D}\!:\forall_xA\vee D\vdash\forall_x(A\vee D)$,
\\*[1ex]
\> $\hat{\theta}^{\exists_x\str}_{A,D}$\>
$=_{df}\gamma^{\exists_x}_{\exists_xA\kon
D}\cirk\exists_x(\iota^{\exists_x}_A\kon\mj_D)\!:\exists_x(A\kon
D)\vdash\exists_xA\kon D$ and
\\[1ex]
\> $({\!\ks\!},Q)\in\{(\vee,\forall),(\kon,\exists)\}$,
\\[1.5ex]
\hspace{2em}\= $(Q1)$\hspace{2.5em}\=
$Q_x\mj_A=\mj_{Q\!_xA}$,\kill \>
${(Q\th{\xi}{}\th{\xi}{})}$\>$\th{\xi}{Q\!_x\rts}_{A,D}\!\!\cirk\!\th{\xi}{Q\!_x\str}_{A,D}=\mj_{Q\!_x(A\kst
D)}$,\quad\quad$\th{\xi}{Q\!_x\str}_{A,D}\!\!\cirk\!\th{\xi}{Q\!_x\rts}_{A,D}=\mj_{Q\!_xA\kst
D}$,
\\[3.5ex]
\hspace{1em}\emph{renaming equations}:\\*[1ex] \hspace{1em}for
$x$, $y$, $z$ and $v$ mutually different variables and\\*
\hspace{1em}$\alpha_{A_1,\ldots,A_n}$ a primitive arrow term
except $\iota^{Q\!_x}_A$,
\\*[1.5ex]
\hspace{2em}\= $\mbox{({\it ren} $\alpha$)}$\quad\=
$[\alpha_{A_1,\ldots,A_n}]^x_y=\alpha_{(A_1)^x_y,\ldots,(A_n)^x_y}$,
\\[2ex]
\> $\mbox{({\it ren} $\cirk\!$)}$\> $[g\cirk
f]^x_y=[g]^x_y\cirk[f]^x_y$,
\\[1.5ex]
\> $\mbox{({\it ren} $\ks\!$)}$\> $[f_1\ks
f_2]^x_y=[f_1]^x_y\ks[f_2]^x_y$,
\\[1.5ex]
\> $\mbox{({\it ren} $Q$)}$\> $[Q_zf]^x_y=Q_z[f]^x_y$,
\\[2.5ex]
\> $\mbox{({\it ren} 1)}$\> $[f]^x_x=f$,
\\*[1.5ex]
\> $\mbox{({\it ren} 2)}$\> $[f]^x_y=f$, \hspace{.5em}if $x$ is
free neither in the source nor in the target of $f$,
\\[1.5ex]
\> $\mbox{({\it ren} 3)}$\> $[[f]^z_v]^x_y=[[f]^x_y]^z_v$,
\\[1.5ex]
\> $\mbox{({\it ren} 4)}$\> $[[f]^z_y]^x_y=[[f]^x_y]^z_y$,
\\[1.5ex]
\> $\mbox{({\it ren} 5)}$\> $[[f]^z_x]^x_y=[[f]^z_y]^x_y$,
\\[1.5ex]
\> $\mbox{({\it ren} 6)}$\> $[[f]^y_x]^x_y=[f]^x_y$.
\end{tabbing}

This concludes the list of axiomatic equations stipulated for
\QDS. To define all the equations of \QDS\ it remains only to say
that the set of these equations is closed under symmetry and
transitivity of equality and under the rules
\[
(\!\cirk\;\mbox{\it cong})\quad \f{f=f'\quad \quad \quad g=g'}
{g\cirk f=g'\cirk f'}\hspace{5em}(\!\ks\;\mbox{\it cong})\quad
\f{f_1=f_1'\quad \quad \quad f_2=f_2'} {f_1\ks f_2=f_1'\ks f_2'}
\]
\[
(Q\;\mbox{\it cong})\quad \f{f=f'}
{Q_xf=Q_xf'}\quad\quad\quad\quad\quad\quad
\]
\[
(\mbox{\it ren cong})\quad \f{f=f'}
{[f]^x_y=[f']^x_y}\quad\quad\quad\quad\quad\quad
\]

On the arrow terms of \QDS\ we impose the equations of \QDS. This
means that an arrow of \QDS\ is an equivalence class of arrow
terms of \QDS\ defined with respect to the smallest equivalence
relation such that the equations of \QDS\ are satisfied (see
\cite{DP04}, Section 2.3, for details).

The equations $(\!\ks\, 1)$, $(\!\ks\, 2)$, $(Q1)$ and $(Q2)$ are
called \emph{functorial} equations. They say that $\!\ks\!$ is a
biendofunctor and $Q_x$ an endofunctor of \QDS\ (i.e.\
2-endo\-functor and 1-endofunctor respectively, in the terminology
of \cite{DP04}, Section 2.4). The equations with ``\emph{nat}'' in
their names are called \emph{naturality} equations. The naturality
equations above say that $\b{\xi}{\str}$, $\c{\xi}$, $d$ and
$\iota^{Q_x}$ are natural transformations ($\b{\xi}{\rts}$ is a
natural transformation too, due to $(\b{\xi}{}\b{\xi}{})$). The
naturality equation $\mbox{($Q\gamma$ {\it nat})}$ says that
$\gamma^{Q_x}$ has some properties of a natural transformation,
but one side of $\mbox{($Q\gamma$ {\it nat})}$ is not always
defined when the other is. We will see later (in Section 1.4)
where $\gamma^{Q\!_x}$ gives rise to a natural transformation. As
$\gamma^{Q\!_x}$, so $\check{\theta}^{\forall\!_x\str}$ and
$\hat{\theta}^{\exists_x\rts}$ have some properties of natural
transformations due to the equations ${(Q\th{\xi}{}\th{\xi}{})}$.
The equations $(\b{\xi}{}\b{\xi}{})$, $(\c{\xi}\c{\xi})$ and
${(Q\th{\xi}{}\th{\xi}{})}$ are equations of isomorphisms.

In spite of the equations \mbox{({\it ren} $\mj$)} and \mbox{({\it
ren} $\cirk\!$)}, which are also functorial equations, renaming
combined with substitution applied to formulae does not give an
endofunctor of \QDS, because we do not have totally defined
functions. If the renaming operations were not assumed as
primitive operations on arrow terms for defining \QDS\ (we will
see in Sections 1.8 and 2.2 that $[\iota^{Q\!_x}_A]^x_y$ could be
assumed instead), then we would have problems in formulating
assumptions that give the following equation of \QDS:
\[
[f]^x_y\cirk[\iota^{\forall\!_x}_A]^x_y=[\iota^{\forall\!_x}_B]^x_y\cirk\forall_xf,
\]
which follows from \mbox{($\forall\iota$ {\it nat})}, (\mbox{\it
ren cong}) and $\mbox{({\it ren} $\cirk\!$)}$; we would not know
what to write for $[f]^x_y$.

The following equation:
\begin{tabbing}
\hspace{7em}$(\forall\gamma\iota)$\hspace{2em}\=
$\gamma^{\forall\!_x}_A\cirk\iota^{\forall\!_x}_A=\mj_{\forall\!_xA}$
\end{tabbing}
holds in \QDS. (By our tacit presupposition, introduced before
presenting the equations of \QDS, the variable $x$ is here not
free in $A$.) This equation is derived as follows:
\begin{tabbing}
\hspace{11em}\=
$\gamma^{\forall\!_x}_A\cirk\iota^{\forall\!_x}_A\,$\=
$=\forall_x\iota^{\forall\!_x}_A\cirk\gamma^{\forall\!_x}_{\forall\!_xA}$,
\hspace{.5em}\=by $\mbox{($\forall\gamma$ {\it nat})}$,
\\*[1ex]
\> \> $=\mj_{\forall\!_xA}$, \>by $(\forall\eta)$.
\end{tabbing}
In an analogous manner, we derive in \QDS\ the equation
\begin{tabbing}
\hspace{7em}\=$(\forall\gamma\iota)$\hspace{2em}\=
$\gamma^{\forall\!_x}_A\cirk\iota^{\forall\!_x}_A=\mj_{\forall\!_xA}$\kill

\>$(\exists\gamma\iota)$\>$
\iota^{\exists_x}_A\cirk\gamma^{\exists_x}_A=\mj_{\exists_xA}$.
\end{tabbing}

With the help of the equations $(Q\gamma\iota)$ we derive easily
the following equations analogous to $(Q1)$:
\begin{tabbing}
\hspace{7em}\=$(\forall\gamma\iota)$\hspace{2em}\=
$\gamma^{\forall\!_x}_A\cirk\iota^{\forall\!_x}_A=\mj_{\forall\!_xA}$\kill

\>$(Q\iota)$\>$Q_x\iota^{Q\!_x}_A\;$\=$=\iota^{Q\!_x}_{Q\!_xA}$,
\\*[1ex]
\>$(Q\gamma)$\>$Q_x\gamma^{Q\!_x}_A$\>$=\gamma^{Q\!_x}_{Q\!_xA}$.
\end{tabbing}
Note that $(Q\iota)$ can replace $(Q\eta)$ for axiomatizing the
equations of \QDS, but $(Q\gamma)$ cannot do so, because for it we
presuppose that $x$ is not free in $A$, which we do not presuppose
for $(Q\eta)$.

Note that $(\forall\eta)$ and $(\exists\eta)$ could be replaced
respectively by the equations
\begin{tabbing}
\hspace{7em}\=$(\forall\gamma\iota)$\hspace{2em}\=
$\gamma^{\forall\!_x}_A\cirk\iota^{\forall\!_x}_A=\mj_{\forall\!_xA}$\kill

\>\mbox{($\forall$ {\it ext})}\>$
\forall_x(\iota^{\forall\!_x}_A\cirk
f)\cirk\gamma^{\forall\!_x}_B\;\:$\=$=f\!:B\vdash\forall_xA$,
\\*[1ex]
\>\mbox{($\exists$ {\it ext})}\>$
\gamma^{\exists_x}_B\cirk\exists_x(g\cirk\iota^{\exists_x}_A)$\>$=g\!:\exists_xA\vdash
B$,
\end{tabbing}
which are easily derived with $\mbox{($Q\gamma$ {\it nat})}$ and
$(Q\eta)$. (In both of these equations we tacitly presuppose that
$x$ is not free in $B$.)

By relying on $(Q\beta)$ and $(Q\gamma\iota)$ we can easily derive
the following equations if $x$ is free neither in $A$ nor in $D$:
\begin{tabbing}
\hspace{7em}\=$(\forall\gamma\iota)$\hspace{2em}\=
$\gamma^{\forall\!_x}_A\cirk\iota^{\forall\!_x}_A=\mj_{\forall\!_xA}$\kill

\>\>$
\check{\theta}^{\forall\!_x\str}_{A,D}\;$\=$=(\gamma^{\forall\!_x}_A\vee\mj_D)\cirk\iota^{\forall\!_x}_{A\vee
D}$,\\*[1.5ex] \>\>$
\hat{\theta}^{\exists_x\rts}_{A,D}$\>$=\iota^{\exists_x}_{A\kon
D}\cirk(\gamma^{\exists_x}_A\kon\mj_D)$.
\end{tabbing}

Note that if $x$ is free in $A$ and not free in $D$, then in \QDS\
we do not have arrows of the types converse to the types of the
following distributivity arrows:
\begin{tabbing}
\hspace{5em}\=$\hat{\theta}^{\forall\!_x\rts}_{A,D}\,$\=$=_{df}\forall_x(\iota^{\forall\!_x}_A\kon\mj_D)\cirk
\gamma^{\forall\!_x}_{\forall\!_xA\kon D}\,$\=$:\forall_xA\kon
D\vdash\forall_x(A\kon D)$,
\\*[1ex]
\>$\check{\theta}^{\exists_x\str}_{A,D}$\>$=_{df}\gamma^{\exists_x}_{\exists_xA\vee
D}\cirk\exists_x(\iota^{\exists_x}_A\vee\mj_D\!)$\>$:\exists_x(A\vee
D)\vdash\exists_xA\vee D$,
\end{tabbing}
which are analogous to the arrows
$\check{\theta}^{\forall\!_x\rts}_{A,D}$ and
$\hat{\theta}^{\exists_x\str}_{A,D}$ respectively. (That these
arrows do not exist in \QDS\ is shown via cut elimination in
\GQDS; see Sections 1.5-9.) So we cannot have a prenex normal form
for formulae, i.e.\ objects.

\subsection{Change of bound variables}
We call \emph{change} of bound variables what could as well be
called \emph{renaming} of bound variables, because we do not want
to confuse this renaming with renaming of free variables. We
define in \QDS\ the following arrows, which formalize change of
bound variables ($\tau$ might come from ``transcribe''):
\begin{tabbing}
\mbox{\hspace{1em}}\= $(Q\tau)$\hspace{3.5em}\=
$\tau^{Q\!_x}_{A,u,v}=\tau^{Q\!_y}_{A^x_z,u,v}$,\kill

\>\>$\tau^{\forall\!_x}_{A,u,v}=_{df}\forall_v[\iota^{\forall\!_u}_{A^x_u}]^u_v
\cirk\gamma^{\forall\!_v}_{\forall\!_uA^x_u}\!:\forall_uA^x_u\vdash
\forall_vA^x_v$,
\\*[1.5ex]
\>\>$\tau^{\exists_x}_{A,v,u}=_{df}\gamma^{\exists_v}_{\exists_uA^x_u}\cirk
\exists_v[\iota^{\exists_u}_{A^x_u}]^u_v\!:\exists_vA^x_v
\vdash\exists_uA^x_u$,
\end{tabbing}
provided $u$ and $v$ are not free in $A$.

Note that $\tau^{Q\!_x}_{A,u,v}$ is the same arrow term as
$\tau^{Q\!_y}_{A^x_y,u,v}$ for $y$ not free in $A$, and a fortiori
for $y$ neither free nor bound in $A$. The variable $x$ in
$\tau^{Q\!_x}_{A,u,v}$ is just a place holder, which can always be
replaced by an arbitrary new variable.

We can derive the following equations of \QDS:
\begin{tabbing}
\mbox{\hspace{1em}}\= $(Q\tau)$\hspace{3.5em}\=
$\tau^{Q\!_x}_{A,u,v}=\tau^{Q\!_y}_{A^x_y,u,v}$, if $y$ is not
free in $A$,\kill

\> $({\mbox{$Q\tau$ {\it ren}}})$\>
$[\tau^{Q\!_x}_{A,u,v}]^y_z=\tau^{Q\!_x}_{A^y_z,u,v}$, if $x$ is
not $y$ or $z$
\end{tabbing}
(if $y$ is free in $A$ and $z$ is $u$ or $v$, then the right-hand
side of $({\mbox{$Q\tau$ {\it ren}}})$ is undefined),
\begin{tabbing}
\mbox{\hspace{1em}}\= $(Q\tau)$\hspace{3.5em}\= \kill \>
$({\mbox{$Q\tau$ {\it nat}}})$\> $Q_v[f]^x_v\cirk
\tau^{Q\!_x}_{A,u,v}=\tau^{Q\!_x}_{B,u,v}\cirk Q_u[f]^x_u$,
\\[1.5ex]
\> $({\mbox{$Q\tau$ {\it ref}}})$\>
$\tau^{Q\!_x}_{A,u,u}=\mj_{Q\!_uA^x_u}$,
\\[1.5ex]
\> $({\mbox{$Q\tau$ {\it sym}}})$\> $\tau^{Q\!_x}_{A,v,u}\cirk
\tau^{Q\!_x}_{A,u,v}\,$\=$=\mj_{Q\!_uA^x_u}$,
\\[1.5ex]
\> $({\mbox{$Q\tau$ {\it trans}}})$\> $\tau^{Q\!_x}_{A,v,w}\cirk
\tau^{Q\!_x}_{A,u,v}$\>$=\tau^{Q\!_x}_{A,u,w}$.
\end{tabbing}
From the equation $({\mbox{$Q\tau$ {\it sym}}})$ we see that
$\tau^{Q\!_x}_{A,u,v}$ and $\tau^{Q\!_x}_{A,v,u}$ are inverse to
each other. We can also derive the following equations of \QDS:
\begin{tabbing}
\mbox{\hspace{1em}}\=$(Q\tau)$\hspace{3.5em}\= \kill

\> $(\forall\tau\iota)$\> $\iota^{\forall\!_v}_{A^x_v}\cirk
\tau^{\forall\!_x}_{A,u,v}\,$\=$=[\iota^{\forall\!_u}_{A^x_u}]^u_v$,
\hspace{4em}\=$(\exists\tau\iota)$\hspace{3em}\=
$\tau^{\exists_x}_{A,v,u}\cirk\iota^{\exists_v}_{A^x_v}\,$\=$=
[\iota^{\exists_u}_{A^x_u}]^u_v$,
\\*[1.5ex]
\>$(\forall\tau\gamma)$\>
$\tau^{\forall\!_x}_{A,u,v}\cirk\gamma^{\forall\!_u}_A$\>$=\gamma^{\forall\!_v}_A$,
\>$(\exists\tau\gamma)$\>$\gamma^{\exists_u}_A\cirk
\tau^{\exists_x}_{A,v,u}$\>$=\gamma^{\exists_v}_A$,
\\[2ex]
\>$(\forall\tau\check{\theta})$\>
$(\tau^{\forall\!_x}_{A,u,v}\vee\mj_D)\cirk\check{\theta}^{\forall\!_u\str}_{A^x_u,D}=
\check{\theta}^{\forall\!_v\str}_{A^x_v,D}\cirk
\tau^{\forall\!_x}_{A\vee D,u,v}$,
\\*[1.5ex]
\>$(\exists\tau\hat{\theta})$\>$\tau^{\exists_x}_{A\kon
D,u,v}\cirk\hat{\theta}^{\exists_u\rts}_{A^x_u,D}=
\hat{\theta}^{\exists_v\rts}_{A^x_v,D}\cirk(\tau^{\exists_x}_{A,u,v}\kon\mj_D)$.
\end{tabbing}
To derive $(\forall\tau\check{\theta})$ we derive
\[
\check{\theta}^{\forall\!_v\rts}_{A^x_v,D}\cirk(\tau^{\forall\!_x}_{A,u,v}\vee\mj_D)=\tau^{\forall\!_x}_{A\vee
D,u,v}\cirk\check{\theta}^{\forall\!_u\rts}_{A^x_u,D}
\]
with the help of $(\forall\tau\iota)$ and $(\forall\tau\gamma)$.
We proceed analogously for $(\exists\tau\hat{\theta})$.

Note that $\tau^{\forall\!_x}$ and
$\check{\theta}^{\forall\!_x\rts}$, as well as $\tau^{\exists_x}$
and $\hat{\theta}^{\exists_x\str}$, have analogous definitions.
The following two equations of \QDS\ are analogous to the
equations $(Q\tau\iota)$:
\begin{tabbing}
\hspace{7em}\=$(\forall\theta\iota)$\hspace{2.5em}\=
$(\iota^{\forall\!_x}_A\vee\mj_D)\cirk\check{\theta}^{\forall\!_x\str}_{A,D}\:$\=$=
\iota^{\forall\!_x}_{A\vee D}$\kill

\>\>$\iota^{\forall\!_x}_{A\vee
D}\cirk\check{\theta}^{\forall\!_x\rts}_{A,D}\;
$\=$=\iota^{\forall\!_x}_A\vee\mj_D$,
\\[1ex]
\>\>$\hat{\theta}^{\exists_x\str}_{A,D}\cirk\iota^{\exists_x}_{A\kon
D}$\>$=\iota^{\exists_x}_A\kon\mj_D$.
\end{tabbing}
As a consequence of these two equations we have
\begin{tabbing}
\hspace{7em}\=$(\forall\check{\theta}\iota)$\hspace{2.5em}\=
$(\iota^{\forall\!_x}_A\vee\mj_D)\cirk\check{\theta}^{\forall\!_x\str}_{A,D}\:$\=$=
\iota^{\forall\!_x}_{A\vee D}$,
\\[1ex]
\>$(\exists\hat{\theta}\iota)$
\>$\hat{\theta}^{\exists_x\rts}_{A,D}\cirk(\iota^{\exists_x}_A\!\kon\mj_D)$\>$=\iota^{\exists_x}_{A\kon
D}$.
\end{tabbing}

\subsection{Quantifiers and adjunction}
Lawvere's presentation of predicate logic in categorial terms (see
\cite{Law69}, \cite{Law70} and \cite{LawR03}, Appendix A.1), and
presentations that follow him more or less closely (see, for
instance, \cite{S83}, \cite{Cub97}, \cite{P96} and \cite{J99},
Chapter 4), are less syntactical than ours. They do not pay close
attention to syntax. If this syntax were to be supplied precisely,
then a language without variables, like Quine's variable-free
language for predicate logic (see \cite{Q81}, and references
therein), called predicate functor logic, would be more
appropriate. Our first-order language is on the contrary quite
standard. It should be mentioned also that Lawvere's approach is
more general, whereas we concentrate on first-order logic.

Lawvere characterized quantifiers in intuitionistic logic through
an adjoint situation. In Lawvere's characterization of
quantifiers, functors from which the universal and existential
quantifiers arise are respectively the right and left adjoints of
a functor that is an instance, involving product types and
projections, of a functor Lawvere calls \emph{substitution}. An
approach in this style to linear predicate logic was first made in
\cite{SEELY89} (Section 2.5, Remark~3).

We will now present two kinds of adjoint situations that involve
the quantifiers of \QDS. These adjunctions are related to
Lawvere's ideas, but, as we said above, our approach is more
syntactical. In this syntactical approach substitution is not
mentioned. (What we call \emph{renaming} plays no role in it.)

Let $\QDS^{-x}$ be the full subcategory of \QDS\ whose objects are
all formulae of \eL\ in which $x$ is not free. From $\QDS^{-x}$ to
\QDS\ there is an obvious inclusion functor, which we call $E$.
(It behaves like identity on objects and on arrows.) The functor
$E$ is full and faithful. By restricting the codomain of the
functors $Q_x$ from \QDS\ to \QDS\ we obtain the functors $Q_x$
from \QDS\ to $\QDS^{-x}$. Then the functor $\forall_x$ is right
adjoint to $E$, and $\exists_x$ is left adjoint to $E$. Consider
first the adjunction involving $\forall_x$ and $E$. In this
adjunction, the arrows $\iota^{\forall\!_x}$ make the counit and
the arrows $\gamma^{\forall\!_x}$ the unit natural transformation,
while $(\forall\beta)$ and $(\forall\eta)$ are the triangular
equations of this adjunction. In the adjunction involving
$\exists_x$ and $E$, the arrows $\iota^{\exists_x}$ make the unit
and the arrows $\gamma^{\exists_x}$ the counit natural
transformation, while $(\exists\beta)$ and $(\exists\eta)$ are the
triangular equations. In other words, the full subcategory
$\QDS^{-x}$ of \QDS\ is both coreflective and reflective in \QDS.
The equations $(Q\gamma\iota)$ of Section 1.2 follow from
Theorem~1 and its dual in \cite{ML71} (Section IV.3). From these
theorems we also obtain that $Q_xA$ and $A$ are isomorphic for
every object $A$ of $\QDS^{-x}$.

Note that these two adjunctions, due to the presence of the
equations $(Q\gamma\iota)$, or $(Q\iota)$, or $(Q\gamma)$, of
Section 1.2, are \emph{trivial adjunctions} in the following
sense. If $f$ and $g$ of the same type are arrow terms of \QDS\
made only of $\mj$, $\cirk$, $Q_x$, $\iota^{Q\!_x}$ and
$\gamma^{Q\!_x}$, then ${f=g}$ in \QDS\ (see \cite{D99}, Sections
4.6.2 and 4.11). Whenever a full subcategory of a category $\cal
C$ is coreflective or reflective in $\cal C$, we have a trivial
adjunction in the same sense. (The notion of trivial adjunction is
closely related to Lambek's notion of idempotent monad of
\cite{L69}, Section~1.)

Note that the adjunctions involving $E$ and $Q_x$ do not deliver
the distributivity arrows $\check{\theta}^{\forall\!_x\str}$ and
$\hat{\theta}^{\exists_x\rts}$. Lawvere was able to define
$\hat{\theta}^{\exists_x\rts}$ in the presence of intuitionistic
implication, while the constant domain arrows
$\check{\theta}^{\forall\!_x\str}$ are not present, and not
desired in intuitionistic logic. In an analogous way, the
adjunctions of product and coproduct with the diagonal functor do
not deliver distributivity isomorphisms of coproduct over product
and of product over coproduct in bicartesian categories, i.e.\
categories that are cartesian and cocartesian. In bicartesian
closed categories, where we have the exponential functor, from
which intuitionistic implication arises, we obtain distributivity
isomorphisms of product over coproduct, but distributivity
isomorphisms of coproduct over product may be missing.

So Lawvere's thesis that logical constants are characterized
completely by adjoint situations should be taken with a grain of
salt. Conjunction, which corresponds to product, is characterized
by right-adjointness to the diagonal functor when it is alone, or
when it is accompanied by intuitionistic implication. When
conjunction and disjunction, which corresponds to coproduct, are
alone, then the two adjunctions with the diagonal functor do not
suffice. Some distributivity arrows, which we would have in the
presence of implication, are missing. The situation is analogous
with quantifiers and the distributivity arrows
$\check{\theta}^{\forall\!_x\str}$ and
$\hat{\theta}^{\exists_x\rts}$ in intuitionistic logic.

The situation is different in classical logic, where duality
reigns. Both of the distributivity arrows
$\check{\theta}^{\forall\!_x\str}$ and
$\hat{\theta}^{\exists_x\rts}$ are definable in the presence of
negation (see Section 2.7; negation yields implication and
``coimplication''). Both, when defined, happen to be isomorphisms
in this paper, and should be such in classical logic, but neither
the distribution of disjunction over conjunction nor the
distribution of conjunction over disjunction should be
isomorphisms in classical logic, as we argued in \cite{DP04}.

The following remark is not about our immediate concerns here, but
it is perhaps worth making once we have raised the issue of the
isomorphism of distribution of conjunction over disjunction, i.e.\
of product over coproduct. It is not clear that a modal
translation based on $S4$ will turn this distribution, which
should not be an isomorphism in classical logic, into an
isomorphism, as it should be in intuitionistic logic. So it is not
clear that in the proof theory of $S4$ based on classical logic we
will be able to represent correctly the proof theory of
intuitionistic logic, if the latter is based on bicartesian closed
categories. Equations between proofs need not be the same. A
similar phenomenon, pointed out in \cite{SEELY89} (to which the
referee brought our attention), is that the modal translation of
intuitionistic logic into linear logic based on Girard's modal
operator ! need not be proof-theoretically correct when
disjunction and the existential quantifier are taken into account.
The Kleisli category of the comonad of a Girard category in the
sense of \cite{SEELY89}, where one expects to find the modal
translation of intuitionistic logic, need not have coproducts and
a left adjoint to Lawvere's substitution functor based on
projection. This category need not be bicartesian closed.

\subsection{The category \GQDS}
In this section we enlarge the results of Section 7.7 of
\cite{DP04}, on which our exposition will heavily rely. We
introduce a category called \GQDS, which extends with quantifiers
the category \GDS\ of \cite{DP04}. In \GQDS\ we will be able to
perform in a manageable manner the \emph{Gentzenization} of \QDS\
(this is where \textbf{G} comes from).

Let a formula of \eL\ be called \emph{diversified} when every
predicate letter occurs in it at most once. A type ${A\vdash B}$
is called \emph{diversified} when $A$ and $B$ are diversified, and
an arrow term is \emph{diversified} when its type is diversified.

It is easy to verify that for every arrow ${f\!:A\vdash B}$ of
\QDS\ there is a diversified arrow term ${f'\!:A'\vdash B'}$ of
\QDS\ such that $f$ is obtained by substituting uniformly
predicate letters for some predicate letters in ${f'\!:A'\vdash
B'}$. Namely, $f$ is a letter-for-letter substitution instance of
$f'$ (cf.\ \cite{DP04}, Sections 3.3 and 7.6).

Our aim is to show that \QDS\ is a \emph{diversified preorder},
which means that if ${f_1,f_2\!:A\vdash B}$ are diversified arrow
terms, then ${f_1=f_2}$ in \QDS. For that purpose we introduce an
auxiliary category \GQDS\ where the $\b{\xi}{\str}$ arrows,
$\b{\xi}{\rts}$ arrows and $\c{\xi}$ arrows are identity arrows.
We will prove that \GQDS\ is a \emph{preorder}, which means that
for all arrow terms $f_1$ and $f_2$ of the same type we have
${f_1=f_2}$ in \GQDS. That \GQDS\ is a preorder will imply that
\QDS\ is a diversified preorder.

From the fact that \QDS\ is a diversified preorder one can infer
that there is a faithful functor $G$ from \QDS\ to the category
\emph{Rel}, whose objects are finite ordinals and whose arrows are
relations between these ordinals (see \cite{DP04}, Sections 2.9
and 7.6). This functor $G$ is defined as for \DS\ in \cite{DP04}
with the understanding that predicate letters now stand for
propositional letters; we have moreover that ${GQ_xA=GA}$ (so that
$GA$ is the number of occurrences of predicate letters in the
formula $A$), the arrow $G\alpha$ for $\alpha$ being
$\iota^{Q\!_x}_A$, $\gamma^{Q\!_x}_A$,
$\check{\theta}^{\forall\!_x\str}_{A,D}$ and
$\hat{\theta}^{\exists_x\rts}_{A,D}$ is an identity arrow, while
${GQ_xf=G[f]^x_y=Gf}$. The theorem that $G$ is a faithful functor
is called \QDS\ \emph{Coherence}.

Let $\eL^{div}$ be the set of diversified formulae of \eL.
Consider the smallest equivalence relation $\equiv$ on $\eL^{div}$
that satisfies
\begin{tabbing}
\hspace{7em}\=$A\ks(B\ks C)\equiv(A\ks B)\ks C$,\\[1ex]\> $A\ks
B\equiv B\ks A$,\\[1ex]\> if ${A_1\equiv B_1}$ and ${A_2\equiv
B_2}$, then ${A_1\ks A_2\equiv B_1\ks B_2}$,\\[1ex]\> if ${A\equiv
B}$, then ${Q_xA\equiv Q_xB}$,
\end{tabbing}
and let $[A]$ be the equivalence class of a diversified formula
$A$ with respect to this equivalence relation. We call $[A]$ a
\emph{form set} (which follows the terminology of \cite{DP04},
Section 7.7). We use $X$, $Y$, $Z,\ldots$, sometimes with indices,
for form sets. It is clear that the form set $[A]$ can be named by
any of the members of the equivalence class $[A]$. In these names
we may delete parentheses tied to $\!\ks\!$ in the immediate scope
of $\!\ks\!$. A \emph{subform set} of a form set $X$ is a form set
$[A]$ for $A$ a subformula of a formula in $X$.

Let the objects of the category \GQDS\ be the form sets we have
just introduced. The arrow terms of \GQDS\ are defined as the
arrow terms of \QDS\ save that their indices are form sets instead
of formulae. The equations of \GQDS\ are defined as those of \QDS\
save that we add the equations
\[
\begin{array}{l}
\b{\xi}{\str}_{X,Y,Z}\;=\;\b{\xi}{\rts}_{X,Y,Z}\;=\mj_{X\kst Y\kst
Z},
\\[1ex]
\c{\xi}_{X,Y}\;=\mj_{X\kst Y}.
\end{array}
\]
This defines the category \GQDS. From the fact that \GQDS\ is a
preorder we infer that \QDS\ is a diversified preorder as in
\cite{DP04} (Sections 3.3, 7.6, beginning of 7.7 and end of 7.8).

We define by induction a set of terms for the arrows of \GQDS,
which we call \emph{Gentzen terms}. They are defined as in
\cite{DP04} (Section 7.7), save that to the \emph{Gentzen
operations} ${cut}_X$, $\kon_{X_1,X_2}$ and $\vee_{X_1,X_2}$ we
add the following Gentzen operations, where $=_{dn}$ is read
``denotes'':

\[
\begin{array}{c}
\f{f\!:X^x_y\kon Z\vdash U} {\forall^L_{x,X}f=_{dn}
f\cirk([\iota^{\forall\!_x}_X]^x_y\kon\mj_Z)\!:\forall_xX\kon
Z\vdash U}
\\*[4ex]
\f{f\!:X^x_y\vdash U} {\forall^L_{x,X}f=_{dn}
f\cirk[\iota^{\forall\!_x}_X]^x_y\!:\forall_xX\vdash U}
\end{array}
\]

\[
\begin{array}{c}
\f{f\!:U\vdash X^x_u\vee Z} {\forall^R_{x,X}f=_{dn}
(\tau^{\forall\!_v}_{X^x_v,u,x}\vee\mj_Z)\cirk\check{\theta}^{\forall\!_u\str}_{X^x_u,Z}
\cirk\forall_uf\cirk\gamma^{\forall\!_u}_U\!:U\vdash\forall_xX\vee
Z}
\\*[4ex]
\f{f\!:U\vdash X^x_u} {\forall^R_{x,X}f=_{dn}
\tau^{\forall\!_v}_{X^x_v,u,x}\cirk\forall_uf
\cirk\gamma^{\forall\!_u}_U\!:U\vdash\forall_xX}
\end{array}
\]

\[
\begin{array}{c}
\f{f\!:U\vdash X^x_y\vee Z} {\exists^R_{x,X}f=_{dn}
([\iota^{\exists_x}_X]^x_y\vee\mj_Z)\cirk
f\!:U\vdash\exists_xX\vee Z}
\\*[4ex]
\f{f\!:U\vdash X^x_y} {\exists^R_{x,X}f=_{dn}
[\iota^{\exists_x}_X]^x_y\cirk f\!:U\vdash\exists_xX}
\end{array}
\]

\[
\begin{array}{c}
\f{f\!:X^x_u\kon Z\vdash U} {\exists^L_{x,X}f=_{dn}
\gamma^{\exists_u}_U\cirk\exists_uf\cirk
\hat{\theta}^{\exists_u\rts}_{X^x_u,Z}\cirk
(\tau^{\exists_v}_{X^x_v,x,u}\kon\mj_Z)\!:\exists_xX\kon Z\vdash
U}\\*[4ex] \f{f\!:X^x_u\vdash U} {\exists^L_{x,X}f=_{dn}
\gamma^{\exists_u}_U\cirk\exists_uf\cirk
\tau^{\exists_v}_{X^x_v,x,u}\!:\exists_xX\vdash U}
\end{array}
\]

\[
\f{f\!:X\vdash Y} {[f]^x_y\!:X^x_y\vdash Y^x_y}
\]
The usual \emph{proviso for the eigenvariable} in connection with
$\forall^R$ and $\exists^L$ is imposed by the tacit provisos
concerning $\gamma^{Q\!_u}_U$,
$\check{\theta}^{\forall\!_u\str}_{X^x_u,Z}$,
$\hat{\theta}^{\exists_u\rts}_{X^x_u,Z}$,
$\tau^{\forall\!_v}_{X^x_v,u,x}$ and
$\tau^{\exists_v}_{X^x_v,x,u}$. This proviso says that $u$, which
is called the \emph{eigenvariable}, is not free in the types of
$\forall^R_{x,X}f$ and $\exists^L_{x,X}f$; i.e., $u$ is free
neither in the sources nor in the targets.

The types of all the subterms of a Gentzen term make a derivation
tree usual in Gentzen systems. (An example may be found in Section
1.10.)

It is easy to show that every arrow of \GQDS\ is denoted by a
Gentzen term. For that we rely on the Gentzenization Lemma of
Section 7.7 of \cite{DP04}, together with the following equations
of \GQDS:
\begin{tabbing}
\hspace{14em}\= $\forall_xf$ \=\kill \> $\iota^{\forall\!_x}_X$\>
$=\forall^L_{x,X}\mj_X$,
\\[1ex]
\> $\gamma^{\forall\!_x}_U$\> $=\forall^R_{x,U}\mj_U$,
\\[1ex]
\> $\forall_xf$\> $=\forall^R_{x,Y}\forall^L_{x,X}f$,
\hspace{.5em}for $f\!:X\vdash Y$,
\end{tabbing}
and the dual equations involving $\exists$ instead of $\forall$.

With the help of the equations $(\forall\tau\iota)$ (together with
renaming), $(\forall\check{\theta}\iota)$, $\mbox{($\forall\iota$
{\it nat})}$ and $(\forall\beta)$ (see Sections 1.2-3) we derive
the following equations of \GQDS:
\begin{tabbing}
\hspace{11em}\=
$([\iota^{\forall\!_x}_X]^x_u\vee\mj_Z)\cirk\forall^R_{x,X}f=f\!:U\vdash
X^x_u\vee Z$,
\\[-.5ex]
\hspace{5em}$\mbox{($\forall\beta$ {\it red})}$
\\[-.5ex]
\> $[\iota^{\forall\!_x}_X]^x_u\cirk\forall^R_{x,X}f=f\!:U\vdash
X^x_u$.
\end{tabbing}
With the help of the naturality of
$\check{\theta}^{\forall\!_u\str}$, $\mbox{($\forall\tau$ {\it
nat})}$, $\mbox{({\it ren} 1)}$, $(\forall\iota)$,
$(\forall\tau\iota)$, $(\forall\check{\theta}\iota)$,
$\mbox{($\forall\iota$ {\it nat})}$ and $(\forall\beta)$ (see
Sections 1.2-3) we derive the following equations of \GQDS:
\begin{tabbing}
\hspace{11em}\=
$\forall^R_{x,X}(([\iota^{\forall\!_x}_X]^x_u\vee\mj_Z)\cirk
g)=g\!:U\vdash\forall_xX\vee Z$,
\\[-.5ex]
\hspace{5em}$\mbox{($\forall\eta$ {\it red})}$
\\[-.5ex]
\> $\forall^R_{x,X}([\iota^{\forall\!_x}_X]^x_u\cirk
g)=g\!:U\vdash\forall_xX$.
\end{tabbing}
We derive analogously the dual equation of \GQDS\ involving
$\exists$ instead of $\forall$, which are called
$\mbox{($\exists\beta$ {\it red})}$ and $\mbox{($\exists\eta$ {\it
red})}$.

\subsection{Variable-purification}
For proving the results of the following sections we need to
replace arbitrary Gentzen terms by Gentzen terms in whose type no
variable is both free and bound. This is the same kind of
condition that Kleene had to satisfy in \cite{K52} (Section 78) in
order to prove cut elimination in the predicate calculus. The
condition is implicit in Gentzen's \cite{G35}, because he did not
use the same letters for free and bound variables.

A variable $x$ is free in the type of ${f\!:X\vdash Y}$ when $x$
is free either in $X$ or in $Y$. We say that $x$
\emph{participates free in} $f$ when $x$ is free in the type of
some subterm of $f$. We have analogous definitions with ``free''
replaced by ``bound''. We say that $x$ \emph{participates in} $f$
when $x$ participates either free or bound in $f$. A Gentzen term
of \GQDS\ is \emph{variable-pure} when no variable participates in
it both free and bound.

By changing only bound variables one can transform an arbitrary
Gentzen term that is not variable-pure into a variable-pure
Gentzen term. (We could as well talk of \emph{renaming} of bound
variables, but, as we said at the beginning of Section 1.3, we do
not want to confuse this renaming with the renaming of free
variables.) The initial term and the resulting term need not be of
the same type, and hence need not be equal, but they will be equal
up to an isomorphism, as we shall see below.

Kleene's purification was done for a sequent where there was no
variable both free and bound, and his aim was to obtain a
derivation for it in which no variable is both free and bound. For
that he could not just change bound variables, but he needed also
to change free variables. Our aim is different, and we can change
only bound variables.

We have the following equations in \GQDS:
\begin{tabbing}
\hspace{11em}\=
$Q^L_{x,X}f=Q^L_{y,X^x_y}f\cirk(\tau^{Q\!_v}_{X^x_v,x,y}\kon\mj_Z)$,
\\[-.8ex]
\hspace{5em}$(Q^L\tau)$
\\[-.8ex]
\> $Q^L_{x,X}f=Q^L_{y,X^x_y}f\cirk \tau^{Q\!_v}_{X^x_v,x,y}$,
\\[3ex]
\> $Q^R_{x,X}f=(\tau^{Q\!_v}_{X^x_v,y,x}\vee\mj_Z)\cirk
Q^R_{y,X^x_y}f$,
\\[-.8ex]
\hspace{5em}$(Q^R\tau)$
\\[-.8ex]
\> $Q^R_{x,X}f=\tau^{Q\!_v}_{X^x_v,y,x}\cirk Q^R_{y,X^x_y}f$.
\end{tabbing}
To prove these equations we use essentially the equations
$(Q\tau\iota)$ and $\mbox{($Q\tau$ {\it trans})}$ of Section 1.3.

We define $\tau$-\emph{terms} inductively with the following
clauses:
\begin{itemize}
\item[] $\tau^{Q\!_x}_{X,u,v}$ is a $\tau$-term;\vspace{-1ex}
\item[] if $f$ is a $\tau$-term and $S$ is a quantifier prefix,
then ${f\ks\mj_Y}$ and $Sf$ are $\tau$-terms.
\end{itemize}
The unique subterm $\tau^{Q\!_x}_{X,u,v}$ of a $\tau$-term is
called its \emph{head}. Then for every $\tau$-term $h$ there is a
$\tau$-term $h'$ such that the following equations hold in \GQDS:
\begin{tabbing}
\hspace{1em}$({\!\ks\!} h)$ \emph{equations}:\\* \hspace{11em}\=
${\!\ks\!}_{X_1,X_2}(f\cirk
h,g)\,$\=$={\!\ks\!}_{X_1',X_2}(f,g)\cirk h'$,
\\*[1ex]
\> ${\!\ks\!}_{X_1,X_2}(h\cirk
f,g)$\>$=h'\cirk{\!\ks\!}_{X_1',X_2}(f,g)$,
\\[1ex]
\hspace{1em}$\mbox{({\it cut} $h$)}$ \emph{equations}:\\* \>
${cut}_X(f\cirk h,g)\,$\=$={cut}_X(f,g)\cirk h'$,
\\*[1ex]
\> ${cut}_X(f,g\cirk h)$\>$={cut}_X(f,g)\cirk h'$,
\\[1ex]
\> ${cut}_X(h\cirk f,g)$\>$=h'\cirk{cut}_X(f,g)$,
\\[1ex]
\> ${cut}_X(f,h\cirk g)$\>$=h'\cirk{cut}_X(f,g)$,
\\[1ex]
\> ${cut}_X(f\cirk h_1,h_2\cirk g)={cut}_{X'}(f,g)$,
\\[2ex]
\hspace{1em}$(Q h)$ \emph{equations}: \hspace{1em}for
$S\in\{L,R\}$,\\*[1ex] \> $Q^S_{x,X}(f\cirk
h)\,$\=$=Q^S_{x,X'}f\cirk h'$,
\\*[1ex]
\> $Q^S_{x,X}(h\cirk f)$\>$=h'\cirk Q^S_{x,X'}f$,
\\[1ex]
\hspace{1em}$\mbox{({\it ren} $h$)}$ \emph{equations}:\\* \>
$[f\cirk h]^x_y\,$\=$=[f]^x_y\cirk h'$,
\\*[1ex]
\> $[h\cirk f]^x_y$\>$=h'\cirk[f]^x_y$.
\end{tabbing}
In these equations $X'$ is either $X$ or a different form set. The
$\tau$-terms $h_1$ and $h_2$ in the last $\mbox{({\it cut} $h$)}$
equation differ in their heads, which are inverse to each other
(see Section 1.3).

To derive the $({\!\ks\!} h)$ equations and the first four
$\mbox{({\it cut} $h$)}$ equations we use essentially functorial
and naturality equations (see \cite{DP04}). For the last
$\mbox{({\it cut} $h$)}$ equation we also use $\mbox{($Q\tau$ {\it
sim})}$, and for the $(Qh)$ and $\mbox{({\it ren} $h$)}$ equations
we use essentially $\mbox{($Q\tau$ {\it ren})}$, $\mbox{($Q\iota$
{\it nat})}$ and $\mbox{($Q\tau$ {\it nat})}$ (see Sections
1.2-3).

By applying the equations of \GQDS\ mentioned in this section, we
can establish the following.

\prop{Variable-Purification Lemma}{For every Gentzen term
${f\!:X\vdash Y}$ there is a variable-pure Gentzen term
${f'\!:X'\vdash Y'}$ such that in \GQDS
\[
f=h_2\cirk f'\cirk h_1
\]
where $h_1$ and $h_2$ are compositions of $\tau$-terms or $\mj_X$
or $\mj_Y$.}

\noindent Let us explain up to a point how we achieve that.

Let $x$ be \emph{new} for $f$ when $x$ does not participate in $f$
(see the beginning of the section) and does not occur as an index
in the Gentzen operations of renaming that occur in $f$. Let
${x_1,\ldots,x_n}$ be all the variables that participate bound in
$f$. Then take the variables ${x_1',\ldots,x_n'}$ all new for $f$,
and apply first the equations $(Q^L\tau)$ and $(Q^R\tau)$ with $x$
being $x_i$ and $y$ being $x_i'$. In Gentzen terms, $\mj_X$ occurs
only with $X$ atomic, and so for every variable that participates
bound in $f$ there is a Gentzen operation by which it was
introduced. It remains then to apply the equations $({\!\ks\!}
h)$, $\mbox{({\it cut} $h$)}$, $(Qh)$ and $\mbox{({\it ren}
$h$)}$. Note that $h_1$ and $h_2$ depend only on the type of $f$
and on the choice of the variables ${x_1',\ldots,x_n'}$.

\subsection{Renaming of eigenvariables}
In this section we prove the equations of \GQDS\ of the following
form:
\begin{tabbing}
\hspace{5em}\mbox{($Q^S$ {\it ren})}\hspace{3em}
$Q^S_{x,X}f=Q^S_{x,X}[f]^u_v$
\end{tabbing}
for ${Q^S\in\{\forall^R,\exists^L\}}$, with $u$ the eigenvariable
and $v$ a variable not free in the type of $f$ (see the beginning
of the preceding section). That $v$ is not free in the type of $f$
is satisfied a fortiori when $v$ is new for the left-hand side
(see the end of the preceding section). The equations
$\mbox{($Q^S$ {\it ren})}$ say that \GQDS\ covers the renaming of
eigenvariables by new variables, which is a technique derived from
\cite{G35} (Section III.3.10). We need the equations $\mbox{($Q^S$
{\it ren})}$ to prove the results of Sections 1.9-10.

We derive now the equation $\mbox{($\forall^R$ {\it ren})}$ for
${f\!:U\vdash X^x_u\vee Z}$:
\begin{tabbing}
\hspace{1em}\= $\forall^R_{x,X}f\,$\=
$=((\tau^{\forall\!_w}_{X^x_w,v,x}\cirk
\tau^{\forall\!_w}_{X^x_w,u,v})\vee\mj_Z)\cirk
\check{\theta}^{\forall\!_u\str}_{X^x_u,Z}\cirk
\forall_u[[f]^u_w]^w_u\cirk\gamma^{\forall\!_u}_U$, by
$\mbox{({\it ren} 6)}$, $\mbox{({\it ren} 2)}$
\\*
\`  and $\mbox{($\forall\tau$ {\it trans})}$,
\\[1ex]
\> \> $=(\tau^{\forall\!_w}_{X^x_w,v,x}\vee\mj_Z)\cirk
\check{\theta}^{\forall\!_v\str}_{X^x_v,Z}\cirk
\forall_w[[f]^u_w]^w_v\cirk \tau^{\forall\!_w}_{U,u,v}\cirk
\gamma^{\forall\!_u}_U$, by $(\forall\tau\check{\theta})$ and
$\mbox{($\forall\tau$ {\it nat})}$,
\\[1ex]
\> \> $=\forall^R_{x,X}[f]^u_v$, by $\mbox{({\it ren} 5)}$,
$\mbox{({\it ren} 2)}$ and $(\forall\tau\gamma)$ (see Sections
1.2-3).
\end{tabbing}
The equation $\mbox{($\exists^L$ {\it ren})}$ is derived
analogously.

We can prove also the equations
\begin{tabbing}
\hspace{5em}\mbox{($Q^T$ {\it ren})}\hspace{3em}
$Q^T_{x,X}f=Q^T_{x,X}[f]^y_z$
\end{tabbing}
for ${Q^T\in\{\forall^L,\exists^R\}}$, with $f$ either of the type
${X^x_y\kon Z\vdash U}$ or ${X^x_y\vdash U}$ or of the type
${U\vdash X^x_y\vee Z}$ or ${U\vdash X^x_y}$, provided $y$ is not
free in the type of $Q^T_{x,X}f$. So, though $y$ is not here an
eigenvariable, it could have been one.

To derive the equation $\mbox{($\forall^L$ {\it ren})}$ for
${f\!:X^x_y\kon Z\vdash U}$ we have
\begin{tabbing}
\hspace{8em}\=
$\forall^L_{x,X}f\,$\=$=[f\cirk([\iota^{\forall\!_x}_X]^x_y\kon
\mj_Z)]^y_z$, \hspace{1em}by $\mbox{({\it ren} 2)}$,\\*[1ex]
 \>\>$=[f]^y_z\cirk([\iota^{\forall\!_x}_X]^x_z\kon
\mj_Z)=\forall^L_{x,X}[f]^y_z$,
\end{tabbing}
by using, together with other renaming equations, $\mbox{({\it
ren} 5)}$ and $\mbox{({\it ren} 2)}$ if $y$ is different from $x$,
since $y$ is then not free in the type of $\iota^{\forall\!_x}_X$,
and by using $\mbox{({\it ren} 1)}$ if $y$ is $x$. The equation
$\mbox{($\exists^R$ {\it ren})}$ is derived analogously.

\subsection{Elimination of renaming}
We can establish the following proposition for \GQDS.

\prop{Renaming Elimination}{For every variable-pure and cut-free
Gentzen term $t$ there is a variable-pure, cut-free and
renaming-free Gentzen term $t'$ such that ${t=t'}$.}

\noindent Here \emph{cut-free} means of course that no instance of
the Gentzen operation ${cut}_X$ occurs in $t$ and $t'$, and
\emph{renaming-free} means that none of the Gentzen operations
${[\;\;]^x_y}$ occurs in $t'$.

The proof of Renaming Elimination is based on the following
equations of \GQDS:
\begin{tabbing}
\hspace{9em}\=
$[{\!\ks\!}_{X_1,X_2}(f,g)]^x_y={\!\ks\!}_{X_1',X_2'}([f]^x_y,[g]^x_y)$,
\\[1.5ex]
\> $[Q^S_{z,X}f]^x_y=Q^S_{z,X'}[f]^x_y$, if $z$ is neither $x$ nor
$y$.
\end{tabbing}
To eliminate all occurrences of renaming we eliminate one by one
innermost occurrences of renaming, i.e.\ occurrences of renaming
within the scope of which there is no renaming. Variable-purity
ensures that the proviso of the second equation is not an
obstacle.

We will use Renaming Elimination for the proof of the
Cut-Elimination Theorem for \GQDS\ in the next section. For that
we need a \emph{strengthened version of Renaming Elimination}, in
which it is specified that the Gentzen term $t'$ is exactly
analogous to $t$: only indices of its identity arrows and of its
Gentzen operations may change.

For ${f\!:X\vdash Y}$ and ${g\!:Y\vdash X}$ such that $x$ is not
free in $X$ in \GQDS\ we have
\begin{tabbing}
\hspace{11em}\=$[f]^x_y\,$\=$={cut}_{\forall\!_xY}
(\forall^R_{x,Y}f,\forall^L_{x,Y}\mj_{Y^x_y})$,\\*[1.5ex]
\>$[g]^x_y$\>$={cut}_{\exists_xY}(\exists^R_{x,Y}\mj_{Y^x_y},\exists^L_{x,X}g)$.
\end{tabbing}
So particular instances of renaming (and Gentzen and Kleene did
not envisage implicitly more than that) can be easily eliminated
provided we want to tolerate cut. (In the presence of implication
we could eliminate all instances of renaming in the presence of
cut, as we mentioned in Section 1.2.) Our aim however is to
eliminate both cut and renaming.

If we delete ``variable-pure'' from Renaming Elimination, then
this proposition cannot be proved. A counterexample, analogous to
a counterexample in \cite{K52} (Section 78, Example~4), is the
following:
\[
[\forall^L_{x,\forall_yRxy}\forall^L_{y,Ruy}\mj_{Ruz}]^u_y\!:
\forall_x\forall_yRxy\vdash Ryz.
\]
From this Gentzen term we can eliminate renaming only by
introducing cut, as above.

Kleene in \cite{K52} also needed variable-purity to eliminate cut.
But his counterexample, mentioned above, which is analogous to our
counterexample, would not be a counterexample in the presence of
renaming.

\subsection{Cut elimination}
Our aim in this section is to establish the following theorem for
\GQDS.

\prop{Cut-Elimination Theorem}{For every variable-pure Gentzen
term $t$ there is a variable-pure and cut-free Gentzen term $t'$
such that ${t=t'}$.}

\noindent The proof of this theorem is obtained by modifying and
expanding the proof of the Cut-Elimination Theorem for \GDS\ in
\cite{DP04} (Section 7.7). We presuppose below the terminology
introduced in this previous proof.

The $Q$-\emph{rank} of ${cut}_{Q\!_xX}(f,g)$ is ${n_1\pl n_2}$
when $f$ has a subterm $Q^R_{x,X}f'$ of depth $n_1$ and $g$ has a
subterm $Q^L_{x,X}g'$ of depth $n_2$. The \emph{rank} of a topmost
cut ${cut}_X(f,g)$ is either its $\kon$-rank, or $\vee$-rank, or
$p$-rank, or $Q$-rank depending on $X$.

The \emph{complexity} of a topmost cut ${cut}_X(f,g)$ is ${(m,n)}$
where ${m\geq 1}$ is the sum of the number of predicate letters
and occurrences of quantifier prefixes in $X$ and ${n\geq 0}$ is
the rank of this cut. Every form set of the form $Q_xX$ is
considered to be both of colour $\kon$ and colour $\vee$.

In the proof we have the following additional cases. We consider
only the most complicated cases, and leave out the remaining
simpler cases, which are dealt with analogously.

\vspace{2ex}

\noindent $(\forall 1)$ \quad If our topmost cut is
\[
{cut}_{\forall\!_xX}(\forall^R_{x,X}f,\forall^L_{x,X}g)\!:U\kon
Y\vdash Z\vee W
\]
for ${f\!:U\vdash X^x_u\vee Z}$ and ${g\!:X^x_v\kon Y\vdash W}$,
with complexity ${(m,0)}$ where ${m>1}$, then we use the equation
\[
{cut}_{\forall\!_xX}(\forall^R_{x,X}f,\forall^L_{x,X}g)={cut}_{X^x_v}([f]^u_v,g),
\]
in which the cut on the right-hand side is of lower complexity
than the topmost cut on the left-hand side. To derive this
equation we use essentially  the equations $\mbox{($\forall\beta$
{\it red})}$ (see Section 1.5) together with naturality and
functorial equations. We proceed analogously when the topmost cut
we start from is
\[
{cut}_{\exists_xX}(\exists^R_{x,X}f,\exists^L_{x,X}g).
\]

Suppose for the cases below that $X$ is of colour $\kon$.
\vspace{2ex}

\noindent $(\forall 2)$ \quad If our topmost cut is
\[
{cut}_X(\forall^R_{x,V}f,g)\!:U\kon Y\vdash \forall_xV\vee Z\vee W
\]
for ${f\!:U\vdash X\vee V^x_u\vee Z}$ and ${g\!:X\kon Y\vdash W}$,
with complexity ${(m,n)}$ where ${m,n\geq 1}$, then we use the
equation
\[
{cut}_X(\forall^R_{x,V}f,g)=\forall^R_{x,V}{cut}_X([f]^u_v,g)
\]
with $v$ being a variable new for the left-hand side. By the
strengthened version of Renaming Elimination from the preceding
section, there is a variable-pure, cut-free and renaming-free
Gentzen term $f'$ such that ${[f]^u_v=f'}$, and the complexity of
${cut}_X(f',g)$ is ${(m,n\mn 1)}$. This equation is derived as
follows:
\begin{tabbing}
\hspace{1em}\= ${cut}_X(\forall^R_{x,V}f,g)$ \=
$={cut}_X(\forall^R_{x,V}[f]^u_v,g)$, by $\mbox{($\forall^R$ {\it
ren})}$ for $v$ new for the left-hand\\*\` side (see Section 1.7),
\\[1ex]
\> \> $=\forall^R_{x,V}(([\iota^{\forall\!_y}_V]^y_v\vee\mj_{Z\vee
W})\cirk{cut}_X(\forall^R_{x,V}[f]^u_v,g))$, by
$\mbox{($\forall\eta$ {\it red})}$
\\*
\`(see Section 1.5),
\\[1ex]
\> \>
$=\forall^R_{x,V}({cut}_X(([\iota^{\forall\!_x}_V]^x_v\vee\mj_Z)\cirk\forall^R_{x,V}[f]^u_v,g)$,
by functorial and
\\*
\` naturality equations,
\\[1ex]
\> \> $=\forall^R_{x,V}({cut}_X([f]^u_v,g))$, by
$\mbox{($\forall\beta$ {\it red})}$.
\end{tabbing}
We needed to rename the eigenvariable $u$ by a new $v$ in order to
ensure that the proviso for the eigenvariable is satisfied in the
second line for the $\forall^R_{x,V}$ operation newly introduced.
\vspace{2ex}

\noindent $(\forall 3)$ \quad If our topmost cut is
\[
{cut}_X(\forall^L_{x,V}f,g)\!:\forall_xV\kon U\kon Y\vdash Z\vee W
\]
for ${f\!:V^x_y\kon U\vdash X\vee Z}$ and ${g\!:X\kon Y\vdash W}$,
then we use the straightforward equation
\begin{tabbing}
\hspace{10em}\=$
{cut}_X(\forall^L_{x,V}f,g)\,$\=$=\forall^L_{x,V}{cut}_X(f,g)$.
\end{tabbing}
We have also the straightforward equation
\begin{tabbing}
\hspace{10em}\=$
{cut}_X(\forall^L_{x,V}f,g)\,$\=$=\forall^L_{x,V}{cut}_X(f,g)$.\kill

\>${cut}_X(\exists^R_{x,V}f,g)$\>$=\exists^R_{x,V}{cut}_X(f,g)$,
\end{tabbing}
and the equation
\begin{tabbing}
\hspace{10em}\=$
{cut}_X(\forall^L_{x,V}f,g)\,$\=$=\forall^L_{x,V}{cut}_X(f,g)$.\kill

\>${cut}_X(\exists^L_{x,V}f,g)$\>$=\exists^L_{x,V}{cut}_X([f]^u_v,g)$,
\end{tabbing}
proved like the analogous equation in case $(\forall 2)$. These
equations enable us to settle the remaining cases when $X$ is of
colour $\kon$. When $X$ is of colour $\vee$, we proceed in a dual
manner.

By Renaming Elimination, we need not consider cases when in our
topmost cut ${cut}_X(f,g)$ either $f$ or $g$ is of the form
$[h]^y_z$.

\subsection{Invertibility in \GQDS}
The results we are going to prove in this section correspond to
inverting rules in derivations, i.e.\ passing from conclusions to
premises. This invertibility is guaranteed by the possibility to
permute rules, i.e.\ change their order in derivations, and we
show for that permuting that it is covered by the equations of
\GQDS. (Permutation of rules is a theme treated in \cite{K52a},
but without considering equations between derivations.)

Besides the equations mentioned in \cite{DP04} (beginning of
Section 7.8) we will need the equations of \GQDS\ of the following
form:
\begin{tabbing}
\hspace{5em}$(\!\ks\!Q^S)$\hspace{2em}
${\!\ks\!}_{X_1,X_2}(Q^S_{x,X}f_1,f_2)=Q^S_{x,X}{\!\ks\!}_{X_1,X_2}(f_1',f_2)$,
\end{tabbing}
for $X^x_y$ being a subform set of the source or target of $f_1$,
and $f_1'$ being $f_1$ when ${Q^S\in\{\forall^L,\exists^R\}}$, and
$[f_1]^y_v$ with $v$ new for the left-hand side when
${Q^S\in\{\forall^R,\exists^L\}}$. These equations are either
straightforward to derive, or when
${Q^S\in\{\forall^R,\exists^L\}}$ we derive them by imitating the
derivation of the equation of case $(\forall 2)$ of the preceding
section, with the help of the equations $\mbox{($Q\beta$ {\it
red})}$ and $\mbox{($Q\eta$ {\it red})}$ (see the end of Section
1.5).

We will also need for the end of the section the following
equations, whose derivations are not difficult to find:
\begin{tabbing}
for ${Q^S\in\{\forall^L,\exists^R\}}$,
\\[1ex]
\hspace{5em}\= $(Q^SQ^S)$\hspace{3em}\=
$Q^S_{y,Y}Q^S_{x,X}f=Q^S_{x,X}Q^S_{y,Y}f$,
\\[1.5ex]
\> $(\exists^R\forall^L)$\>
$\exists^R_{y,Y}\forall^L_{x,X}f=\forall^L_{x,X}\exists^R_{y,Y}f$,
\\[2ex]
for $Q\in\{\forall,\exists\}$, $S\in\{L,R\}$, and the proviso for
the eigenvariable being satisfied,
\\[1ex]
\> $(Q^RQ^L)$\> $Q^R_{y,Y}Q^L_{x,X}f=Q^L_{x,X}Q^R_{y,Y}f$,
\\[1.5ex]
\> $(\exists^S\forall^S)$\>
$\exists^S_{y,Y}\forall^S_{x,X}f=\forall^S_{x,X}\exists^S_{y,Y}f$.
\end{tabbing}

The Invertibility Lemmata for $\kon$ and $\vee$ are formulated as
in \cite{DP04} (Section 7.8). Only ${let}(X)$ is the set of
\emph{predicate} letters occurring in the form set $X$. These
lemmata hold also when we replace throughout ``cut-free Gentzen
term'' by ``variable-pure, cut-free and renaming-free Gentzen
term''. They are proved as in \cite{DP04}, with additional cases
covered by the equations $(\!\ks\!Q^S)$.

The following invertibility lemmata are easy consequences of the
equations $\mbox{($\forall\eta$ {\it red})}$ and
$\mbox{($\exists\eta$ {\it red})}$ (see the end of Section 1.5).

\prop{Invertibility Lemma for $\forall^R$}{If $f$ is a
variable-pure Gentzen term of the type ${U\vdash\forall_xX\vee Z}$
or ${U\vdash\forall_xX}$, then there is a variable-pure Gentzen
term $f'$ of the type ${U\vdash X^x_u\vee Z}$ or ${U\vdash X^x_u}$
respectively such that ${\forall^R_{x,X}f'=f}$.}

\vspace{-2ex}

\prop{Invertibility Lemma for $\exists^L$}{If $f$ is a
variable-pure Gentzen term of the type ${\exists_xX\kon Z\vdash
U}$ or ${\exists_xX\vdash U}$, then there is a variable-pure
Gentzen term $f'$ of the type ${X^x_u\kon Z\vdash U}$ or
${X^x_u\vdash U}$ respectively such that ${\exists^L_{x,X}f'=f}$.}

Before formulating the remaining invertibility lemmata for
$\forall^L$ and $\exists^R$ we must introduce a number of notions
concerning occurrences of variables within the types of
variable-pure, cut-free and renaming-free Gentzen terms. Although
many of the notions introduced make sense also for other Gentzen
terms, we need these notions only in the context of variable-pure,
cut-free and renaming-free Gentzen terms. The essential assertions
using these notions, which we need for our results, need not hold
for all Gentzen terms.

Let $\alpha$ and $\beta$, sometimes with indices, stand for
occurrences of individual variables in a form set, and let
$\gamma$, sometimes with indices, stand for an occurrence of a
quantifier prefix in a form set. Let $\alpha$ and $\beta$ be
different occurrences of the variable $x$ in the form set $X$, and
let $\gamma$ be an occurrence of $Q_x$ in $X$. Then we say that
$\alpha$ and $\beta$ are \emph{simultaneously bound by} $\gamma$
when $X$ has a subform set $\gamma Y$ such that $\alpha$ and
$\beta$ are free in $Y$.

We say that a predicate letter $P$ \emph{occurs in the type} of
the Gentzen term ${f\!:X_1\vdash X_2}$ when it occurs in $X_1$ or
$X_2$. Because of diversification (see the beginning of Section
1.5), every $n$-ary predicate letter $P$ that occurs in the type
of $f$ occurs exactly once in $X_i$ in a subform set
${P\alpha^i_1\ldots\alpha^i_n}$ of $X_i$, for ${i\in\{1,2\}}$.
Here $\alpha^1_j$ and $\alpha^2_j$, for ${1\leq j\leq n}$, are not
necessarily occurrences of the same variable.

We say that the pair
${(P\alpha^1_1\ldots\alpha^1_n,P\alpha^2_1\ldots\alpha^2_n)}$ is a
\emph{formula couple} of $f$, and we say that
${(\alpha^1_j,\alpha^2_j)}$ is a \emph{couple} of $f$. We say that
it is the $P_j$-\emph{couple} of $f$ when we want to stress from
which formula couple and from which place $j$ in it originates. If
${(\alpha^1,\alpha^2)}$ is a couple of $f$, then $\alpha^i$, for
${i\in\{1,2\}}$, is \emph{free} in ${(\alpha^1,\alpha^2)}$ when it
is free in $X_i$, and analogously with ``free'' replaced by
``bound'', ``universally bound'' and ``existentially bound''. For
example, if $f$ is a Gentzen term of type
\[
\forall_xRux\kon Py\vdash\exists_z(Ruz\kon Pz),
\]
then $(Rux,Ruz)$ and $(Py,Pz)$ are the formula couples of $f$, and
the $R_1$-couple of $f$ consists of the occurrences of $u$ in
$Rux$ and $Ruz$, the $R_2$-couple of $f$ consists of the
occurrence of $x$ in $Rux$ and the occurrence of $z$ in $Ruz$, and
the $P_1$-couple of $f$ consists of the occurrence of $y$ in $Py$
and the occurrence of $z$ in $Pz$ in the type of $f$.

A \emph{left bridge} between the different couples
${(\alpha^1,\alpha^2)}$ and ${(\beta^1,\beta^2)}$ of $f$ such that
$\alpha^1$ and $\beta^1$ are occurrences of the same variable $x$
is an occurrence $\gamma$ of a quantifier prefix $Q_x$ in $X_1$
such that $\alpha^1$ and $\beta^1$ are simultaneously bound by
$\gamma$. We define analogously a \emph{right bridge} by replacing
$\alpha^1$, $\beta^1$ and $X_1$ with $\alpha^2$, $\beta^2$ and
$X_2$ respectively. A \emph{bridge} is a left bridge or a right
bridge. For $f$ as above, we have that the occurrence of
$\exists_z$ in the type of $f$ is the right bridge between the
$R_2$-couple and the $P_1$-couple of $f$, while the occurrence of
$\forall_x$ is not a bridge at all.

Two different couples of $f$ are \emph{bridgeable} when there is a
bridge between them (there might be both bridges, as, for example,
between the $R_1$-couple and $R_2$-couple of a Gentzen term of
type $\forall_xRxx\vdash\exists_yRyy$). For ${n\geq 2}$, a
sequence of couples
\[
(\alpha^1_1,\alpha^2_1),(\alpha^1_2,\alpha^2_2),\ldots,(\alpha^1_n,\alpha^2_n)
\]
of $f$ such that ${(\alpha^1_j,\alpha^2_j)}$ and
${(\alpha^1_{j+1},\alpha^2_{j+1})}$, where ${1\leq j\leq n\mn 1}$,
are bridgeable is called a \emph{bridgeable chain of couples}. For
every bridgeable chain of couples we can find at least one
sequence of bridges ${\gamma_1,\ldots,\gamma_{n-1}}$ that ensure
its bridgeability.

We say that ${(\alpha^1_1,\alpha^2_1)}$ and
${(\alpha^1_n,\alpha^2_n)}$ are \emph{clustered} when there is a
bridgeable chain of couples in which ${(\alpha^1_1,\alpha^2_1)}$
and ${(\alpha^1_n,\alpha^2_n)}$ are respectively the first and
last member. A sequence of bridges ensuring the bridgeability of
this bridgeable chain of couples is said to \emph{ensure the
clustering} of ${(\alpha^1_1,\alpha^2_1)}$ and
${(\alpha^1_n,\alpha^2_n)}$.

A set $C$ of couples of $f$ is a \emph{cluster} of $f$ when there
is a couple ${(\alpha^1,\alpha^2)}$  in $C$ such that for every
couple ${(\beta^1,\beta^2)}$ of $f$ different from
${(\alpha^1,\alpha^2)}$ we have that ${(\beta^1,\beta^2)\in C}$
iff ${(\beta^1,\beta^2)}$ is clustered with
${(\alpha^1,\alpha^2)}$. If $f$ is of type $\forall_xRux\kon
Py\vdash\exists_z(Ruz\kon Pz)$, then we have two clusters of $f$:
one singleton consisting of the $R_1$-couple, and another whose
elements are the $R_2$-couple and $P_1$-couple of $f$. If $f$ is
of type $\forall_xRxx\kon Py\vdash\exists_z(Ryz\kon Pz)$, then we
have just one cluster of $f$, since the $R_1$-couple, $R_2$-couple
and $P_1$-couple of $f$ make a bridgeable chain of couples.

Let $P$ be a predicate letter occurring in the type of a cut-free
Gentzen term $f$, and let ${Sub}_P(f)$ be the set of subterms of
$f$ in whose type $P$ occurs. For every member $f'$ of
${Sub}_P(f)$ there is a formula couple
${(P\alpha^1_1\ldots\alpha^1_n,P\alpha^2_1\ldots\alpha^2_n)}$ of
$f'$. The set of all the couples ${(\alpha^1,\alpha^2)}$ such that
there is a member $f'$ of ${Sub}_P(f)$ with
${(\alpha^1,\alpha^2)}$ the $P_j$-couple of $f'$ is called an
\emph{arc} of $f$. For example, in the following picture:

\begin{center}
\begin{picture}(250,90)
\put(125,40){\makebox(0,0){$\fS{\fS{\fS{\mj_{Ruy}\!:Ru\encircle{y}\vdash
Ru\encircle{y}}{\forall^L_{x,Rux}\mj_{Ruy}\!:\forall_{x}Ru\encircle{x}\vdash
Ru\encircle{y}}\quad\quad\quad\afrac{
\mj_{Py}\!:P\encircle{y}\vdash
P\encircle{y}}}{\kon_{Ruy,Py}(\forall^L_{x,Rux}\mj_{Ruy},\mj_{Py})\!:\forall_{x}Ru\encircle{x}\kon
P\encircle{y}\vdash Ru\encircle{y}\kon
P\encircle{y}}}{\exists^R_{z,Ruz\kon
Pz}\kon_{Ruy,Py}(\forall^L_{x,Rux}\mj_{Ruy},\mj_{Py})\!:
\forall_{x}Ru\encircle{x}\kon P\encircle{y}\vdash \exists_{z}
(Ru\encircle{z}\kon P \encircle{z})}$}}

 \put(75.5,67.5){\line(1,-1){12}}
 \put(113.5,67.5){\line(1,-1){12}}
 \put(97,46){\line(3,-1){44}}
 \put(135,46){\line(5,-1){74}}
 \put(151.3,22){\line(2,-3){7}}
 \put(220,22){\line(2,-1){21}}
 \put(91,80){\oval(38,20)[t]}
 \put(200.5,46){\line(-3,-2){20}}
 \put(200,46){\line(-3,-2){20}}
 \put(199.5,46){\line(-3,-2){20}}
 \put(240,46){\line(2,-5){5}}
 \put(240.3,46){\line(2,-5){5}}
 \put(239.7,46){\line(2,-5){5}}
 \put(182.7,22){\line(2,-3){7}}
 \put(183,22){\line(2,-3){7}}
 \put(183.3,22){\line(2,-3){7}}
 \put(251,22){\line(2,-1){20}}
 \put(251.5,22){\line(2,-1){20}}
 \put(250.5,22){\line(2,-1){20}}
 \put(220,58){\oval(33,20)[t]}

\put(220,58){\oval(32,19)[t]}

\put(220,58){\oval(33.5,20.5)[t]}

\end{picture}
\end{center}
the encircled occurrences of variables connected by lines make two
arcs of the variable-pure, cut-free and renaming-free Gentzen term
in the last line.

We say that an arc is the $P_j$-\emph{arc} when we want to stress
from which formula couples and from which place $j$ in them it
originates. In our example above, we have drawn the $R_2$-arc and
the $P_1$-arc. The \emph{bottom} of the $P_j$-arc of $f$ is the
$j$-th coordinate ${(\alpha^1_j,\alpha^2_j)}$ of the formula
couple
${(P\alpha^1_1\ldots\alpha^1_n,P\alpha^2_1\ldots\alpha^2_n)}$ of
$f$. In the example above, the bottom of the $R_2$-arc is
${(x,z)}$ in the last line and the bottom of the $P_1$-arc is
${(y,z)}$ in the last line.

Two arcs of $f$ are \emph{clustered} when their bottoms are
clustered. In the example above, the two arcs are clustered. If we
delete the last line, then we obtain two arcs of
\[
\kon_{Ruy,Py}(\forall^L_{x,Rux}\mj_{Ruy},\mj_{Py})
\]
that are not clustered.

A set $\cal{A}$ of arcs of $f$ is an \emph{arc-cluster} of $f$
when there is an arc $a$ in $\cal{A}$ such that for every arc $b$
of $f$ different from $a$ we have that ${b\in\cal{A}}$ iff $b$ is
clustered with $a$. The bottoms of the arcs in an arc-cluster of
$f$ make a cluster of $f$. We call this cluster the \emph{bottom
cluster} of the arc-cluster. In our example, we have an
arc-cluster, whose bottom cluster is made of the bottom of the
$R_2$-arc and the bottom of the $P_1$-arc.

All occurrences of variables that are free in the couples of an
arc, or of an arc-cluster, of a cut-free and renaming-free Gentzen
term $f$ are occurrences of the same variable. (Here cut-freedom
and renaming-freedom is essential.) This variable is called the
\emph{free variable of} the arc, or of the arc-cluster. The free
variable of the arc-cluster in our example is $y$.

For ${S\in\{L,R\}}$, suppose we have a cut-free and renaming-free
Gentzen term $f$ that has a subterm $Q^S_{x,X}g$ for $x$ free in
$X$ and $X^x_y$ occurring in the type of $g$. We say that
$Q^S_{x,X}g$ \emph{belongs to an arc-cluster} of $f$ when the
occurrences of $y$ in $X^x_y$ in the type of $g$ that have
replaced $x$ in $X$ belong to couples in this arc-cluster. We say
that $Q^S_{x,X}g$ \emph{belongs to a cluster} of $f$ when it
belongs to an arc-cluster whose bottom cluster is this cluster. In
our example above, the Gentzen term in the last line:
\[
\exists^R_{z,Ruz\kon
Pz}\kon_{Ruy,Py}(\forall^L_{x,Rux}\mj_{Ruy},\mj_{Py})
\]
belongs to the arc-cluster consisting of the $R_2$-arc and
$P_1$-arc, and the Gentzen term $\forall^L_{x,Rux}\mj_{Ruy}$ in
the second line  belongs to the arc-cluster consisting of the
$R_1$-arc, which is not drawn in the picture.

The subterms of $f$ that belong to an arc-cluster of $f$ are
called the \emph{gates} of that arc-cluster, and analogously with
``arc-cluster'' replaced by ``cluster''. A gate belonging to a
cluster may correspond to bridges in bridgeable chains of couples
in this cluster, but it need not correspond to such a bridge. In
our example above, the Gentzen term in the last line corresponds
to a bridge, but $\forall^L_{x,Rux}\mj_{Ruy}$ does not.

A gate is called an \emph{eigengate} when it is either of the
$\forall^R$ or of the $\exists^L$ type. As a consequence of the
proviso for the eigenvariable, we obtain that if a cluster has an
eigengate $Q^S_{x,X}g$, then every other gate of that cluster is a
subterm of $g$. This implies that every cluster has at most one
eigengate. As another consequence of the proviso for the
eigenvariable, we have the following remark.

\prop{Eigengate Remark}{If we have an eigengate in a cluster, then
for every couple ${(\alpha^1,\alpha^2)}$ in this cluster both
$\alpha^1$ and $\alpha^2$ are bound.}

A couple ${(\alpha^1,\alpha^2)}$ can be of the following six
kinds, depending on whether the occurrences of variables in it are
universally bound, free or existentially bound:
\begin{center}
\begin{tabular}{c|c|c}
$\alpha^1$ & $\alpha^2$ & name of kind
\\
\hline universally bound & universally bound & $(\forall,\forall)$
\\
universally bound & free & $(\forall,\emptyset)$
\\
universally bound & existentially bound & $(\forall,\exists)$
\\
free & free & $(\emptyset,\emptyset)$
\\
free & existentially bound & $(\emptyset,\exists)$
\\
existentially bound & existentially bound & $(\exists,\exists)$
\end{tabular}
\end{center}
The kinds not mentioned---namely, $(\emptyset,\forall)$,
$(\exists,\forall)$ and $(\exists,\emptyset)$---are not possible.

If a $(\emptyset,\emptyset)$ couple occurs in a cluster, then this
cluster is a singleton. If a $(\forall,\forall)$ couple or a
$(\exists,\exists)$ couple occurs in a cluster, then this cluster
has an eigengate. Together with $(\forall,\forall)$ couples in a
cluster we can find only $(\forall,\forall)$ couples and
$(\forall,\exists)$ couples, and analogously together with
$(\exists,\exists)$ couples in a cluster we can find only
$(\exists,\exists)$ couples and $(\forall,\exists)$ couples. This
is a consequence of the Eigengate Remark and of the fact that a
cluster can have only one eigengate. Couples of the
$(\forall,\exists)$, $(\forall,\emptyset)$ and
$(\emptyset,\exists)$ kind, can be joined together in a cluster
without eigengate. In the bottom cluster in our example above, we
have a ${(\forall,\exists)}$ couple $(x,z)$ and a
${(\emptyset,\exists)}$ couple $(y,z)$. This cluster has no
eigengate.

A Gentzen term $f$ is \emph{eigendiversified} when it is
variable-pure, cut-free and renaming-free, and, moreover, for
every arc-cluster of $f$ that has an eigengate the free variable
of this arc-cluster is different from the free variable of any
other arc-cluster of $f$. (Eigendiversification is inspired by
\cite{G35}, Section III.3.10.) We have the following for \GQDS.

{\samepage \prop{Eigendiversification Lemma}{For every
variable-pure, cut-free and renaming-free Gentzen term $f$ there
is an eigendiversified Gentzen term $f'$ such that ${f=f'}$.}}

\dkz This lemma is proved by replacing the free variable of an
arc-cluster of $f$ that has an eigengate by a variable new for
$f$. By doing that for every arc-cluster of $f$ that has an
eigengate we obtain $f'$, which differs from $f$ just in the
indices of identity arrows and of Gentzen operations. The
equations $\mbox{($Q^S$ {\it ren})}$ of Section 1.7 guarantee that
${f'=f}$.

Take for example a subterm of $f$ of the form
${\forall^R_{x,X}g\!:Y\vdash \forall_x X\vee Z}$ for $g\!:Y\vdash
X^x_u\vee Z$, and suppose that in $f'$ we have instead at the same
place a subterm $\forall^R_{x,X}g'\!:Y\vdash \forall_x X\vee Z$
for ${g'\!:Y\vdash X^x_{u'}\vee Z}$ where $u'$ is new for $f$.
Since $u$ is not free in the type of ${\forall^R_{x,X}g}$, so $u$
is not free in the type of $g'$, and by $\mbox{($\forall^R$ {\it
ren})}$ we have that
${\forall^R_{x,X}g'=\forall^R_{x,X}[g']^{u'}_u}$. By the
strengthened version of Renaming Elimination (see Section 1.8), we
obtain that ${[g']^{u'}_u=g}$. \qed

\vspace{2ex}

Next we give the following inductive definition of the notion of
\emph{subform} of a form set, which extends the notion of subform
set (see the beginning of Section 1.5):
\begin{itemize}
 \item[] $X$ is a subform of $X$;\vspace{-1ex}
 \item[] if $Y$ is a subform of $X^x_y$,
 then $Y$ is a subform of $Q_xX$;\vspace{-1ex}
 \item[] if $X$ and $Y$ are subforms of $X'$ and $Y'$ respectively,
 then ${X\ks Y}$ is a subform of ${X'\ks Y'}$;\vspace{-1ex}
 \item[] if $X$ is a subform of $Y$, then $X$ is a subform of ${Y\ks Z}$.
\end{itemize}
(Note that ${Y\ks Z}$ is the same form set as ${Z\ks Y}$.)

We say that a Gentzen term is $\forall_xX$-\emph{regular} when it
does not have subterms of one of the following two forms:
\begin{itemize}
\item[(\emph{a})] ${\forall^R_{z,Z}\forall^L_{x,X}g}$ for $g$ of
the type ${X^x_y\kon U\vdash Z^z_y\vee V}$, or of one of the three
types obtained by omitting ${\kon\;U}$ or
${\vee\:V}$,\vspace{-.5ex} \item[(\emph{b})]
${\exists^L_{z,Z}\forall^L_{x,X}g}$ for $g$ of the type
${X^x_y\kon Z^z_y\kon Y\vdash V}$ or ${X^x_y\kon Z^z_y\vdash V}$.
\end{itemize}

We can then prove the following lemma.

\prop{Lemma $\forall^L$}{If ${f_1\!:\forall_xX\kon U_1\vdash Z_1}$
is an eigendiversified Gentzen term and ${f_2\!:X^x_u\kon
U_2\vdash Z_2}$ is a variable-pure Gentzen term such that $U_1$
and $Z_1$ are subforms of $U_2$ and $Z_2$ respectively, then $f_1$
is $\forall_xX$-regular. The same holds if in all the types above
we omit ${\kon\;U_1}$ and ${\kon\;U_2}$, or just ${\kon\;U_1}$.}

\dkz Suppose $f_1$ is not $\forall_xX$-regular. We will consider
only the case when $f_1$ has a subterm of the form
${\forall^R_{z,Z}\forall^L_{x,X}g}$ for ${g\!:X^x_y\kon U\vdash
Z^z_y\vee V}$. When $f_1$ has a subterm of the form mentioned in
the remaining cases of (\emph{a}) or in case (\emph{b}), we
proceed analogously. Let $Z'$ be the subform set of the target
$Z_1$ of $f_1$ containing exactly the same predicate letters as
$Z$, and let $\gamma$ be the occurrence of $\forall_x$ at the
beginning of $\forall_xX$ in the source of $f_1$. By the
assumption that $f_1$ is eigendiversified, for $\alpha$ an
occurrence of $x$ in $X$ and $\beta$ an occurrence of $z$ in $Z'$,
either
\begin{itemize}
\item[(1)] we have a couple ${(\alpha,\beta)}$ of $f_1$,
or\vspace{-1ex} \item[(2)] we have two clustered couples
${(\alpha,\alpha')}$ and ${(\beta',\beta)}$ of $f_1$ with a
sequence of bridges ${\gamma_1,\ldots,\gamma_{n-1}}$ different
from $\gamma$ that ensure their clustering.
\end{itemize}

If we have clustered couples as in (2), but $\gamma$ occurs in
${\gamma_1,\ldots,\gamma_{n-1}}$, then let $\gamma_j$, for ${1\leq
j\leq n\mn 1}$, be the rightmost occurrence of $\gamma$ in
${\gamma_1,\ldots,\gamma_{n-1}}$. The bridge $\gamma_j$ is between
${(\alpha_j,\beta_j)}$ and ${(\alpha_{j+1},\beta_{j+1})}$, and
$\alpha_{j+1}$ is an occurrence of $x$ in $X$ in the source of
$f_1$. If ${j=n\mn 1}$, then we have (1), and if ${j< n\mn 1}$,
then we have (2).

By Renaming Elimination and the Cut-Elimination Theorem (see
Sections 1.8-9), we may assume that $f_2$ is cut-free and
renaming-free. If we have (1), then we should have a
$(\emptyset,\forall)$ couple of $f_2$, which is impossible. If we
have (2), then in ${(\alpha,\alpha')}$ we have that $\alpha$ is
universally bound. Since $U_1$ and $Z_1$ are subforms of $U_2$ and
$Z_2$ respectively, there should be a bridgeable chain of couples
of $f_2$ whose first member is of the kind $(\emptyset,\emptyset)$
or $(\emptyset,\exists)$, and whose last member is of the kind
$(\forall,\forall)$. The bridges ensuring the bridgeability of
this chain of couples of length $n$ correspond to
${\gamma_1,\ldots,\gamma_{n-1}}$. However, a bridgeable chain of
couples of the kind above cannot exist, as we said after the
Eigengate Remark. (All couples in a bridgeable chain of couples
belong to the same cluster.) \qed

\vspace{2ex}

There is an analogous lemma that should be called \emph{Lemma}
$\exists^R$. It involves $\exists^R$ instead of $\forall^L$ (which
engenders the notion of $\exists_xX$-regularity). We can now
finally state the following lemmata.

\prop{Invertibility Lemma for $\forall^L$}{If
${f_1\!:\forall_xX\kon U_1\vdash Z_1}$ is an eigendiversified
Gentzen term, and there is a variable-pure Gentzen term
$f_2\!:X^x_y\kon U_2\vdash Z_2$ where $U_1$ and $Z_1$ are subforms
of $U_2$ and $Z_2$ respectively, then there is an eigendiversified
Gentzen term ${f_1'\!:X^x_y\kon U_1\vdash Z_1}$ such that
${\forall^L_{x,X}f_1'=f_1}$. The same holds if in all the types
above we omit ${\kon\;U_1}$ and ${\kon\;U_2}$, or just
${\kon\;U_1}$.}

\vspace{-2ex}

\prop{Invertibility Lemma for $\exists^R$}{If
${f_1\!:Z_1\vdash\exists_xX\vee U_1}$ is an eigendiversified
Gentzen term, and there is a variable-pure Gentzen term
$f_2\!:Z_2\vdash X^x_y\vee U_2$ where $U_1$ and $Z_1$ are subforms
of $U_2$ and $Z_2$ respectively, then there is an eigendiversified
Gentzen term ${f_1'\!:Z_1\vdash X^x_y\vee U_1}$ such that
${\exists^R_{x,X}f_1'=f_1}$. The same holds if in all the types
above we omit ${\vee\;U_1}$ and ${\vee\;U_2}$, or just
${\vee\;U_1}$.}

These two lemmata are proved by induction on the complexity of
$f_1$ with the help of the equations $({\!\ks\!}Q^S)$, $(Q^SQ^S)$,
$(\exists^R\forall^L)$, $(Q^RQ^L)$ and $(\exists^S\forall^S)$,
from the beginning of the section, and Lemmata $\forall^L$ and
$\exists^R$. Without the Lemmata $\forall^L$ and $\exists^R$ we
would not be able to apply the equations $(Q^RQ^L)$ and
$(\exists^S\forall^S)$. If we end up with ${f_1'}$ of the type
${X^x_{y'}\kon U_1\vdash Z_1}$, or ${Z_1\vdash X^x_{y'}\vee U_1}$
respectively, where $y'$ is different from $y$, then we proceed in
the spirit of the proof of the Eigendiversification Lemma by using
the equations \mbox{($Q^T$ {\it ren})} of Section 1.7.

\subsection{Proof of \QDS\ Coherence}
We are now ready to prove that \GQDS\ is a preorder (see the
beginning of Section 1.5). Suppose we have two arrow terms of
\GQDS\ of the same type. These arrow terms are equal to two
Gentzen terms $f_1$ and $f_2$ by the new Gentzenization Lemma of
Section 1.5. Let ${x_1,\ldots,x_k}$ be the variables occurring
bound in the types of $f_1$ and $f_2$. By the
Variable-Purification Lemma of Section 1.6, we have that
\[
\quad\quad\quad\quad\quad\quad f_i=h^i_2\cirk f_i'\cirk
h^i_1,\quad {\mbox {\rm for }} i\in\{1,2\},
\]
where $f_i'$ is variable-pure, while $h^i_1$ and $h^i_2$ are
isomorphisms of \GQDS. By choosing the same new variables
${x_1',\ldots,x_k'}$ both for $f_1$ and $f_2$, we obtain that
${h^1_j=h^2_j}$ for ${j\in\{1,2\}}$. So if ${f_1'=f_2'}$, we will
be able to derive ${f_1=f_2}$.

The Gentzen terms $f_1'$ and $f_2'$ are variable-pure, and hence
by Renaming Elimination and the Cut-Elimination Theorem (see
Sections 1.8-9) we can assume that they are cut-free and
renaming-free. By the Eigendiversification Lemma of the preceding
section, we can assume that they are moreover eigendiversified.

Let the \emph{quantity} of a Gentzen term of \GQDS\ be the sum of
the number of predicate letters in its source (which is equal to
the number of predicate letters in its target) with the number of
occurrences of quantifier prefixes in its source and target. Then
we proceed by induction on the quantity of $f_2'$, which is equal
to the quantity of $f_1'$, in order to show that ${f_1'=f_2'}$. In
the basis of this induction, if ${n=1}$, then ${f_1'=f_2'=\mj_A}$,
where $A$ is atomic. In the induction step we apply the
invertibility lemmata of the preceding section (cf.\ \cite{DP04},
end of Section 7.9).

So \GQDS\ is a preorder. And, as we explained at the beginning of
Section 1.5, we have proved thereby \QDS\ Coherence.

\section{Coherence of \QPNN}
\subsection{The categories \QPNN\ and \QPN}
In this section we introduce the category \QPNN\ (here \textbf{PN}
comes from ``proof net''), which corresponds to the multiplicative
fragment without propositional constants of classical linear
first-order predicate logic without mix. This category extends
with quantifiers the propositional category $\PN^\neg$ of
\cite{DP07} (Section 2.2).

The category \QPNN\ is defined as the category \QDS\ in Section
1.2 save that we make the following additions and changes. Instead
of the language \eL\ of Section 1.1 we have the language \eLn,
which differs from \eL\ by having the additional unary connective
$\neg$. So in the definition of formula we have the additional
clause
\begin{itemize}
\item[] if $A$ is a formula, then $\neg A$ is a formula.
\end{itemize}
The objects of the category \QPNN\ are the formulae of \eLn.

Let $\mX_n$ stand for the sequence ${x_1,\ldots,x_n}$ when ${n\geq
1}$, and for the empty sequence when ${n=0}$. For $A$ a formula,
let $A^{\mX_n}_{\mY_n}$ stand for
${A^{x_1\ldots\,x_n}_{y_1\ldots\,y_n}}$ when $n\geq 1$, and for
$A$ when $n=0$; so $A^{\mX_n}_{\mY_n}$ is the result of a series
of $n$ substitutions. We use $Q_{\mX_n}$ as an abbreviation for
${Q_{x_n}\!\ldots Q_{x_1}}$ when $n\geq 1$, and for the empty
sequence when $n=0$.

When $A$ is a formula containing free exactly the mutually
different variables $\mX_n$ in order of first occurrence counting
from the left, we say that $\mX_n$ is the \emph{free-variable
sequence of} $A$. For example, the free-variable sequence of
$\forall_y(Pyx\kon \exists_xRzxz)$ is $x,z$ (provided $x$, $y$ and
$z$ are all mutually different).

To define the arrow terms of \QPNN, in the inductive definition we
had for the arrow terms of \QDS\ we replace \eL\ by \eLn\ and
assume in addition that for all formulae $A$ and $B$ of \eLn, and
for $\mX_n$ being the free-variable sequence of $B$, the following
\emph{primitive arrow terms}:
\begin{tabbing}
\hspace{10em}\=$\Delta^{\forall}_{B,A}$\=$:A\vdash A\kon
\forall_{\mX_n}\!(\neg B\vee B)$,\\*[1ex]
\>$\Sigma^{\exists}_{B,A}$\>$:\exists_{\mX_n}\!(B\kon\neg B)\vee
A\vdash A$
\end{tabbing}
are arrow terms of \QPNN. In other words, ${\forall_{\mX_n}\!(\neg
B\vee B)}$ is the universal closure of ${\neg B\vee B}$, and
${\exists_{\mX_n}\!(B\kon\neg B)}$ is the existential closure of
${B\kon\neg B}$. We assume throughout the remaining text that
${\Xi\in\{\Delta,\Sigma\}}$.

We call the first index $B$ of $\Delta^{\forall}_{B,A}$ and
$\Sigma^{\exists}_{B,A}$ the {\it crown} index, and the second
index $A$ the {\it stem} index. (We need the stem $A$ because we
lack propositional constants.) The right conjunct
${\forall_{\mX_n}\!(\neg B\vee B)}$ in the target of
${\Delta^{\forall}_{B,A}}$ is the {\it crown} of
${\Delta^{\forall}_{B,A}}$, and the left disjunct
${\exists_{\mX_n}\!(B\kon\neg B)}$ in the source of
${\Sigma^{\exists}_{B,A}}$ is the {\it crown} of
$\Sigma^{\exists}_{B,A}$. We have analogous definitions of crown
and stem indices, and crowns, for $\Sigma^{\forall}$,
$\Delta'^{\forall}$, $\Sigma'^{\forall}$, $\Delta^{\exists}$,
$\Sigma'^{\exists}$ and $\Delta'^{\exists}$, which will be
introduced later. (The symbol $\Delta$ should be associated with
the Latin \emph{dexter}, because in ${\Delta^{\forall}_{B,A}}$,
${\Delta'^{\forall}_{B,A}}$, ${\Delta^{\exists}_{B,A}}$ and
${\Delta'^{\exists}_{B,A}}$ the crown is on the right-hand side of
the stem; analogously, $\Sigma$ should be associated with
\emph{sinister}.)

Before we define the arrows of \QPNN, we introduce a number of
abbreviations:
\begin{tabbing}
\hspace{1em}\=for $n=0$ and ${\alpha\in\{\iota,\gamma\}}$,
\hspace{2em}\=$\alpha^{Q_{\mX_n}}_B=_{df}\mj_B$,\\[1ex]
\>for $n>0$,\\*[.5ex] \>for
$(\alpha,Q)\in\{(\iota,\forall),(\gamma,\exists)\}$, \>
$\alpha^{Q_{\mX_n}}_B=_{df}\alpha^{Q_{\mX_{n\mn
1}}}_B\!\cirk\alpha^{Q_{x_n}}_{Q_{\mX_{n\mn
1}}B}$\=$:Q_{\mX_n}B\vdash B$,\\*[1ex] \>for
$(\alpha,Q)\in\{(\iota,\exists),(\gamma,\forall)\}$, \>
$\alpha^{Q_{\mX_n}}_B=_{df}\alpha^{Q_{x_n}}_{Q_{\mX_{n\mn
1}}B}\cirk\alpha^{Q_{\mX_{n\mn 1}}}_B$\>$:B\vdash
Q_{\mX_n}B$,\\[2ex]
\>for $n=0$,\>$[f]^{\mX_n}_{\mY_n}=_{df}f$,\\*[1ex] \>for
$n>0$,\>$[f]^{\mX_n}_{\mY_n}=_{df}
[[f]^{\mX_{n\mn 1}}_{\mY_{n\mn 1}}]^{x_n}_{y_n}$,\\[2ex]
\>for $\mU_n$ and $\mV_n$ not free in $A$,\\*[1ex]
\hspace{7em}\=$\tau^{\forall_{\mX_n}}_{A,\mU_n,\mV_n}\!=_{df}\forall_{\mV_n}
[\iota^{\forall_{\mU_n}}_{A^{\mX_n}_{\mU_n}}]^{\mU_n}_{\mV_n}\cirk
\gamma^{\forall_{\mV_n}}_{\forall_{\mU_n}A^{\mX_n}_{\mU_n}}$\=$:
\forall_{\mU_n}A^{\mX_n}_{\mU_n}\,$\=$\vdash\forall_{\mV_n}A^{\mX_n}_{\mV_n}$,\\*[1ex]
\>$\tau^{\exists_{\mX_n}}_{A,\mV_n,\mU_n}\!=_{df}
\gamma^{\exists_{\mV_n}}_{\exists_{\mU_n}A^{\mX_n}_{\mU_n}}
\!\cirk\exists_{\mV_n}
[\iota^{\exists_{\mU_n}}_{A^{\mX_n}_{\mU_n}}]^{\mU_n}_{\mV_n}
$\>$:\exists_{\mV_n}A^{\mX_n}_{\mV_n}$\>$\vdash\exists_{\mU_n}A^{\mX_n}_{\mU_n}$.
\end{tabbing}
In \QDS, for $n=0$ we have
${\tau^{Q_{\mX_n}}_{A,\mU_n,\mV_n}\!\!=\mj_A}$, and for $n>0$ we
have
\begin{tabbing}
\hspace{7em}\=$\tau^{\forall_{\mX_n}}_{A,\mU_n,\mV_n}\!=
\tau^{\forall\!_{x_n}}_{\forall_{\mV_{n\mn 1}}A,u_n,v_n} \!\cirk
\forall_{u_n}\tau^{\forall_{\mX_{n\mn 1}}}_{A,\mU_{n\mn
1},\mV_{n\mn 1}}$,\\*[1ex]
\>$\tau^{\exists_{\mX_n}}_{A,\mV_n,\mU_n}\!=
\exists_{u_n}\tau^{\exists_{\mX_{n\mn 1}}}_{A,\mV_{n\mn
1},\mU_{n\mn 1}}\!\cirk \tau^{\exists_{x_n}}_{\exists_{\mV_{n\mn
1}}A,v_n,u_n}\!$.
\end{tabbing}

For $\mX_n$ being the free-variable sequence of $B$, we have also
the abbreviations
\begin{tabbing}
\hspace{7em}\=$\Sigma^{\forall}_{B,A}\,$\=$=_{df}\hat{c}_{A,\forall_{\mX_n}\!(\neg
B\vee B)}\cirk\Delta^{\forall}_{B,A}\;$\=$:A\vdash
\forall_{\mX_n}\!(\neg B\vee B)\kon A$,\kill

\>$\Sigma^{\forall}_{B,A}$\>$=_{df}\hat{c}_{A,\forall_{\mX_n}\!(\neg
B\vee B)}\cirk\Delta^{\forall}_{B,A}$\>$:A\vdash
\forall_{\mX_n}\!(\neg B\vee B)\kon A$,\\*[1ex]
\>$\Delta^{\exists}_{B,A}$\>$=_{df}\Sigma^{\exists}_{B,A}
\cirk\check{c}_{\exists_{\mX_n}\!(B\kon\neg
B),A}$\>$:A\vee\exists_{\mX_n}\!(B\kon\neg B)\vdash A$,\\[2ex]
\>$\Dk_{B,A}$\>$=_{df}(\mj_A\kon\iota^{\forall_{\mX_n}}_{\neg
B\vee B})\cirk\Delta^{\forall}_{B,A}$\>$:A\vdash A\kon(\neg B\vee
B)$,\\*[1ex]
\>$\Sd_{B,A}$\>$=_{df}\Sigma^{\exists}_{B,A}\cirk(\iota^{\exists_{\mX_n}}_{B\kon\neg
B}\vee\mj_A)$\>$:(B\kon\neg B)\vee A\vdash A$,\\[1.5ex]
\>$\Dkp_{B,A}$\>$=_{df}(\mj_A\kon \check{c}_{B,\neg
B})\cirk\Dk_{B,A}$\>$:A\vdash A\kon(B\vee\neg B)$,
\\*[1ex]
\>$\Sdp_{B,A}$\>$=_{df}\Sd_{B,A}\cirk(\hat{c}_{\neg
B,B}\vee\mj_A)$\>$:(\neg B\kon B)\vee A\vdash A$.
\end{tabbing}

To define the arrows of \QPNN\ we assume in the inductive
definition we had for the equations of \QDS\ the following
additional axiomatic equations:
\begin{tabbing}

\mbox{\hspace{2em}}\=$\mbox{({\it ren}
$\Xi^Q$)}$\quad\=$[\Xi^Q_{B,A}]^x_y=\Xi^Q_{B,A^x_y}$\kill

\>\mbox{($\Delta^{\forall}$ {\it
nat})}\>$(f\kon\mj_{\forall_{\mX_n}\!(\neg B\vee
B)})\cirk\Delta^{\forall}_{B,A}=\Delta^{\forall}_{B,D}\cirk f$,
\\*[1ex]
\>\mbox{($\Sigma^{\exists}$ {\it nat})}\> $f\cirk
\Sigma^{\exists}_{B,A}=\Sigma^{\exists}_{B,D}\cirk(\mj_{\exists_{\mX_n}\!(B\kon
\neg B)}\vee f)$,
\\[2ex]
\>\mbox{($\hat{b}\Delta^{\forall}$)}\>
$\hat{b}^{\rts}_{A,B,\forall_{\mX_n}\!(\neg C\vee
C)}\cirk\Delta^{\forall}_{C,A\kon
B}\,$\=$=\mj_A\kon\Delta^{\forall}_{C,B}$,
\\*[1ex]
\>\mbox{($\check{b}\Sigma^{\exists}$)}\>
$\Sigma^{\exists}_{C,B\vee
A}\cirk\check{b}^{\rts}_{\exists_{\mX_n}(C\kon\neg
C),B,A}$\>$=\Sigma^{\exists}_{C,B}\vee\mj_A$,
\\[2ex]
\>\mbox{($d\Sigma^{\forall}$)}\> $d_{\forall_{\mX_n}\!(\neg A\vee
A),B,C}\cirk\Sigma^{\forall}_{A,B\vee
C}$\>$=\Sigma^{\forall}_{A,B}\vee\mj_C$,
\\*[1ex]
\>\mbox{($d\Delta^{\exists}$)}\>$\Delta^{\exists}_{A,C\kon B}\cirk
d_{C,B,\exists_{\mX_n}\!(A\kon\neg
A)}$\>$=\mj_C\kon\Delta^{\exists}_{A,B}$,
\\[2ex]
\>\mbox{($\Sd\Dk$)}\>$\Sd_{A,A}\cirk d_{A,\neg A,A}
\cirk\Dk_{A,A}=\mj_A$,
\\*[1ex]
\>\mbox{($\Sdp\Dkp$)}\>$\Sdp_{A,\neg A}\cirk d_{\neg A,A,\neg
A}\cirk\Dkp_{A,\neg A}=\mj_{\neg A}$,\\[2.5ex]
\>\mbox{({\it ren}
$\Xi^Q$)}\>$[\Xi^Q_{B,A}]^x_y=\Xi^Q_{B,A^x_y}$,\quad for
$\Xi^Q\in\{\Delta^{\forall},\Sigma^\exists\}$,\\[2.5ex]
\>\mbox{($\Delta\tau$)}\>$\Delta^{\forall}_{B^{\mX_n}_{\mV_n},A}\,$\=$=(\mj_A\kon
\tau^{\forall_{\mX_n}}_{\neg B\vee
B,\mU_n,\mV_n})\cirk\Delta^{\forall}_{B^{\mX_n}_{\mU_n},A}$,\\*[1ex]
\>\mbox{($\Sigma\tau$)}\>$\Sigma^{\exists}_{B^{\mX_n}_{\mV_n},A}$\>$=\Sigma^{\exists}_{B^{\mX_n}_{\mU_n},A}
\cirk(\tau^{\exists_{\mX_n}}_{B\kon\neg
B,\mV_n,\mU_n}\!\!\vee\mj_A)$.
\end{tabbing}
The equation \mbox{({\it ren} $\alpha$)} of Section 1.2 does not
hold when $\alpha$ is $\Delta^{\forall}$ or $\Sigma^{\exists}$,
but instead we have the equations \mbox{({\it ren} $\Xi^Q$)}
above. This defines the category \QPNN.

In this list of axiomatic equations the equations
\mbox{($\Sd\Dk$)} and \mbox{($\Sdp\Dkp$)} are taken as they stand
from \cite{DP07} (Section 2.2), where they were used to axiomatize
the category $\PN^{\neg}$. The preceding first six axiomatic
equations of $\PN^{\neg}$ are obtained from the first six
axiomatic equations of \QPNN\ above by replacing
$\Delta^{\forall}$ and $\Sigma^{\exists}$ with $\Dk$ and $\Sd$
respectively, and by deleting quantifier prefixes. It is clear
that we can derive these axiomatic equations of $\PN^{\neg}$ in
\QPNN, and hence we have in \QPNN\ all the equations of
$\PN^{\neg}$, with $A$, $B$, $C,\ldots$ being formulae of the
language \eLn. The really new axiomatic equations of \QPNN\ are
only the last displayed \mbox{({\it ren} $\Xi^Q$)} and
\mbox{($\Xi\tau$)}.

We have in \QPNN\ the additional abbreviations
\begin{tabbing}
\hspace{5em}\=$\Delta'^{\forall}_{B,A}$\=$=_{df}(\mj_A\kon\forall_{\mX_n}\check{c}_{B,\neg
B})\cirk\Delta^{\forall}_{B,A}$\=$:A\vdash
A\kon\forall_{\mX_n}\!(B\vee\neg B)$,\\*[1ex]
\>$\Sigma'^{\exists}_{B,A}$\>$=_{df}\Sigma^{\exists}_{B,A}\cirk(\exists_{\mX_n}\hat{c}_{\neg
B,B}\vee\mj_A)$\>$:\exists_{\mX_n}\!(\neg B\kon B)\vee A\vdash
A$,\\[1.5ex]
\>$\Sigma'^{\forall}_{B,A}\,$\>$=_{df}\hat{c}_{A,\forall_{\mX_n}\!(B\vee\neg
B)}\cirk\Delta'^{\forall}_{B,A}$\=$
:A\vdash\forall_{\mX_n}\!(B\vee\neg B)\kon A$,
\\*[1ex]
\>$\Delta'^{\exists}_{B,A}$\>$=_{df}\Sigma'^{\exists}_{B,A}\cirk\check{c}_{\exists_{\mX_n}\!(\neg
B\kon B),A}$\>$:A\vee\exists_{\mX_n}\!(\neg B\kon B)\vdash A$,
\end{tabbing}
and as in \cite{DP07} (see Section 2.2) we have also the
abbreviations
\begin{tabbing}
\hspace{5em}\=$\Delta'^{\forall}_{B,A}$\=$=_{df}(\mj_A\kon\forall_{\mX_n}\check{c}_{B,\neg
B})\cirk\Delta^{\forall}_{B,A}$\=$:A\vdash
A\kon\forall_{\mX_n}\!(B\vee\neg B)$,\kill

\>$\Sk_{B,A}$\>$=_{df}\hat{c}_{A,\neg B\vee B}\cirk\Dk_{B,A}$\=$:
A\vdash (\neg B\vee B)\kon A$,\\*[1ex]
\>$\Dd_{B,A}$\>$=_{df}\Sd_{B,A}\cirk\check{c}_{B\kon\neg
B,A}$\>$:A\vee(B\kon\neg B)\vdash A$,\\[1.5ex]
\>$\Skp_{B,A}\,$\>$=_{df}\hat{c}_{A,B\vee\neg
B}\cirk\Dkp_{B,A}$\>$:A\vdash (B\vee\neg B)\kon A$,
\\*[1ex]
\>$\Ddp_{B,A}$\>$=_{df}\Sdp_{B,A}\cirk\check{c}_{\neg B\kon
B,A}$\>$:A\vee(\neg B\kon B)\vdash A$.
\end{tabbing}

Note that the equations \mbox{($\Xi\tau$)} say that we could
define our arrows $\Delta^{\forall}_{B,A}$ and
$\Sigma^{\exists}_{B,A}$ in terms of such arrows with the proviso
that the set of variables free in the stem index $A$ and the set
variables free in the crown index $B$ are disjoint.

In \QDS\ and \QPNN\ for ${f\!:B\vdash\forall_{\mX_n}A}$ and
${g\!:\exists_{\mX_n}A\vdash B}$ such that the variables $\mX_n$
are not free in $B$ we have the equations
\begin{tabbing}
\mbox{\hspace{6em}}\=\mbox{({\it ren}
$\Xi^Q$)}\quad\=$[\Xi^Q_{B,A}]^x_y=\Xi^Q_{B,A^x_y}$\kill

\>\>$\forall_{\mX_n}\!(\iota^{\forall_{\mX_n}}_A\!\cirk
f)\cirk\gamma^{\forall_{\mX_n}}_B$\=$=f$,\\*[1ex]
\>\>$\gamma^{\exists_{\mX_n}}_B\!\cirk\exists_{\mX_n}\!(g\cirk\iota^{\exists_{\mX_n}}_A\!)$\>$=g$,
\end{tabbing}
which generalize the equations \mbox{($\forall$ {\it ext})} and
\mbox{($\exists$ {\it ext})} of Section 1.2. These equations,
together with equations of \QDS\ analogous to the equations
\mbox{$(Q\tau\iota)$} of Section 1.3 and the isomorphism of
$\tau^{Q_{\mU_n}}\!$, entail the following cancellation
implications:
\begin{tabbing}
\mbox{\hspace{2em}}\=\mbox{({\it ren}
$\Xi^Q$)}\quad\=$[\Xi^Q_{B,A}]^x_y=\Xi^Q_{B,A^x_y}$\kill

\>\mbox{($\forall\iota$ \emph{canc})}\>if
$[\iota^{\forall_{\mX_n}}_A]^{\mX_n}_{\mY_n}\!\cirk
f_1\,$\=$=[\iota^{\forall_{\mX_n}}_A]^{\mX_n}_{\mY_n}\!\cirk f_2$,
\= then $f_1\,$\=$=f_2$,\\*[1ex] \>\mbox{($\exists\iota$
\emph{canc})}\>if
$g_1\cirk[\iota^{\exists_{\mX_n}}_A]^{\mX_n}_{\mY_n}$\>$=g_2\cirk
[\iota^{\exists_{\mX_n}}_A]^{\mX_n}_{\mY_n}\!$,\> then
$g_1$\>$=g_2$,
\end{tabbing}
provided the variables $\mY_n$ are not free in the source of $f_1$
and $f_2$ and in the target of $g_1$ and $g_2$.

In \QPNN\ we have \emph{stem-increasing} equations analogous to
the stem-increasing equations of \cite{DP07} (Section 2.5;
$\hat{\Delta}$ and $\check{\Sigma}$ are replaced by
$\Delta^{\forall}$ and $\Sigma^{\exists}$ respectively, which
entails further adjustments). The equations
\mbox{($\hat{b}\Delta^{\forall}$)},
\mbox{($\check{b}\Sigma^{\exists}$)}, \mbox{($d\Sigma^{\forall}$)}
and \mbox{($d\Delta^{\exists}$)} are such stem-increasing
equations (when read from right to left), and there are further
such equations for all the arrows $\Xi^Q_{B,A}$ and
${\Xi'}^Q_{B,A}$.

We have in \QPNN\ the following additional stem-increasing
equations:
\begin{tabbing}
\hspace{1em}\=\mbox{({\it ren}
$\Xi^Q$)}\=$[\Xi^Q_{B,A}]^x_y=\Xi^Q_{B,A^x_y}$\kill

\>\mbox{($\forall\Delta^{\forall}$)}\>$\forall_x\Delta^{\forall}_{B,A}\,$\=
$=\forall_x(\iota^{\forall\!_x}_A\kon\mj_{\forall_{\mX_n}\!(\neg
B\vee
B)})\cirk\gamma^{\forall\!_x}_{\forall\!_xA\kon\forall_{\mX_n}\!(\neg
B\vee B)}\cirk\Delta^{\forall}_{B,\forall\!_xA}$,\\[1ex]
\>\mbox{($\exists\Delta^{\forall}$)}\>$\exists_x\Delta^{\forall}_{B,A}\,$\>
$=\hat{\theta}^{\exists_x\rts}_{A,\forall_{\mX_n}\!(\neg
B\vee B)}\cirk\Delta^{\forall}_{B,\exists_xA}$,\\[1ex]
\>\mbox{($\forall\Sigma^{\exists}$)}\>$\forall_x\Sigma^{\exists}_{B,A}\,$\>
$=\Sigma^{\exists}_{B,\forall\!_xA}\cirk\check{c}_{\exists_{\mX_n}\!(B\kon\neg
B),\forall\!_xA}\cirk\check{\theta}^{\forall\!_x\str}_{A,\exists_{\mX_n}\!(B\kon\neg
B)}\cirk\forall_x\check{c}_{A,\exists_{\mX_n}\!(B\kon\neg
B)}$,\\[1ex]
\>\mbox{($\exists\Sigma^{\exists}$)}\>$\exists_x\Sigma^{\exists}_{B,A}\,$\>
$=\Sigma^{\exists}_{B,\exists_xA}\cirk\gamma^{\exists_x}_{\exists_{\mX_n}\!(B\kon\neg
B)\vee\exists_xA}\cirk\exists_x(\mj_{\exists_{\mX_n}\!(B\kon\neg
B)}\vee\iota^{\exists_x}_A)$,
\end{tabbing}
which are derived with the help of the implications
\mbox{($Q\iota$ \emph{canc})}, \QDS\ equations, \QDS\ Coherence
and the naturality of $\Delta^{\forall}$ and $\Sigma^{\exists}$ in
their stem index.

We introduce next a category called \QPN, for which we will
establish in Section 2.6 that it is equivalent to the category
\QPNN. The category \QPN\ is for us an auxiliary category (though
it is closer to the formulation of linear logic in \cite{G87}). We
prove coherence for this category in Section 2.5, and from that
and the equivalence of \QPNN\ and \QPN\ we infer coherence for
\QPNN\ in Section 2.7. The category \QPN\ is very much like \QPNN\
save that in its objects the negation connective $\neg$ is
prefixed only to atomic formulae. The arrow terms
$\Delta^{\forall}_{B,A}$ and $\Sigma^{\exists}_{B,A}$ are
primitive only for the crown index $B$ being an atomic formula.
Here is a more formal definition of \QPN.

For \Pe\ being the set of letters that we used to generate $\eL$
and $\eL_{\neg}$ in Sections 1.1 and 2.1, let $\Pe^{\neg}$ be the
set of predicate letters ${\{\neg P\mid P\in\Pe\}}$. The arity of
the new predicate letter $\neg P$ is the same as the arity of $P$.
The objects of \QPN\ are the formulae of the first-order language
${\eL^{\neg P}}$ generated from ${\Pe\cup\Pe^{\neg}}$ in the same
way as \eL\ was generated from \Pe\ in Section 1.1.

To define the arrow terms of \QPN, in the inductive definition we
had for the arrow terms of \QDS\ we replace \eL\ by ${\eL^{\neg
P}}$, and we assume in addition that for every formula $A$ of
${\eL^{\neg P}}$, for every predicate letter $P\in\Pe$ of arity
$n$, and for $\mX'_{n'}$ being the free-variable sequence of
$P\mX_n$,
\begin{tabbing}
\hspace{9em}\=$\Delta^{\forall}_{P\mX_n,A}$\=$:A\vdash
A\kon\forall_{\mX'_{n'}}\!(\neg P\mX_n\vee P\mX_n)$,
\\*[1ex]
\>$\Sigma^{\exists}_{P\mX_n,A}$\>$:\exists_{\mX'_{n'}}\!(P\mX_n\kon\neg
P\mX_n)\vee A\vdash A$
\end{tabbing}
are primitive arrow terms of \QPN.

To define the arrows of \QPN, we assume as additional axiomatic
equations in the inductive definition we had for the equations of
\QDS\ all the additional axiomatic equations assumed above for
\QPN, but restricted to the arrow terms
$\Delta^{\forall}_{P\mX_n,A}$ and $\Sigma^{\exists}_{P\mX_n,A}$
whose crown index is atomic. This defines the category \QPN.

\subsection{Development for \QDS, \QPNN\ and \QPN}
If $\beta$ is a primitive arrow term of \QPNN\ except $\mj_B$,
then we call $\beta$-{\it terms} of \QPNN\ the set of arrow terms
defined inductively as follows: $\beta$ is a $\beta$-term; if $f$
is a $\beta$-term, then for every $A$ in ${\eL_\neg}$ and all
variables $x$ and $y$ we have that ${\mj_A\ks f}$, ${f\ks \mj_A}$,
$Q_xf$ and $[f]^x_y$ are $\beta$-terms, provided $[f]^x_y$ is
defined.

In a $\beta$-term the subterm $\beta$ is called the {\it head} of
this $\beta$-term. For example, the head of the
$\Delta^{\forall}_{B,C}$-term
${\mj_A\kon\forall\!_x([\Delta^{\forall}_{B,C}]^y_z\vee\:\mj_E)}$
is $\Delta^{\forall}_{B,C}$.

We define $\mj$-\emph{terms} like $\beta$-terms; we just replace
$\beta$ in the definition above by $\mj_B$. So $\mj$-terms are
headless.

An arrow term of the form ${f_n \cirk \ldots \cirk f_1}$, where
$n\geq 1$, with parentheses tied to $\cirk$ associated
arbitrarily, such that for every ${i\in\{1,\ldots,n\}}$ we have
that $f_i$ is composition-free is called {\it factorized}. In a
factorized arrow term ${f_n\cirk\ldots\cirk f_1}$ the arrow terms
$f_i$ are called {\it factors}. A factor that is a $\beta$-term
for some $\beta$ is called a {\it headed} factor. A factorized
arrow term is called {\it headed} when each of its factors is
either headed or a $\mj$-term. A factorized arrow term
${f_n\cirk\ldots\cirk f_1}$ is called {\it developed} when $f_1$
is a $\mj$-term and if $n>1$, then every factor of
${f_n\cirk\ldots\cirk f_2}$ is headed. Analogous definitions of
$\beta$-term and developed arrow term can be given for \QDS.

We have the following lemma for \QDS.

\prop{Development Lemma}{For every arrow term $f$ there is a
developed arrow term $f'$ such that $f=f'$.}

\dkz This lemma would be easy to prove by using the categorial and
functorial equations together with the equation $\mbox{({\it ren}
$\cirk\!$)}$ of Section 1.2 if for $\mbox{({\it ren} $\cirk$)}$,
as for the other of these equations, we had that the right-hand
side is defined whenever the left-hand side is defined. Since this
need not be the case, we must first eliminate renaming. This can
be achieved by relying on the category \GQDS, cut elimination and
the result of Section 1.8.

We adapt the definitions of $\beta$-term and developed arrow term
to the category \GQDS\ of Section 1.5. For every arrow term $g$ of
\GQDS\ there is a Gentzen term $g'$ denoting the arrow $g$. By the
Variable-Purification Lemma of Section 1.6, we have that $g'$ is
equal to ${h_2\cirk g''\cirk h_2}$ where $h_1$ and $h_2$ are
compositions of $\tau$-terms and $g''$ is a variable-pure Gentzen
term, which by Renaming Elimination and the Cut-Elimination
Theorem (see Sections 1.8-9) we may assume to be cut-free and
renaming-free. Then by using the categorial and functorial
equations it is easy to obtain from ${h_2\cirk g''\cirk h_2}$ a
developed arrow term $g'''$ of \GQDS\ equal to the initial arrow
term $g$.

For an arbitrary arrow term ${f\!:A\vdash B}$ of \QDS\ we find a
diversified arrow term ${f'\!:A'\vdash B'}$ of \QDS\ such that $f$
is a letter-for-letter substitution instance of $f'$ (see the
beginning of Section 1.5). As in \cite{DP04} (Sections 3.2-3) we
pass by a functor $H_{\cal G}$ from $f$ to the arrow term $H_{\cal
G}f$ of \GQDS, which, as we have shown above, is equal to a
developed arrow term $(H_{\cal G}f)'''$ of \GQDS. By applying a
functor $H$ in the opposite direction we obtain a developed arrow
term $H((H_{\cal G}f)''')$ of \QDS, which we call $h$. The type of
$h$ is ${A''\vdash B''}$, where $A''$ and $B''$ belong to the same
form sets as $A'$ and $B'$ respectively. So by \QDS\ Coherence we
have that ${f'=j_2\cirk h\cirk j_1}$, where $j_1$ and $j_2$ are
headed factorized arrow terms of \QDS\ whose heads are of the
$\b{\xi}{}$ and $\c{\xi}$ kind. We obtain the arrow term of \QDS\
equal to $f$ as a letter-for-letter substitution instance of
${j_2\cirk h\cirk j_1}$. \hfill\mbox{\hspace{1em}}\qed

\vspace{1.7ex}

By relying on various renaming equations of \QDS, we can prove a
\emph{Refined Development Lemma} for \QDS, which differs from the
Development Lemma by requiring that in the developed arrow term
$f'$ renaming occurs only in subterms of the form
$[\iota^{Q\!_x}_A]^x_y$ for $x$ different from $y$ and free in
$A$. With the help of the Refined Development Lemma for \QDS, the
stem-increasing equations of \QPNN\ (see the preceding section),
together with the naturality of $\Delta^{\forall}$ and
$\Sigma^{\exists}$ in their stem index and the equations
\mbox{({\it ren} $\Xi^Q$)}, we can prove the Refined Development
Lemma, and hence also the Development Lemma, for the categories
\QPNN\ and \QPN\ too.

The Refined Development Lemma is not only important because of the
applications it will find latter in this paper. It is also
important because we can conclude from it that renaming, except in
$[\iota^{Q\!_x}_A]^x_y$ for $x$ different from $y$ and free in
$A$, is eliminable in \QDS, \QPNN\ and \QPN. This elimination of
renaming is not straightforward, but is achieved in a roundabout
way, involving cut elimination. The eliminability of renaming may
perhaps serve to explain why it is neglected as a primitive rule
of inference in logic.

\subsection{Some properties of \QDS}
In this section we establish some results concerning the category
\QDS\ of Section 1.2, which we will use to prove coherence for
\QPN\ and \QPNN. First we introduce a definition.

Suppose $X$ is the \mbox{$n$-th} occurrence of a predicate letter
(counting from the left) in a formula $A$ of \eL, and $Y$ is the
\mbox{$m$-th} occurrence of the same predicate letter in a formula
$B$ of \eL. Then we say that $X$ and $Y$ are {\it tied} in an
arrow ${f\!:A\vdash B}$ of \QDS\ when ${(n\mn 1,m\mn 1)\in Gf}$
(see Section 1.5; note that to find the \mbox{$n$-th} occurrence
we count starting from 1, but the ordinal ${n>0}$ is
${\{0,\ldots,n\mn 1\}}$). It is easy to establish that every
occurrence of a predicate letter in $A$ is tied to exactly one
occurrence of the same letter in $B$, and vice versa. This is
related to matters about diversification mentioned at the
beginning of Section 1.5.

For the lemma below, let $X$ in $A$ and $Y$ in $B$ be occurrences
of the same predicate letter tied in an arrow ${f\!:A\vdash B}$ of
\QDS, and let $S^A$ and $S^B$ be two finite (possibly empty)
sequences of quantifier prefixes. Then by an easy induction on the
complexity of $f$ we can prove the following, which generalizes
Lemma~2 of Section 2.4 of \cite{DP07}.

\prop{\mbox{$\kon\vee$ Lemma}}{It is impossible that $A$ has a
subformula ${S^AX\mX_n\!\kon A'}$ or ${A'\!\kon S^AX\mX_n}$ while
$B$ has a subformula ${S^BY\mY_n\!\vee B'}$ or ${B'\vee
S^BY\mY_n}$.}

For the next lemma, for ${i\in\{1,2\}}$ let $X_i$ in $A$ and $Y_i$
in $B$ be occurrences of the predicate letter $P_i$ tied in an
arrow ${f\!:A\vdash B}$ of \QDS\ (here $P_1$ and $P_2$ may also be
the same predicate letter).

\prop{\mbox{$\vee\kon$ Lemma}}{For every ${i,j\in\{1,2\}}$, it is
impossible that $A$ has a subformula ${X_i\mY_n\!\vee
X_{3-i}\mZ_m}$ while $B$ has a subformula ${Y_j\mU_k\!\kon
Y_{3-j}\mV_l}$.}

\noindent This lemma, exactly analogous to Lemma~3 of Section 2.4
of \cite{DP07}, is a corollary of lemmata exactly analogous to
Lemmata 3D and 3C of Section 2.4 of \cite{DP07}, which are easily
proved by induction on the complexity of the arrow term $f$.

As a matter of fact, the $\kon\vee$ and $\vee\kon$ Lemmata above
could be proved by supposing the contrary and deleting quantifiers
and individual variables together with arrow terms and operations
on arrow terms involving them, which would yield arrow terms
contradicting Lemma~2 and Lemma~3 respectively of Section 2.4 of
\cite{DP07}. The $\kon\vee$ Lemma is related to the acyclicity
condition of proof nets, while the $\vee\kon$ Lemma is related to
the connectedness condition (see \cite{DP07}, Sections 2.4, 7.1,
and references therein).

Next we can prove the following lemma.

\prop{\emph{P-Q-R} Lemma}{Let ${f\!:A\vdash B}$ be an arrow of
\QDS, let $X_i$ for ${i\in\{1,2,3\}}$ be occurrences of the
predicate letters $P$, $Q$ and $R$, respectively, in $A$, and let
$Y_i$ be occurrences of $P$, $Q$ and $R$, respectively, in $B$,
such that $X_i$ and $Y_i$ are tied in $f$. Let, moreover,
${X_2\mX^2_q\vee X_3\mX^3_r}$ be a subformula of $A$ and
${Y_1\mY^1_p\kon Y_2\mY^2_q}$ a subformula of $B$. Then there is a
$d_{P\mZ^1_p,Q\mZ^2_q,R\mZ^3_r}$-term ${h\!:A'\vdash B'}$ such
that $X_i'$ are occurrences of $P$, $Q$ and $R$, respectively, in
the source ${P\mZ^1_p\kon(Q\mZ^2_q\vee R\mZ^3_r)}$ of the head of
$h$ and $Y_i'$ are occurrences of $P$, $Q$ and $R$, respectively,
in the target ${(P\mZ^1_p\kon Q\mZ^2_q)\vee R\mZ^3_r}$ of the head
of $h$, such that for some arrows ${f_x\!:A\vdash A'}$ and
${f_y\!:B'\vdash B}$ of \QDS\ we have ${f=f_y\cirk h\cirk f_x}$ in
\QDS, and $X_i$ is tied to $X_i'$ in $f_x$, while $Y_i'$ is tied
to $Y_i$ in $f_y$.}

This lemma, exactly analogous to the \mbox{\emph{p-q-r}} Lemma of
\cite{DP07} (Section 2.4), is proved like this previous lemma by
relying on the Gentzenization of \GQDS.

\subsection{Some properties of \QPN}
In this section, by relying on the results of the preceding
section, we establish some results concerning the category \QPN\
introduced at the end of Section 2.1, which we will find useful
for calculations later on. For these results we need to introduce
the following.

Let \QDSN\ be the category defined as \QDS\ save that it is
generated not by \Pe, but by ${\Pe\cup\Pe^{\neg}}$ (see the end of
Section 2.1). So the objects of \QDSN\ are the formulae of
${\eL^{\neg P}}$, i.e.\ the objects of \QPN. For $A$ and $B$
formulae of ${\eL^{\neg P}}$, we define when an occurrence of the
predicate letter $P$ in $A$ is tied to an occurrence of $P$ in $B$
in an arrow ${f\!:A\vdash B}$ of \QDSN\ analogously to what we had
at the beginning of the preceding section.

We say that a finite (possibly empty) sequence $S$ of quantifier
prefixes is \emph{foreign} to a formula $B$ when the set of
variables occurring in $S$ is disjoint from the set of free
variables of $B$. If $S$ is foreign to $B$, then there is an
isomorphism ${j^\str_{S,B}\!:B\vdash SB}$ of \QDSN; defined in
terms of $\gamma^{\forall}$ and $\iota^{\exists}$. The inverse of
$j^\str_{S,B}$ is the arrow ${j^\rts_{S,B}\!:SB\vdash B}$ of
\QDSN\ defined in terms of $\iota^{\forall}$ and
$\gamma^{\exists}$. By \QDS\ Coherence, these isomorphisms are
unique.

We introduce next a generalization of the arrow terms
$\Xi^Q_{B,A}$ obtained by letting the arrow terms $\Xi^Q_{B,A}$
``absorb'' various \QDS\ arrow terms. These generalized terms have
the right form for the first two technical results below---the
$\Xi^Q$-Permutation Lemmata.

For $S$ and $S^\neg$ two independent finite sequences of
quantifier prefixes both foreign to $P\mX_n$, and for $I$ being
the sequence of indices ${P\mX_n,A,\mY_m,S,S^\neg}$, we have
\begin{tabbing}
$\Delta^{\forall}_I\,$\=$=_{df}(\mj_A\kon(\forall_{\mY_m}\!((j^\str_{S^\neg,\neg
P\mX_n}\!\!\!\vee
j^\str_{S,P\mX_n})\cirk\iota^{\forall_{\mX_n}}_{\neg P\mX_n\vee
P\mX_n})\cirk\gamma^{\forall_{\mY_m}}_{\forall_{\mX_n}\!(\neg
P\mX_n\!\vee
P\mX_n)}))\!\cirk\!\Delta^{\forall}_{P\mX_n,A}\!:$\\*[1ex]
\`$A\vdash A\kon\forall_{\mY_m}\!(S^\neg\neg P\mX_n\vee
SP\mX_n)$,\\[2ex]
$\Sigma^{\exists}_I\,$\>$=_{df}\Sigma^{\exists}_{P\mX_n,A}\cirk
((\gamma^{\exists_{\mY_m}}_{\exists_{\mX_n}\!(P\mX_n\!\kon \neg
P\mX_n)}\cirk\exists_{\mY_m}\!(\iota^{\exists_{\mX_n}}_{P\mX_n\kon\neg
P\mX_n}\!\!\cirk(j^\rts_{S,P\mX_n}\!\!\kon j^\rts_{S^\neg,\neg
P\mX_n})))\vee\mj_A)\!:$\\*[1ex] \`$\exists_{\mY_m}\!(SP\mX_n\kon
S^\neg\neg P\mX_n)\vee A\vdash A$.
\end{tabbing}
The analogous abbreviations
\begin{tabbing}
\hspace{9em}\=$\Sigma^{\forall}_I\!:$\=$\;\;A\vdash
\forall_{\mY_m}\!(S^\neg\neg P\mX_n\vee SP\mX_n)\kon A$,\\*[1ex]
\>$\Delta^{\exists}_I\!:$\>$\;\;A\vee\exists_{\mY_m}\!(SP\mX_n\kon
S^\neg\neg P\mX_n)\vdash A$,\\[1.5ex]
\>$\Delta'^{\forall}_I\!\!:$\>$\;\;A\vdash
A\kon\forall_{\mY_m}\!(SP\mX_n\vee S^\neg\neg P\mX_n)$,\\*[1ex]
\>$\Sigma'^{\exists}_I\!:$\>$\;\;\exists_{\mY_m}\!(S^\neg\neg
P\mX_n\kon SP\mX_n)\vee A\vdash A$,\\[1.5ex]
\>$\Sigma'^{\forall}_I\!:$\>$\;\;A\vdash
\forall_{\mY_m}\!(SP\mX_n\vee S^\neg\neg P\mX_n)\kon A$,\\*[1ex]
\>$\Delta'^{\exists}_I\!:$\>$\;\;A\vee\exists_{\mY_m}\!(S^\neg\neg
P\mX_n\kon SP\mX_n)\vdash A$
\end{tabbing}
are defined in terms of $\Delta^{\forall}$ and $\Sigma^{\exists}$
like the analogous abbreviations of Section 2.1. The right
conjunct $\forall_{\mY_m}\!(S^\neg\neg P\mX_n\vee SP\mX_n)$ in the
target of $\Delta^{\forall}_I$ is the \emph{crown} of
$\Delta^{\forall}_I$, and analogously for $\Sigma^{\exists}_I$ and
the other abbreviations, replacing ``right'' by ``left'',
``conjunct'' by ``disjunct'', and ``target'' by ``source'', as
appropriate. Note that here $\mY_m$ is an arbitrary sequence of
variables, and not necessarily the free-variable sequence of
$P\mX_n$ as in $\Xi^Q_{P\mX_n,A}$ (see Section 2.1).

The definition of $\Delta^{\forall}_I$ above is of the form
$f\cirk\Delta^{\forall}_{P\mX_n,A}$. By \QDS\ Coherence, instead
of the arrow term $f$ of \QDSN\ we could have used for this
definition any other arrow term $g$ of \QDSN\ of the same type as
$f$ such that ${Gf=Gg}$, since we have ${g=f}$, and analogously
for $\Sigma^{\exists}_I$, etc.

Let ${\Xi,\Theta\in\{\Delta,\Delta',\Sigma,\Sigma'\}}$, and let a
$\Xi^Q_I$-term be defined as a $\beta$-term in Section 2.2 save
that $\beta$ is replaced by $\Xi^Q_I$, and the clause ``if $f$ is
a $\Xi^Q_I$-term, then $[f]^x_y$ is a $\Xi^Q_I$-term'' is omitted.
Then we have the following analogue of the $\hat{\Xi}$-Permutation
Lemma of \cite{DP07} (Section 2.5).

\prop{$\Xi^\forall$-Permutation Lemma}{Let ${g\!:C\vdash D}$ be a
$\Xi^\forall_{P\mX_n,A,\mY_m,S,S^\neg}$-term of \QPN\ such that
$X_1$ and $\neg X_2$ are respectively the occurrences within $D$
of the predicate letters $P$ and $\neg P$ in the crown of the head
${\Xi^\forall_{P\mX_n,A,\mY_m,S,S^\neg}}$ of $g$, and let
${f\!:D\vdash E}$ be an arrow term of \QDSN\ such that we have an
occurrence $Y_1$ of $P$ and an occurrence $\neg Y_2$ of $\neg P$
within a subformula of $E$ of the form
${{\forall_{\mY'_{m'}}\!(S'Y_1\mX'_n\!\vee {S^\neg}'\neg
Y_2\mX'_n\!)}}$ or ${{\forall_{\mY'_{m'}}\!({S^\neg}'\neg
Y_2\mX'_n\!\vee S'Y_1\mX'_n\!)}}$, for $S'$ and ${S^\neg}'$ finite
sequences of quantifier prefixes, and $X_i$ is tied to $Y_i$ for
${i\in\{1,2\}}$ in $f$. Then there is a
${\Theta^\forall_{P\mX'_n,A',\mY'_{m'},S',{S^\neg}'}}$-term
${g'\!:D'\vdash E}$ of \QPN\ the crown of whose head is
${{\forall_{\mY'_{m'}}\!(S'Y_1\mX'_n\!\vee {S^\neg}'\neg
Y_2\mX'_n\!)}}$ or ${{\forall_{\mY'_{m'}}\!({S^\neg}'\neg
Y_2\mX'_n\!\vee S'Y_1\mX'_n\!)}}$, and there is an arrow term
${f'\!:C\vdash D'}$ of \QDSN\ such that in \QPN\ we have ${f\cirk
g=g'\cirk f'}$.}

\dkz We proceed in principle as for the proof of the
$\hat{\Xi}$-Permutation Lemma in \cite{DP07}, with some
adjustments and additions. We appeal to the Refined Development
Lemma for \QPNN\ (see Section 2.2), and we use the $\kon\vee$
Lemma of the preceding section to ascertain that cases involving
``problematic'' ${d_{A,SP\mX_n,S^\neg \neg P\mX_n}}\!$-terms or
${d_{A,S^\neg \neg P\mX_n,SP\mX_n}}\!$-terms in the developed
arrow term $f$ are excluded.

We rely then on equations analogous to the equations mentioned in
the proof of the $\hat{\Xi}$-Permutation Lemma, where
$\hat{\Xi}_{p,A}$ is replaced by $\Xi^\forall_I$, which entails
further adjustments. Such equations, which are either
stem-increasing, or related to the stem-increasing equations, or
are simply consequences of definitions, are established with the
help of the implications \mbox{($Q\iota$ \emph{canc})} together
with the equations \mbox{($\Xi\tau$)} (see Section 2.1) and \QDS\
Coherence. We rely also on the remark we made before the lemma
concerning the alternative definitions of
$\Xi^\forall_I$.\hfill\mbox{\hspace{2em}}\qed

\vspace{2ex}

We have a dual lemma, called the $\Xi^\exists$-\emph{Permutation
Lemma}, analogous to the $\check{\Xi}$-Permutation Lemma of
\cite{DP07} (Section 2.5), which involves $\Xi^\exists_I$-terms
instead of $\Xi^\forall_I$-terms.

Next we have a lemma analogous to the \mbox{\emph{p-$\neg$p-p}}
Lemma of \cite{DP07} (Section 2.5).

\vspace{1ex}

\prop{\mbox{\emph{P-$\neg$P-P}} Lemma}{Let $X_1$, $\neg X_2$ and
$X_3$ be occurrences of the predicate letters $P$, $\neg P$ and
$P$, respectively, in a formula $A$ of ${\eL^{\neg p}}$, and let
$Y_1$, $\neg Y_2$ and $Y_3$ be occurrences of $P$, $\neg P$ and
$P$, respectively, in a formula $B$ of ${\eL^{\neg p}}$. Let
${g_1\!:A'\vdash A}$ be a
$\Xi^\forall_{P\mX_n,A,\mY_m,S,S^\neg}$-term of \QPN\ such that
$\forall_{\mY_m}\!(S^\neg\neg X_2\mX_n\!\vee SX_3\mX_n\!)$ or
$\forall_{\mY_m}\!(SX_3\mX_n\!\vee S^\neg\neg X_2\mX_n\!)$ is the
crown of the head of $g_1$, let ${g_2\!:B\vdash B'}$ be a
$\Theta^\exists_{P\mX'_n,A',\mY'_{m'},S',{S^\neg}'}$-term of \QPN\
such that $\exists_{\mY'_{m'}}\!(S'Y_1\mX'_n\!\kon {S^\neg}'\neg
Y_2\mX'_n\!)$ or $\exists_{\mY'_{m'}}\!({S^\neg}'\neg
Y_2\mX'_n\!\kon S'Y_1\mX'_n\!)$ is the crown of the head of $g_2$,
and let ${f\!:A\vdash B}$ be an arrow term of \QDSN\ such that
$X_i$ and $Y_i$ are tied in $f$ for ${i\in\{1,2,3\}}$. Then
${g_2\cirk f\cirk g_1}$ is equal in \QPN\ to an arrow term of
\QDSN.}

\noindent The proof of this lemma is analogous to the proof of the
\mbox{\emph{p-$\neg$p-p}} Lemma in \cite{DP07}. We use the
\mbox{\emph{P-Q-R}} Lemma of the preceding section and the
$\Xi^Q$-Permutation Lemmata instead of the \mbox{\emph{p-q-r}}
Lemma and the $\stackrel{\raisebox{-2pt}{\mbox{\tiny $\xi$}}}
{\raisebox{0pt}{$\Xi$}}$-Permutation Lemmata, and we apply the
equation \mbox{($\Sd\Dk$)} of Section 2.1.

We establish in the same manner the
\mbox{\emph{$\neg$P-P-$\neg$P}} \emph{Lemma}, analogous to the
\mbox{\emph{$\neg$p-p-$\neg$p}} Lemma of \cite{DP07} (Section
2.5). The formulation of the \mbox{\emph{$\neg$P-P-$\neg$P}} Lemma
is obtained from that of the \mbox{\emph{P-$\neg$P-P}} Lemma by
replacing the sequence $P,\neg P,P$ by the sequence $\neg P,P,\neg
P$, which entails that $S$ and $S^\neg$, as well as $S'$ and
${S^\neg}'$, are permuted. The \mbox{\emph{$\neg$P-P-$\neg$P}}
Lemma is proved by applying the \mbox{\emph{P-Q-R}} Lemma, the
$\Xi^Q$-Permutation Lemmata and the equation \mbox{($\Sdp\Dkp$)}
of Section 2.1.

\subsection{\QPN\ Coherence}
In \cite{DP07} (Section 2.3) one can find a detailed definition of
a category called \emph{Br}, whose objects are finite ordinals,
and whose arrows are graphs sometimes called \emph{Kelly-Mac Lane
graphs} (because of \cite{KML71}). These graphs may also be found
in \cite{B37} (from whose author the name of \emph{Br} is
derived). We define a functor $G$ from \QPNN\ or \QPN\ into
\emph{Br} as we defined in \cite{DP07} (Section 2.3) an
identically named functor from the categories $\PN^\neg$ and \PN\
into \emph{Br}, without paying attention to variables and
quantifier prefixes. This means that ${GQ_xA=GA}$ (so that $GA$ is
the number of occurrences of predicate letters in the formula
$A$), the arrow $G\alpha$ for $\alpha$ being $\iota^{Q\!_x}_A$,
$\gamma^{Q\!_x}_A$, $\check{\theta}^{\forall\!_x\str}_{A,D}$ and
$\hat{\theta}^{\exists_x\rts}_{A,D}$ is an identity arrow,
${GQ_xf=G[f]^x_y=Gf}$, while $G\Delta^{\forall}_{B,A}$ and
$G\Sigma^{\exists}_{B,A}$ are like $G\Dk_{B,A}$ and $G\Sd_{B,A}$
respectively. The category \emph{Rel} mentioned in Section 1.5 is
a subcategory of \emph{Br}, and $G$ restricted to the \QDS\ part
of \QPNN\ and \QPN\ coincides with the functor $G$ from \QDS\ to
\emph{Rel}.

The theorems that the functors $G$ from \QPNN\ or \QPN\ into
\emph{Br} are faithful functors are called \QPNN\ \emph{Coherence}
and \QPN\ \emph{Coherence} respectively. We establish first \QPN\
Coherence, and \QPNN\ Coherence will be derived from it in Section
2.7.

We prove \QPN\ Coherence by proceeding as for the proof of \PN\
Coherence in \cite{DP07} (Section 2.7), through lemmata analogous
to the Confrontation and Purification Lemmata. Roughly speaking,
the analogue of the Confrontation Lemma says that a
$\Delta^{\forall}_{P\mX_n,A}$-term, called a
$\Delta^{\forall}$-\emph{factor}, and a
$\Sigma^{\exists}_{P\mY_n,B}$-term, called a
$\Sigma^{\exists}$-\emph{factor}, mutually tied in a direct manner
through the crowns, which are called \emph{confronted} factors,
can be permuted with the help of stem-increasing and naturality
equations so that they are ready to get eliminated by applying the
\mbox{\emph{P-$\neg$P-P}} and \mbox{\emph{$\neg$P-P-$\neg$P}}
Lemmata of Section 2.4. The analogue of the Purification Lemma
states that this elimination can be pursued until we obtain an
arrow term without confronted factors, such an arrow term being
called \emph{pure}.

For the proof of these analogues of the Confrontation and
Purification Lemmata we need the Refined Development Lemma for
\QPN\ of Section 2.2. We also need the stem-increasing equations
for $\Delta^{\forall}$ and $\Sigma^{\exists}$ (see Section 2.1)
and the naturality of $\Delta^{\forall}$ and $\Sigma^{\exists}$ in
the stem index. Where in the proof of the Purification Lemma in
\cite{DP07} (Section 2.7) we appealed to Lemma~3, we now appeal to
the $\vee\kon$ Lemma of Section 2.3. Instead of the
\mbox{\emph{p-$\neg$p-p}} and \mbox{\emph{$\neg$p-p-$\neg$p}}
Lemmata we now have the \mbox{\emph{P-$\neg$P-P}} and
\mbox{\emph{$\neg$P-P-$\neg$P}} Lemmata.

The equation \mbox{($\Delta\tau$)} of Section 2.1 is essential,
together with the stem-increas\-ing equations and the naturality
of $\Delta^{\forall}$ and $\Sigma^{\exists}$ in the stem index, to
guarantee that if there is a $\Delta^{\forall}$-factor in a pure
arrow term $f$, then
\[
f=f'\cirk\Delta^{\forall}_{P\mX_n,A}
\]
for any sequence of variables $\mX_n$. The equation
\mbox{($\Sigma\tau$)} of Section 2.1 is needed to establish an
analogous equation for $\Sigma^{\exists}$-factors. By pushing in
this manner to the extreme right the $\Delta^{\forall}$-factors
remaining in a pure arrow term, and to the extreme left the
remaining $\Sigma^{\exists}$-factors, and by relying on \QDS\
Coherence, we establish \QPN\ Coherence.

\subsection{The equivalence of \QPNN\ and \QPN}
To prove that \QPNN\ and \QPN\ are equivalent categories we
proceed as in \cite{DP07} (Section 2.6), with the following
adjustments and additions.

When we define the functor $F$ from \QPNN\ to \QPN\ we have the
following new clauses on objects:
\begin{tabbing}
\mbox{\hspace{7em}}\= $FA=A$,\quad for $A$ of the form $P\mX_n$ or
$\neg P\mX_n$,
\\[1ex]
\>$FQ_xA=Q_xFA$,
\\[1ex]
\>$F\neg\forall_xA=\exists_xF\neg A$,
\\[1ex]
\>$F\neg\exists_xA=\forall_xF\neg A$.
\end{tabbing}
On arrows we have first new clauses analogous to the old clauses
where $\alpha$ is $\iota^{Q\!_x}$, $\gamma^{Q\!_x}$,
$\check{\theta}^{\forall\!_x\str}$ and
$\hat{\theta}^{\exists_x\rts}$, while $\Dk$ and $\Sd$ are replaced
by $\Delta^{\forall}$ and $\Sigma^{\exists}$, the letter $p$ is
replaced by $P\mX_n$, and some further adjustments are made. We
have moreover the following new clauses:
\begin{tabbing}
\hspace{1em}if $x$ is free in $B$,\\*[1ex]
\hspace{2em}\=$F\Delta^{\forall}_{\forall\!_xB,A}\,$\=$=(\mj_A\kon(\forall_{\mX_{n\mn
1}}\!(\check{c}_{\exists_xF\neg
B,\forall\!_xFB}\cirk\check{\theta}^{\forall\!_x\str}_{FB,\exists_xF\neg
B}\cirk\forall_x(\check{c}_{FB,\exists_xF\neg B}\cirk$\\*[1ex]
\`$(\iota^{\exists_x}_{F\neg B}\vee\mj_{FB})))\cirk h))\cirk
F\Delta^{\forall}_{B,A}$,
\end{tabbing}
where
\[
h\!:\forall_{x_{n\mn
1}}\ldots\forall_{x_{i+1}}\forall_x\forall_{x_i}\ldots\forall_{x_1}(F\neg
B\vee FB)\vdash\forall_{\mX_{n\mn 1}}\forall_x(F\neg B\vee FB)
\]
is an isomorphism of \QDSN\ (see the preceding section)
generalizing isomorphisms of the type
$\forall_x\forall_yC\vdash\forall_y\forall_xC$,
\begin{tabbing}
\hspace{2em}\=$F\Delta^{\forall}_{\forall\!_xB,A}\,$\=$=$\kill

\hspace{1em}if $x$ is not free in $B$,\\*[1ex]
\>$F\Delta^{\forall}_{\forall\!_xB,A}\,$\>$=(\mj_A\kon\forall_{\mX_n}\!(\iota^{\exists_x}_{F\neg
B}\vee\gamma^{\forall\!_x}_{FB}))\cirk
F\Delta^{\forall}_{B,A}$,\\[2ex]
\hspace{1em}if $x$ is free in $B$, for $h$ as above,\\*[1ex]
\>$F\Delta^{\forall}_{\exists_xB,A}\,$\>$=(\mj_A\kon(\forall_{\mX_{n\mn
1}}\!(\check{\theta}^{\forall\!_x\str}_{F\neg
B,\exists_xFB}\cirk\forall_x(\mj_{F\neg
B}\vee\iota^{\exists_x}_{FB}))\cirk h))\cirk
F\Delta^{\forall}_{B,A}$,\\[2ex]
\hspace{1em}if $x$ is not free in $B$,\\*[1ex]
\>$F\Delta^{\forall}_{\exists_xB,A}\,$\>$=(\mj_A\kon\forall_{\mX_n}\!(\gamma^{\forall\!_x}_{F\neg
B}\vee\iota^{\exists_x}_{FB}))\cirk F\Delta^{\forall}_{B,A}$,
\end{tabbing}
and dual clauses for $F\Sigma^{\exists}_{\forall\!_xB,A}$ and
$F\Sigma^{\exists}_{\exists_xB,A}$,
\begin{tabbing}
\hspace{2em}\=$FQ_xf\,$\=$=Q_xFf$,\\[1.5ex]
\>$F[f]^x_y$\>$=[Ff]^x_y$.
\end{tabbing}
This defines the functor $F$.

For $f$ an arrow term of \QPNN\ we have that $GFf$ coincides with
$Gf$, where $G$ in $GFf$ is the functor $G$ from \QPN\ to
\emph{Br}, and $G$ in $Gf$ is the functor $G$ from \QPNN\ to
\emph{Br} (see the beginning of the preceding section). To show
that, it is essential to check that $GF\Delta^{\forall}_{B,A}$ and
$GF\Sigma^{\exists}_{B,A}$ coincide with $G\Delta^{\forall}_{B,A}$
and $G\Sigma^{\exists}_{B,A}$ respectively, which is done by
induction on the complexity of the crown index $B$.

Then we can easily verify that $F$, as defined above, is indeed a
functor. If ${f=g}$ in \QPNN, then ${Gf=Gg}$, and hence, as we
have just seen, ${GFf=GFg}$. By \QPN\ Coherence of the preceding
section, we conclude that ${Ff=Fg}$ in \QPN. (To verify that the
functor $F$ from $\PN^\neg$ to \PN\ in Section 2.6 of \cite{DP07}
is a functor we could have proceeded analogously, by establishing
\PN\ Coherence first, before introducing the functor $F$. We did
not need the functor $F$ to prove \PN\ Coherence. This would make
the exposition in \cite{DP07} somewhat simpler, and better
organized.)

We define a functor $F^\neg$ from \QPN\ to \QPNN\ by stipulating
that ${F^\neg A=A}$ and ${F^\neg f=f}$. To show that \QPNN\ and
\QPN\ are equivalent categories via the functors $F$ and $F^\neg$
we proceed as in \cite{DP07} (Section 2.6) with the following
additions. We have the following auxiliary definitions in \QPNN,
for $\mX_n$ being the free-variable sequence of $A$ (see Section
2.1), and $\mY_m$ being this sequence with $x$ omitted (if $x$ is
free in $A$, then ${m=n\mn 1}$; otherwise, $\mX_n$ is $\mY_m$ and
${m=n}$):
\begin{tabbing}
\hspace{1.5em}\=$q^{\forall\!_x\str}_A\,$\=$=_{df}\Sigma'^{\exists}_{\forall\!_xA,\exists_x\neg
A}\cirk(\iota^{\exists_{\mY_m}}_{\neg\forall\!_xA\kon\forall\!_xA}\vee\mj_{\exists_x\neg
A})\cirk d_{\neg\forall\!_xA,\forall\!_xA,\exists_x\neg A}\cirk
$\\*[1ex]
\`$(\mj_{\neg\forall\!_xA}\kon(\check{\theta}^{\forall\!_x\str}_{A,\exists_x\neg
A}\cirk\forall_x((\mj_A\vee\iota^{\exists_x}_{\neg
A})\cirk\iota^{\forall_{\mX_n}}_{A\vee\neg
A})\cirk\gamma^{\forall\!_x}_{\forall_{\mX_n}\!(A\vee\neg
A)}))\cirk\Delta'^{\forall}_{A,\neg\forall\!_xA}\!\!:$\\*[1ex]
\`$\neg\forall_xA\vdash\exists_x\neg A$,\\[2ex]
\>$q^{\forall\!_x\rts}_A$\>$=_{df}\Sigma'^{\exists}_{A,\neg\forall\!_xA}\cirk
((\gamma^{\exists_x}_{\exists_{\mX_n}\!(\neg A\kon
A)}\cirk\exists_x(\iota^{\exists_{\mX_n}}_{\neg A\kon
A}\cirk(\mj_{\neg
A}\kon\iota^{\forall\!_x}_A))\cirk\hat{\theta}^{\exists_x\rts}_{\neg
A,\forall\!_xA})\vee\mj_{\neg\forall\!_xA})$
\\*[1ex]
\`$\cirk d_{\exists_x\neg
A,\forall\!_xA,\neg\forall\!_xA}\cirk(\mj_{\exists_x\neg
A}\kon\iota^{\forall_{\mY_m}}_{\forall\!_xA\vee\neg\forall\!_xA})\cirk
\Delta'^{\forall}_{\forall\!_xA,\exists_x\neg A}\!\!:$\\*[1ex]
\`$\exists_x\neg A\vdash\neg\forall_xA$,
\end{tabbing}
and we have analogous definitions of
\begin{tabbing}
\hspace{12em}\=$q^{\exists_x\str}_A$\=$:\neg\exists_xA\,\,$\=$\vdash\forall_x\neg
A$,\\[1.5ex]
\>$q^{\exists_x\rts}_A$\=$:\forall_x\neg
A$\>$\vdash\neg\exists_xA$.
\end{tabbing}
It can be shown that $q^{Q\!_x\str}_A$ is an isomorphism, with
inverse $q^{Q\!_x\rts}_A$.

Next in the inductive definitions of the isomorphisms
${i_A\!:A\vdash FA}$ and ${i}^{-1}_A\!:FA\vdash A$ we have the
following clauses in addition to clauses in \cite{DP07} (Section
2.6):
\begin{tabbing}
\hspace{12em}$i_A={i}^{-1}_A=\mj_A$, \hspace{.5em} if $A$ is
$P_{\mX_n}$ or $\neg P_{\mX_n}$,
\\[1.5ex]
\hspace{5em}$i_{Q\!_xA}=Q_xi_A$,\hspace{9em}${i}^{-1}_{Q\!_xA}=Q_x{i}^{-1}_A$,
\\[1ex]
\hspace{5em}$i_{\neg\forall\!_xA}\,$\=$=\exists_xi_{\neg A}\cirk
q^{\forall\!_x\str}_A$,\hspace{5.6em}\=
${i}^{-1}_{\neg\forall\!_xA}\,$\=$=q^{\forall\!_x\rts}_A\cirk\exists_x{i}^{-1}_{\neg
A}$,
\\[1ex]
\hspace{5em}$i_{\neg\exists_xA}$\>$=\forall_xi_{\neg A}\cirk
q^{\exists_x\str}_A$,\>${i}^{-1}_{\neg\exists_xA}$\>$=q^{\exists_x\rts}_A\cirk\forall_x{i}^{-1}_{\neg
A}$.
\end{tabbing}

We can then extend the proof of the Auxiliary Lemma of Section 2.6
of \cite{DP07} in order to establish that for ${f\!:A\vdash B}$ we
have in \QPNN\ the equation ${f={i}^{-1}_B\cirk F\!f\cirk i_A}$.
In this extended proof, for the isomorphism
${n^{\rts}_B\!:B\vdash\neg\neg B}$ of \QPNN, we need the following
equation of \QPNN:
\begin{tabbing}
\mbox{\hspace{2em}}\= \mbox{($\Delta^{\forall} n$)}\hspace{4em}\=
$\Delta^{\forall}_{\neg
B,A}\;=(\mj_A\kon\forall_{\mX_n}\!(n^{\rts}_B\vee\mj_{\neg
B}))\cirk\Delta'^{\forall}_{B,A}$,
\end{tabbing}
analogous to the equation \mbox{($\Dk\,n$)} of \cite{DP07}
(Section 2.6, Proof of the Auxiliary Lemma). To derive
\mbox{($\Delta^{\forall} n$)} we use, analogously to what we had
before for the derivation of \mbox{($\Dk\,n$)}, the
stem-increasing equation \mbox{($\forall\Delta^{\forall}$)} of
Section 2.1, the naturality of $\Delta^{\forall}$ in the stem
index, the \mbox{\emph{$\neg$P-P-$\neg$P}} Lemma of the preceding
section and \QDS\ Coherence. (In the derivation of
\mbox{($\Dk\,n$)} in the printed text of \cite{DP07}, Section 2.6,
Proof of the Auxiliary Lemma, ``(with $p$ replaced by $A$)'' is a
misprint for ``(with $p$ replaced by $B$)''.) We derive similarly
an equation analogous to the equation \mbox{($\Dk\,r$)} of
\cite{DP07} (Section 2.6, Proof of the Auxiliary Lemma) involving
$\Delta^{\forall}$.

We need also the following equations of \QPNN, analogous to the
clauses defining $F\Delta^{\forall}_{\forall\!_xB,A}$ above:
\begin{tabbing}
\hspace{1em}if $x$ is free in $B$,\\*[1ex]
\hspace{2em}\=$\Delta^{\forall}_{\forall\!_xB,A}\,$\=$=(\mj_A\kon(\forall_{\mX_{n\mn
1}}\!((q^{\forall\!_x\rts}_B\vee\mj_{\forall\!_xB})\cirk\check{c}_{\exists_x\neg
B,\forall\!_xB}\cirk\check{\theta}^{\forall\!_x\str}_{B,\exists_x\neg
B}\cirk\forall_x(\check{c}_{B,\exists_x\neg B}\cirk$\\*[1ex]
\`$(\iota^{\exists_x}_{\neg B}\vee\mj_{B})))\cirk h))\cirk
\Delta^{\forall}_{B,A}$,
\end{tabbing}
where
\[
h\!:\forall_{x_{n\mn
1}}\ldots\forall_{x_{i+1}}\forall_x\forall_{x_i}\ldots\forall_{x_1}(\neg
B\vee B)\vdash\forall_{\mX_{n\mn 1}}\forall_x(\neg B\vee B)
\]
is an isomorphism of \QDSN,
\begin{tabbing}
\hspace{2em}\=$\Delta^{\forall}_{\forall\!_xB,A}\,$\=$=$\kill

\hspace{1em}if $x$ is not free in $B$,\\*[1ex]
\>$\Delta^{\forall}_{\forall\!_xB,A}\,$\>$=(\mj_A\kon\forall_{\mX_n}\!
((q^{\forall\!_x\rts}_B\vee\mj_{\forall\!_xB})\cirk(\iota^{\exists_x}_{\neg
B}\vee\gamma^{\forall\!_x}_{B})))\cirk \Delta^{\forall}_{B,A}$.
\end{tabbing}
The idea for the derivation of these equations is the same as the
idea for the derivation of \mbox{($\Delta^{\forall} n$)} above. We
need also equations analogous to the clauses defining
$F\Delta^{\forall}_{\exists_xB,A}$ above.

To show that in \QPNN\
\[
[f]^x_y={i}^{-1}_{B^x_y}\cirk F[f]^x_y\cirk i_{A^x_y}
\]
we need the equation ${[i_A]^x_y=i_{A^x_y}}$ of \QPNN, which is
established by induction on the complexity of $A$. For this
induction we use the equation
\[
[\Dk_{B,A}]^x_y=\Dk_{B^x_y,A^x_y}
\]
and analogous equations of \QPNN. The last displayed equation is
established with the help of the equations \mbox{({\it ren}
$\Delta^{\forall}$)} and \mbox{($\Delta\tau$)} of Section 2.1. We
need also the equation
\[
[q^{\forall\!_z\str}_A]^x_y=q^{\forall\!_z\str}_{A^x_y}
\]
of \QPNN, for which we use the equations \mbox{($\Xi\tau$)} (see
Section 2.1) and \QDS\ Coherence. This suffices to establish that
the categories \QPNN\ and \QPN\ are equivalent.

\subsection{\QPNN\ Coherence}
As we said at the beginning of Section 2.5, \QPNN\
\emph{Coherence} is the theorem that the functor $G$ from \QPNN\
to the category \emph{Br} is faithful. We can then prove \QPNN\
Coherence as follows.

\vspace{2ex}

\noindent {\sc Proof of \QPNN\ Coherence.} Suppose that for $f$
and $g$ arrows of \QPNN\ of the same type we have ${Gf=Gg}$. Then,
as we noted after the definition of the functor $F$ from \QPNN\ to
\QPN\ in the preceding section, we have ${GFf=GFg}$, and hence
${Ff=Fg}$ in \QPN\ by \QPN\ Coherence of Section 2.5. It follows
that ${f=g}$ in \QPNN\ by the equivalence of the categories \QPNN\
and \QPN\ established in the preceding section. \qed

\vspace{2ex}

With \QPNN\ Coherence we can establish easily equations of \QPNN\
whose derivation may otherwise be quite demanding. We have, for
example, the following equations in \QPNN:
\begin{tabbing}
\hspace{1em}\=$\check{\theta}^{\forall\!_x\str}_{A,D}\:$\=$=((\forall_x(\Dd_{D,A}\cirk
d^R_{A,D,\neg D})\cirk\hat{\theta}^{\forall\!_x\rts}_{A\vee D,\neg
D})\vee\mj_D)\cirk d_{\forall\!_x(A\vee D),\neg D,
D}\cirk\Dk_{D,\forall\!_x(A\vee D)}$,\\*[1ex]

\>$\hat{\theta}^{\exists_x\rts}_{A,D}$\>$=\Ddp_{D,\exists_x(A\kon
D)}\cirk d^R_{\exists_x(A\kon D),\neg D,
D}\cirk((\check{\theta}^{\exists_x\str}_{A\kon D,\neg
D}\cirk\exists_x(d_{A,D,\neg D}\cirk\Dkp_{D,A}))\kon\mj_D)$
\end{tabbing}
(see the end Section 1.2 for the definitions of
$\hat{\theta}^{\forall\!_x\rts}_{A\vee D,\neg D}$ and
$\check{\theta}^{\exists_x\str}_{A\kon D,\neg D}$). These
equations say that the distributivity arrow terms
$\check{\theta}^{\forall\!_x\str}_{A,D}$ and
$\hat{\theta}^{\exists_x\rts}_{A,D}$ are definable in \QPNN\ in
terms of the remaining primitive arrow terms and operations on
arrow terms. If these distributivity arrow terms are taken as
defined when we introduce \QPNN, then the equations
${(Q\th{\xi}{}\th{\xi}{})}$ of Section 1.2 become superfluous as
axioms---they can be derived from the remaining axiomatic
equations.

We define a contravariant endofunctor of \QPNN, i.e.\ a functor
from \QPNN\ to $\QPN^{\neg op}$, in the following manner, for
${f\!:A\vdash B}$:
\[
\neg f =_{df}\;\Sdp_{B,\neg A}\!\cirk d_{\neg
B,B,A}\cirk(\mj_{\neg B}\kon(f\vee\mj_{\neg
A}))\cirk\!\Dkp_{A,\neg B}: \neg B\vdash\neg A,
\]
and we verify that this is indeed a contravariant functor by
proceeding as in \cite{DP07} (Section 2.8, where there is also an
alternative definition of $\neg f$). In the course of this
verification, we establish easily with the help of \QPNN\
Coherence that $\Xi^Q$ is a dinatural transformation in the crown
index (see \cite{ML71}, Section IX.4, for the notion of dinatural
transformation).

\section{Coherence of \QMDS\ and \QMPNN}
\subsection{\QMDS\ Coherence}
The category \QMDS\ is defined as the category \QDS\ in Section
1.2 save that we have the additional primitive arrow terms
\[
m_{A,B}\!:A\kon B\vdash A\vee B
\]
for all formulae $A$ and $B$ of \eL, and we assume the following
additional equations:
\begin{tabbing}
\hspace{2em}\=($m$~{\it nat})\quad\=$m_{A\kon B,C}\cirk
\hat{b}^{\str}_{A,B,C}\,$\=$=m_{D,E}\cirk(f\kon g)$\kill

\>($m$~{\it nat})\>$ (f\vee g)\cirk m_{A,B} = m_{D,E}\cirk(f\kon
g)$, \quad for $f\!:A\vdash D$ and $g\!:B\vdash E$,
\\[1.5ex]
\>$(\hat{b}\,m)$\>$m_{A\kon B,C}\cirk \hat{b}^{\str}_{A,B,C}$\>
$=d_{A,B,C}\cirk (\mj_A\kon m_{B,C})$,
\\*[1ex]
\>$(\check{b}\,m)$ \> $\check{b}^{\str}_{C,B,A}\cirk m_{C,B\vee
A}$\>$=(m_{C,B}\vee \mj_A)\cirk d_{C,B,A}$,
\\*[1.5ex]
\>$(c\,m)$\>$m_{B,A}\cirk \hat{c}_{A,B}=\check{c}_{B,A}\cirk
m_{A,B}$.
\end{tabbing}
The proof-theoretical principle underlying $m_{A,B}$ is called
\emph{mix} (see the Gentzen operation below, and \cite{DP04},
Section 8.1, where references are given).

To obtain the functor $G$ from \QMDS\ to the category \emph{Rel}
(see Section 1.5), or to the category \emph{Br} (see Section 2.5),
we extend the definition of the functor $G$ from \QDS\ to
\emph{Rel} by adding the clause that says that ${Gm_{A,B}}$ is an
identity arrow. To prove that this functor $G$ is faithful---this
result is called \QMDS\ \emph{Coherence}---we extend the proof of
\QDS\ Coherence of the first part of this paper.

The Gentzenization of \QMDS\ is obtained with the category \GQMDS,
which has an additional Gentzen operation
\[
\f{f\!:U\vdash Z \quad\quad\quad\quad g\!:Y\vdash W}
{\mix(f,g)=_{dn}(f\vee g)\cirk m_{U,Y}\!:U\kon Y\vdash Z\vee W}
\]
The Cut-Elimination Theorem is proved for \GQMDS\ by enlarging the
proof we had for \GQDS\ in Section 1.9 with an additional case
dealt with in \cite{DP04} (Section 8.4). The preparation for this
Cut-Elimination Theorem involving variable-purity is not impeded
by the presence of $\mix$.

To prove the invertibility lemmata we need for \GQMDS\ we rely on
the following equations of \GQMDS:
\begin{tabbing}
\hspace{2em}\=($m$~{\it nat})\hspace{5em}\=$m_{A\kon B,C}\cirk
\hat{b}^{\str}_{A,B,C}\,$\=$=m_{D,E}\cirk(f\kon g)$\kill

\>\mbox{({\it
mix}~$Q^S$)}\>$\mix(Q^S_{x,X}f_1,f_2)=Q^S_{x,X}\mix(f'_1,f_2)$,
\end{tabbing}
for $f'_1$ being as for \mbox{($\!\ks\!Q^S$)} in Section 1.10 and
${S\in\{L,R\}}$. These equations are either straightforward to
derive, or when ${Q^S\in\{\forall^R,\exists^L\}}$ we derive them
by imitating the derivation of the equation of case ($\forall$2)
of Section 1.9, with the help of the equations
($Q\beta$~\emph{red}) and ($Q\eta$~\emph{red}) (see the end of
Section 1.5). To prove the new Invertibility Lemmata for $\kon$
and $\vee$ we enlarge the proofs of such invertibility lemmata we
had for \GQDS\ in Section 1.10 with cases involving $\mix$ covered
by the remarks preceding the Invertibility Lemma for $\mix$ in
\cite{DP04} (Section 8.4). The proofs of the new Invertibility
Lemmata for $\forall^R$ and $\exists^L$ are taken over unchanged.
To prove the new Invertibility Lemmata for $\forall^L$ and
$\exists^R$ we use in addition the equations \mbox{({\it
mix}~$\forall^L$)} and \mbox{({\it mix}~$\exists^R$)}
respectively.

We need moreover a new Invertibility Lemma for $\mix$, analogous
to the lemma with the same name in \cite{DP04} (Section 8.4). The
proof of this new lemma is based on the proof in \cite{DP04} and
on the equations \mbox{({\it mix}~$Q^S$)}. This suffices to
establish \QMDS\ Coherence.

\subsection{\QMPNN\ Coherence}
We introduce now the category \QMPNN, which corresponds to the
multiplicative fragment without propositional constants of
classical linear first-order predicate logic with mix. The
category \QMPNN\ is defined as the category \QPNN\ in Section 2.1
save that we have the additional primitive arrow terms
${m_{A,B}\!:A\kon B\vdash A\vee B}$ for all formulae $A$ and $B$
of \eLn, and we assume as additional equations ($m$~{\it nat}),
$(\hat{b}\,m)$, $(\check{b}\,m)$ and $(c\,m)$ of the preceding
section. To obtain the functor $G$ from \QMPNN\ to the category
\emph{Br} we extend what we had for the functor $G$ from \QPNN\ to
\emph{Br} (see Section 2.5) with the clause that says that
${Gm_{A,B}}$ is an identity arrow. The theorem asserting that this
functor is faithful is called \QMPNN\ Coherence.

The category \QMPN\ is defined as the category \QPN\ at the end of
Section 2.1 save that we have the additional primitive arrow terms
${m_{A,B}}$ for all objects $A$ and $B$ of \QPN, and we assume the
additional equations ($m$~{\it nat}), $(\hat{b}\,m)$,
$(\check{b}\,m)$ and $(c\,m)$. We can prove that \QMPNN\ and
\QMPN\ are equivalent categories as in Section 2.6, with trivial
additions.

The proof of \QMPNN\ Coherence is then reduced to the proof of
\QMPN\ Coherence, and the latter proof can be obtained quite
analogously to what we have in \cite{DP07} (Sections 6.1-2). Here
are some remarks concerning additions and changes.

The problem here is that the $\kon\vee$ Lemma of Section 2.3,
which was used for proving the $\Xi^\forall$-Permutation Lemma for
\QPN\ in Section 2.4, does not hold for \QMDS\ (the $\vee\kon$
Lemma of Section 2.3 holds for \QMDS). We can nevertheless prove a
modified version of the $\Xi^\forall$-Permutation Lemma, where we
assume that $Y_1$ and $Y_2$ occur within a subformula of $E$ of
the form ${\neg P\mX'_n\!\kon(Y_1\mX'_n\!\vee \neg Y_2\mX'_n\!)}$
or ${P\mX'_n\!\kon(\neg Y_2\mX'_n\!\vee Y_1\mX'_n\!)}$. For the
proof of this modified version of the $\Xi^\forall$-Permutation
Lemma we rely on some auxiliary results, which we will now
consider.

Let us call \emph{quasi-atomic} formulae of \eL\ all formulae of
the form $SP\mX_n$ for $S$ a finite sequence of quantifier
prefixes, i.e.\ formulae in which $\kon$ and $\vee$ do not occur.
For $X$ a particular occurrence of a predicate letter in a formula
$A$ such that there is a subformula of the form $B\ks C$ or $C\ks
B$ of $A$ where $C$ is a quasi-atomic formula in which $X$ occurs,
let $A^{-X}$ be obtained from $A$ by replacing the particular
subformula $B\ks C$ or $C\ks B$ by $B$. When $X$ occurs in $A$ as
we have just said we say that $X$ is \emph{deletable} from $A$.

For ${i\in\{1,2\}}$, let $A_i$ be a formula of \eL, let $X_i$ be
an occurrence of the predicate letter $P$ deletable in $A_i$, and
let $X_1$ and $X_2$ be tied in the arrow ${f\!:A_1\vdash A_2}$ of
\QMDS\ (see the beginning of Section 2.3 for the meaning of
``tied''). The new version of Lemma~1 of Section 6.1 of
\cite{DP07} then says that there is an arrow term
${f^{-P}\!:A^{-X_1}_1\vdash A^{-X_2}_2}$ of \QMDS\ such that
$Gf^{-P}$ is obtained from $Gf$ by deleting the pair corresponding
to ${(X_1,X_2)}$. In the proof of Lemma~1 of Section 6.1 in the
printed version of \cite{DP07} there is an omission. The last
sentence of the first paragraph should be replaced by: ``If $x_i$
is not a proper subformula of the subformula $B_j$, then
$d_{B_1,q,B_3}^{-q}$ is $m_{B_1,B_3}$ or $f^{-q}$ is
$\mj_{A^{-x_i}_i}$.'' The proof of the new version of Lemma~1 is
then analogous to the old proof with the addition in the induction
step that when $f$ is $g\ks h$ or $h\ks g$, then $f^{-P}$ is $g$
not only when $h$ is equal to $\mj_{P\mX_n}$, but also when it is
of a type ${B_1\vdash B_2}$ such that $B_i$ is a quasi-atomic
subformula of $A_i$ in which $X_i$ occurs. Note that this deleting
lemma does not hold for \QDS, because we cannot cover
$d_{B_1,SP\mX_n,B_3}^{-P}$.

A \emph{context} $Z$ is obtained from a formula of \eL\ by
replacing a particular occurrence of an atomic subformula with a
place holder \koc. We write $Z(A)$ for the formula obtained by
putting the formula $A$ at the place of \koc\ in $Z$, and we write
$Z(f)$ for the arrow term obtained by putting the arrow term $f$
at the place of \koc\ in $Z$ and $\mj_B$ at the place of every
atomic formula $B$ in $Z$. For $X$ and $Y$ contexts, let
${f\!:X(P\mX_n\!)\kon \forall_{\mY_m}B\vdash Y(P\mX_n\!\kon B)}$
be an arrow term of \QMDS\ such that $\mY_m$ are all the variables
free in $B$, the displayed occurrences of $P$ in the source and
target are tied in $f$, and the same holds for the $k$-th
occurrence of predicate letter (counting from the left) in the
displayed occurrences of $B$ in the source and target, for every
$k$. Then, by successive applications of the new version of
Lemma~1, we obtain the arrow term ${f^{-B}\!:X(P\mX_n\!)\vdash
Y(P\mX_n\!)}$ of \QMDS\ such that the displayed occurrences of $P$
in the source and target are tied in $f$.

Let ${f^{\dag}\!:X(P\mX_n\!\kon \forall_{\mY_m}B)\vdash
Y(P\mX_n\!\kon \forall_{\mY_m}B)}$ be the arrow term of \QMDS\
obtained from $f^{-B}$ by replacing $P\mX_n$ by ${P\mX_n\!\kon
\forall_{\mY_m}B}$ in the indices of the primitive arrow terms of
$f^{-B}$ at places corresponding to the occurrences displayed in
the source and target (see \cite{DP07}, Section 6.1, for an
example). The new version of Lemma 2$\kon$ of Section 6.1 of
\cite{DP07} should state the following:

\vspace{2ex}

{\it Let $f$ and $f^{\dag}$ be as above. Then there is an arrow
term}
\[
{h_X\!\!:X(P\mX_n\!)\kon \forall_{\mY_m}B\vdash X(P\mX_n\!\kon
\forall_{\mY_m}B)}
\]

{\it of \QDS\ such that
${f=Y(\mj_{P\mX_n}\!\kon\iota^{\forall_{\mY_m}}_B\!)\cirk
f^{\dag}\cirk h_X}$ in \QMDS.}

\vspace{2ex}

\noindent In the proof of this new lemma, when we define
inductively $h_X$, besides clauses analogous to the old clauses,
we should have the additional clauses
\begin{tabbing}
\hspace{11em}\=$h_{\forall\!_xZ}\;\,$\=$=\forall_xh_Z\cirk\hat{\theta}^{\forall\!_x\rts}_{Z(P\mX_n\!),
\forall_{\mY_m}\!B}$,\\[1.5ex]
\>$h_{\exists_xZ}$\>$=\exists_xh_Z\cirk\hat{\theta}^{\exists_x\rts}_{Z(P\mX_n\!),\forall_{\mY_m}\!B}$
\end{tabbing}
(see the end Section 1.2 for the definition of
$\hat{\theta}^{\forall\!_x\rts}_{Z(P\mX_n\!),\forall_{\mY_m}\!B}$).

There is an analogous new version of Lemma 2$\vee$ of Section 6.1
of \cite{DP07}. The proof of \QMPN\ Coherence then proceeds as in
Section 6.2 of \cite{DP07}. This suffices to establish \QMPNN\
Coherence.

\section{Concluding remarks}

In this paper we have not dealt with the multiplicative
propositional constants because, as we said in the Introduction,
they raise problems for coherence understood as the existence of a
faithful functor into the category \emph{Br}. Star-autonomous
categories have however unit objects corresponding to the
multiplicative propositional constants, and it would be
interesting to define a notion of star-autonomous category with
quantifiers, i.e.\ with functors corresponding to quantifiers. It
would be desirable that our category \QPNN\ be isomorphic to a
full subcategory of a category equivalent with the freely
generated star-autonomous category with quantifiers. (We are not
looking for a full subcategory of the freely generated
star-autonomous category with quantifiers, but for a full
subcategory of a category \emph{equivalent} with this category,
because of a difference in language; see \cite{DP07}, Chapter~3).

We have proved an analogous result, which concerns propositional
linear logic, in \cite{DP07} (Chapter~4); we have proved namely
that the category $\PN^{\neg}$ is isomorphic to a full subcategory
of a category equivalent with the freely generated star-autonomous
category. This shows that our notion of category for which
$\PN^{\neg}$ is the freely generated one is the right notion of
star-autonomous category without units.

Besides \cite{DP07}, a systematic work devoted to star-autonomous
categories without units is \cite{H07}. It introduces a
differently defined notion, for which it is supposed that it is
equivalent to ours (see also \cite{HHS}).

The freely generated star-autonomous category with quantifiers
should be obtained by extending the freely generated
star-autonomous category with assumptions involving quantifiers
like those introduced in this paper. We suppose that the approach
of \cite{DP07} could be extended to prove for this category the
result mentioned in the first paragraph of this section.

\vspace{5ex}

\noindent {\small {\it Acknowledgement$\,$}. We are grateful to
Phil Scott and to an anonymous referee for helpful suggestions.
Work on this paper was supported by the Ministry of Science of
Serbia (Grants 144013 and 144029).}

\end{document}